\pdfoutput=1
\RequirePackage{ifpdf}
\ifpdf 
\documentclass[pdftex]{sigma}
\else
\documentclass{sigma}
\fi

\usepackage{mathrsfs}
\usepackage{accents}
\usepackage[all]{xy}

\numberwithin{equation}{section}

\newtheorem{Theorem}{Theorem}[section]
\newtheorem*{Theorem*}{Theorem}
\newtheorem{Corollary}[Theorem]{Corollary}
\newtheorem{Lemma}[Theorem]{Lemma}
\newtheorem{Proposition}[Theorem]{Proposition}
 { \theoremstyle{definition}
\newtheorem{Definition}[Theorem]{Definition}

\newtheorem{Example}[Theorem]{Example}
\newtheorem{Remark}[Theorem]{Remark} }

\newcommand{\C}{\mathbb{C}}
\newcommand{\R}{\mathbb{R}}
\newcommand{\Z}{\mathbb{Z}}
\newcommand{\N}{\mathbb{N}}
\newcommand{\Q}{\mathbb{Q}}

\newcommand{\id}{\mathop{\mathrm{id}}\nolimits}

\newcommand{\SCS}{\mathscr S}

\newcommand{\CO}{\mathcal O}

\renewcommand{\tilde}{\widetilde}
\renewcommand{\setminus}{\smallsetminus}
\newcommand{\nin}{/\kern-2.1ex\in}
\newcommand{\abs}[1]{\lvert#1\rvert}
\newcommand{\norm}[1]{\lVert#1\rVert}
\newcommand{\vardbtilde}[1]{\tilde{\raisebox{0pt}[0.85\height]{$\tilde{#1}$}}}

\def\t{\mathfrak t}
\def\dev{\operatorname{dev}}
\def\re{\operatorname{Re}}
\def\im{\operatorname{Im}}
\def\<{\left\langle}
\def\>{\right\rangle}
\def\End{\operatorname{End}}

\def\Aut{\operatorname{Aut}}
\def\Hom{\operatorname{Hom}}

\def\GL{\operatorname{GL}}
\def\J{\mathrm{J}}

\begin{document}

\allowdisplaybreaks

\newcommand{\arXivNumber}{1904.04076}

\renewcommand{\PaperNumber}{065}

\FirstPageHeading

\ShortArticleName{Adiabatic Limit, Theta Function, and Geometric Quantization}

\ArticleName{Adiabatic Limit, Theta Function,\\ and Geometric Quantization}

\Author{Takahiko YOSHIDA}

\AuthorNameForHeading{T.~Yoshida}

\Address{Department of Mathematics, School of Science and Technology, Meiji University,\\ 1-1-1 Higashimita, Tama-ku, Kawasaki, 214-8571, Japan}
\Email{\href{mailto:takahiko@meiji.ac.jp}{takahiko@meiji.ac.jp}}

\ArticleDates{Received March 20, 2023, in final form July 06, 2024; Published online July 19, 2024}

\Abstract{Let $\pi\colon (M,\omega)\to B$ be a non-singular Lagrangian torus fibration on a~complete base $B$ with prequantum line bundle $\bigl(L,\nabla^L\bigr)\to (M,\omega)$. Compactness on~$M$ is not assumed. For a positive integer $N$ and a compatible almost complex structure~$J$ on $(M,\omega)$ invariant along the fiber of $\pi$, let $D$ be the associated Spin${}^c$ Dirac operator with coefficients in $L^{\otimes N}$. First, in the case where $J$ is integrable, under certain technical condition on $J$, we give a complete orthogonal system $\{ \vartheta_b\}_{b\in B_{\rm BS}}$ of the space of holomorphic $L^2$-sections of $L^{\otimes N}$ indexed by the Bohr--Sommerfeld points~$B_{\rm BS}$ such that each $\vartheta_b$ converges to a delta-function section supported on the corresponding Bohr--Sommerfeld fiber $\pi^{-1}(b)$ by the adiaba\-tic\mbox{(-type)} limit. We also explain the relation of $\vartheta_b$ with Jacobi's theta functions when $(M,\omega)$ is~$T^{2n}$. Second, in the case where $J$ is not integrable, we give an orthogonal family \smash{$\big\{ {\tilde \vartheta}_b\big\}_{b\in B_{\rm BS}}$} of $L^2$-sections of $L^{\otimes N}$ indexed by $B_{\rm BS}$ which has the same property as above, and show that each $D{\tilde \vartheta}_b$ converges to $0$ by the adiabatic(-type) limit with respect to the $L^2$-norm.}


\Classification{53D50; 58H15; 58J05}

\section{Introduction}
The purpose of this paper is to investigate the relationship between Spin${}^c$ quantization and real quantization from the viewpoint of the adiabatic(-type) limit for Lagrangian torus fibrations on complete bases. In this paper, a Lagrangian torus fibration is assumed to be non-singular, but its total space is not assumed to be compact unless otherwise stated.

\subsection{Background and motivation}
First let us explain the motivation which comes from geometric quantization. For geometric quantization, see \cite{Hall, Ki2,Sn,Woodhouse}. In physics, quantization is the procedure for building quantum mechanics starting from classical mechanics. In the mathematical context, it is often thought of as a representation of the Poisson algebra consisting of certain functions on a~symplectic manifold to some Hilbert space, so called the quantum Hilbert space. When a~symplectic manifold~$(M,\omega)$ and a prequantum line bundle $\bigl(L,\nabla^L\bigr)\to (M,\omega)$ are given, geometric quantization provides us with a method to construct the quantum Hilbert space and the representation from these data
geometrically. In the theory of geometric quantization by Kostant and Souriau~\cite{Kos1,Sou2, Sou1}, we need an additional structure which is called a polarization to obtain the quantum Hilbert space. By definition, a polarization is an integrable Lagrangian distribution~$P$ of the complexified tangent bundle $TM\otimes \C$ of $(M,\omega)$. When a polarization~$P$ is given, the quantum Hilbert space is naively given as the space of $L^2$-sections of $\bigl(L,\nabla^L\bigr)$ covariant constant along~$P$.\looseness=1

A common example is the K\"ahler polarization. When $(M,\omega)$ is K\"ahler (i.e., $(M,\omega)$ has a~compatible complex structure) and $\bigl(L,\nabla^L\bigr)$ is a holomorphic line bundle with canonical connection, we can take $T^{0,1}M$ as a polarization, and the obtained quantum Hilbert space is nothing but the space of holomorphic $L^2$-sections. This polarization is called the K\"ahler polarization and the quantization procedure is called the K\"ahler quantization. Note that when $M$ is compact and the Kodaira vanishing theorem holds, the quantum Hilbert space is $H^0(M,\CO_L)$. In particular, its dimension is equal to the index of the Dolbeault operator with coefficients in $L$.

Another example is a real polarization. Suppose $(M,\omega)$ admits a structure of a Lagrangian torus fibration $\pi\colon (M,\omega)\to B$. For each point $b\in B$ of the base manifold $B$, the restriction \smash{$\bigl(L,\nabla^L\bigr)|_{\pi^{-1}(b)}$} of $\bigl(L,\nabla^L\bigr)$ to the fiber $\pi^{-1}(b)$ is a flat line bundle. Let \smash{$H^0 \bigl( \pi^{-1}(b); \bigl(L,\nabla^L\bigr)|_{\pi^{-1}(b)}\bigr)$} be the space of covariant constant sections of \smash{$\bigl(L,\nabla^L\bigr)|_{\pi^{-1}(b)}$}. Then, a point $b\in B$ is said to be Bohr--Sommerfeld if $H^0\bigl(\pi^{-1}(b); \bigl(L,\nabla^L\bigr)|_{\pi^{-1}(b)}\bigr)\not=\{ 0\}$. It is well known that Bohr--Sommerfeld points appear discretely. We denote by $B_{\rm BS}$ the set of Bohr--Sommerfeld points. In this case, we can take $T_\pi M\otimes \C $, the complexified tangent bundle along the fiber of $\pi$ as a polarization, and if $M$ is compact, the quantum Hilbert space is given by $\oplus_{b\in B_{\rm BS}}H^0\bigl(\pi^{-1}(b); \bigl(L,\nabla^L\bigr)|_{\pi^{-1}(b)}\bigr)$. See~\cite{Sn} for more details. In this paper, we call this quantization the real quantization.

When $(M,\omega)$ has a structure of a Lagrangian torus fibration $\pi\colon (M,\omega)\to B$ as well as a~K\"ahler structure,
it is natural to ask whether K\"ahler and real quantizations give the same results. This paper focuses on the quantum Hilbert spaces obtained by both quantizations.
It is easy to see that for any compatible almost complex structure $J$ of $(M,\omega)$, $(TM,J)$ admits a flat connection as a complex vector bundle. So, if $(M,\omega ,J)$ is closed K\"ahler, by \cite[Exemples~12.5.2\,(ii)]{MR0278248}, as a manifold, $M$ is finitely covered by a torus. A typical example of this case is an abelian variety, whose geometric quantization is well understood in the context of the theory of theta functions. For example, see \cite{BMN}. Moreover, even in the non-compact case, with an appropriate choice of the quantum Hilbert space for the real quantization, the relationship between these quantizations has been investigated. For the cotangent bundle of the Lie group of compact type, these are related by the generalized Segal--Bargmann transform \cite{MR1892462}.

A completely integrable system can be thought of as a Lagrangian fibration with singular fibers. As one of such examples, for the moment map of a smooth projective toric variety, Danilov showed in~\cite{Dan} that $H^0(M,\CO_L)$ has the irreducible decomposition \smash{$H^0(M,\CO_L)=\bigoplus_{m \in\Delta \cap \t^*_{\Z}}\C_m$} as a~compact torus representation, where $\Delta$ is the moment polytope, $\t^*_{\Z}$ is the weight lattice, and~$\C_m$ is the irreducible representation of the torus with weight $m$. In this case, since singular fibers of the moment map are isotropic tori, the Bohr--Sommerfeld condition is still meaningful for singular fibers, and $\Delta \cap \t^*_{\Z}$ is identified with the set of Bohr--Sommerfeld points. This implies the dimensions of the quantum Hilbert spaces obtained by the above two kinds of quantizations agree with each other. It has been shown that similar results hold for the Gelfand--Cetlin system on the flag variety~\cite{GuS4}, the Goldman system on the moduli space of flat ${\rm SU}(2)$ connections on a surface~\cite{JW1}, and the Kapovich--Millson system on the polygon space~\cite{Kamiyama}.

Moreover, for smooth projective toric varieties, in their paper~\cite{BFMN}, Baier--Florentino--Mour\~ao--Nunes developed a geometric approach to understand the relationship between K\"ahler and real quantizations. Namely, they gave a one-parameter family of complex structures $\{J_t\}_{t\in[0,\infty)}$ and a~basis $\{s_t^m\}_{m\in \Delta \cap \t^*_{\Z}}$ of the space of holomorphic sections associated with the complex structure~$J_t$ for each $t$ such that each section $s_t^m$ converges to a delta function section supported on the corresponding Bohr--Sommerfeld fiber as $t$ goes to $\infty$. The similar results have been obtained for flag manifolds~\cite{Hamilton-Konno} and smooth irreducible complex algebraic varieties~\cite{Hamilton-Harada-Kaveh}. But in \cite{Hamilton-Konno} and~\cite{Hamilton-Harada-Kaveh} the convergence was showed only for the non-singular Bohr--Sommerfeld fibers whereas in~\cite{BFMN} it was showed for all Bohr--Sommerfeld fibers.

The K\"ahler quantization can be generalized to a non-integrable compatible almost complex structure on a closed $(M,\omega)$. When a compatible almost complex structure $J$ on $(M,\omega)$ is given, we can consider the associated Spin${}^c$ Dirac operator $D$ acting on $\Gamma \bigl(\wedge^\bullet T^*M^{0,1}\otimes L\bigr)$. It is well known that $D$ is a formally self-adjoint, first order, elliptic differential operator of degree-one, and if $J$ is integrable, $D$ agrees with the Dolbeault operator up to constant. If $J$ is not integrable, $T^{0,1}M$ is no more polarization. But, even in this case, since $D$ is Fredholm, we can still take the element of the K-theory of a point
\[
\ker ( D|_{\wedge^{\rm even}T^*M^{0,1}\otimes L}) -\ker ( D|_{\wedge^{\rm odd}T^*M^{0,1}\otimes L}) \in K(pt)
\]
as a (virtual) quantum Hilbert space. Its virtual dimension is equal to the index of $D$. We call this quantization the Spin${}^c$ quantization. It has been showed in~\cite{Andersen, FFY1,Kubota} that the equality between dimensions of two quantum Hilbert spaces still holds by replacing K\"ahler quantization with Spin${}^c$ quantization in terms of the index theory.

\subsection{Main theorems}
In this paper, we apply the approach taken in \cite{BFMN} to both of K\"ahler and Spin${}^c$ quantizations of Lagrangian torus fibrations. Our setting is as follows. Let $\pi\colon (M,\omega)\to B$ be a Lagrangian torus fibration on a complete base $B$ with prequantum line bundle $\bigl(L,\nabla^L\bigr)\to (M,\omega)$. We do not assume $M$ is compact. Let $J$ be
a compatible almost complex structure of $(M,\omega)$ invariant along the fiber of $\pi$ in the sense of Lemma~\ref{invariant J}. For $J$, in Section~\ref{adiabatic limit}, we give a one-parameter family $\big\{ J^t\big\}_{t>0}$ of compatible almost complex structures of $(M,\omega)$ with $J^1=J$ so that the fiber shrinks as $t$ goes to $\infty$ with respect to the associated Riemannian metrics. We also show that $J$ is integrable if and only if every $J^t$ is integrable. For $t>0$ and a positive integer $N$, let $D^t$ be the Spin${}^c$ Dirac operator with coefficients in $L^{\otimes N}$ associated with $J^t$.

Firstly, let us consider the case where $J$ is integrable. In this case, we show the following Theorem, which is obtained by putting Corollary~\ref{cor_main2} and Theorem~\ref{main3} together.
\begin{Theorem}\label{main integrable case}
Under the above setting, assume that $J$ is integrable and satisfies certain technical condition. Then, for each $t>0$, 
there exists a complete orthogonal system \smash{$\big\{\vartheta^t_b\big\}_{b\in B_{\rm BS}}$} of~holomorphic $L^2$-sections of $L^{\otimes N}\to \bigl(M,N\omega ,J^t\bigr)$ indexed by the Bohr--Sommerfeld points $B_{\rm BS}$ such that each $\vartheta^t_b$ converges as a delta-function section supported on $\pi^{-1}(b)$ as $t\to \infty$ in the following sense, for any $L^2$-section $s$ of $L^{\otimes N}$, we have
\begin{equation*}
\lim_{t\to \infty}\int_M\left\langle s, \frac{\vartheta^t_b}{\norm{\vartheta^t_b}_{L^1}} \right\rangle_{L^{\otimes N}} (-1)^{\frac{n(n-1)}{2}}\dfrac{\omega^n}{n!}=\int_{\pi^{-1}(b)}\left\langle s,\delta_b\right\rangle_{L^{\otimes N}} \abs{{\rm d}y},
\end{equation*}
where $\langle \, ,\, \rangle_{L^{\otimes N}}$ is the Hermitian metric of $L^{\otimes N}$, $\delta_b$ is the covariant constant section of \linebreak $\bigl(L, \smash{\nabla^L\bigr)^{\otimes N}|_{\pi^{-1}(b)}}$
defined by \eqref{delta_b}, and $\abs{{\rm d}y}$ is the natural one-density on $\pi^{-1}(b)$.
\end{Theorem}
We give $\big\{\vartheta^t_b\big\}_{b\in B_{\rm BS}}$ explicitly in Section~\ref{the integrable case}.

One of examples of the total space of a Lagrangian torus fibration with complete base is the abelian variety. In this case, we show that, under an appropriate choice of $J$, each $\vartheta_b$ coincides with Jacobi's theta function up to function on the base space (see Theorem~\ref{relation with theta}). For the theta functions, see \cite{Mumford1,Mumford3}.

We remark that there are several works which deal with theta functions from the viewpoint of geometric quantization of Lagrangian fibrations, for example, \cite{BMN,Mumford3,Nohara1, Nohara2}. In~\cite{BU}, Borthwick--Uribe have introduced another approach to generalize the K\"ahler quantization to non-integrable almost complex structures by using the metric Laplacian of the connection on the prequantum line bundle instead of Spin${}^c$ Dirac operator. Their approach is called the almost K\"ahler quantization. In the almost K\"ahler quantization of the Kodaira--Thurston manifold, Kirwin--Uribe and Egorov separately constructed an analog of the theta function as an element of the quantum Hilbert space~\cite{Egorov2,Kirwin--Uribe}. In \cite{Egorov1}, Egorov also showed the similar result for Lagrangian $T^2$-fibrations on $T^2$ with zero Euler class.

Secondly, let us consider the case where $J$ is not integrable. Even in this case, we obtain the following theorem which is a combination of Theorems~\ref{main4} and~\ref{main5}.
\begin{Theorem}\label{main non-integrable case}
Under the above setting, for each $t>0$, there exists an orthogonal family~\smash{$\big\{{\tilde \vartheta}^t_b\big\}_{b\in B_{\rm BS}}$} of $L^2$-sections on $L^{\otimes N}$ indexed by $B_{\rm BS}$ such that
\begin{enumerate}\itemsep=0pt
\item[$(1)$] each ${\tilde \vartheta}^t_b$ converges as a delta-function section supported on $\pi^{-1}(b)$ as $t\to \infty$ in the above sense, and
\item[$(2)$] $\displaystyle \lim_{t\to \infty}\big\|D^t {\tilde \vartheta}^t_b\big\|_{L^2}=0$.\label{3}
\end{enumerate}
\end{Theorem}
We also give \smash{$\big\{\tilde \vartheta^t_b\big\}_{b\in B_{\rm BS}}$} explicitly in Section~\ref{approximation}.

When $M$ is compact, the index of the Spin${}^c$ Dirac operator $D:=D^1$ can be considered and it coincides with the number of Bohr--Sommerfeld points. See \cite{FFY1}. Moreover, by the Spin${}^c$ Dirac vanishing theorem due to Borthwick--Uribe~\cite{BU}, $\ker ( D|_{\wedge^{\rm odd}T^*M^{0,1}\otimes L^{\otimes N}})$ vanishes for a sufficiently large $N$. So, \eqref{3} in Theorem~\ref{main non-integrable case} suggests that we can interpret the Hilbert space generated by \smash{$\big\{{\tilde \vartheta}^t_b\big\}_{b\in B_{\rm BS}}$} as an approximation to the quantum Hilbert space in the Spin${}^c$ quantization for a sufficiently large $N$.

In Kostant--Souriau's formulation of geometric quantization, there is a systematic method to associate an operator on the space of $L^2$-sections of $L^{\otimes N}$, called the prequantum operator, with a~smooth function on $(M,\omega)$. When $(M,\omega)$ is not compact, this induces a nontrivial representation of the Poisson algebra consisting of certain functions to the quantum Hilbert space.
One of the advantages of our setting is that it enable us to deal with this representation of the Poisson algebra concretely by using the complete orthogonal system of the quantum Hilbert space given in this paper, which we will discuss somewhere.

The idea used in this paper is quite simple, and consists of two key facts.
The one is about the symmetry by the fundamental group and the other is about the integrability of almost complex structures.
Namely, the first key fact is Corollary~\ref{complete base} which claims that any Lagrangian torus fibration $\pi\colon (M,\omega)\to B$ with complete base $B$ and a prequantum line bundle~${\bigl(L,\nabla^L\bigr)\to (M,\omega )}$ can be obtained as the quotient of a $\pi_1(B)$-action on the standard Lagrangian fibration $(M_0, \omega_0):=( \R^n\times T^n, \sum_{i=1}^n{\rm d}x_i\wedge {\rm d}y_i)\to \R^n$ with standard prequantum line bundle $\bigl( L_0, \nabla^{L_0}\bigr):=\bigl(\R^n\times T^n\times \C , d-2\pi\sqrt{-1}\sum_{i=1}^nx_i {\rm d}y_i\bigr)$, where $(x_1,\ldots ,x_n)$ and $(y_1,\ldots ,y_n)$ are the coordinates of $\R^n$ and $T^n$, respectively. In particular, any compatible almost complex structure $J$ on $(M,\omega)$ is induced from a $\pi_1(B)$-equivariant one on $(M_0, \omega_0)$. Since the set of compatible almost complex structures on $(M_0, \omega_0)$ corresponds one-to-one to the set of smooth maps from $M_0$ to the Siegel upper half space, it enables us to describe the Spin${}^c$ Dirac operator associated with $J$ in terms of the corresponding map. We show that there exists a $\pi_1(B)$-invariant compatible almost complex structure whose corresponding map is invariant along the fiber (see Lemma~\ref{invariant J}). For the Spin${}^c$ Dirac operator $D$ associated with such an almost complex structure $J$, we consider the problem on existence of nontrivial sections of $L_0^{\otimes N}$ contained in the kernel of $D$. By taking the Fourier series expansion of a section $s$ of $L_0^{\otimes N}$ with respect to the fiber coordinates, the equation $Ds=0$ can be reduced to a system of partial differential equations on $\R^n$.

The other key fact is Proposition~\ref{prop commutativity2} which gives a necessary and sufficient condition in order that the system of partial differential equations has nontrivial solutions and also shows that it is equal to the integrability condition for $J$, i.e., $\left(M_0, \omega_0, J\right)$ is K\"ahler. Moreover, in this case, we give a family of $\pi_1(B)$-equivariant solutions of $Ds=0$ indexed by the Bohr--Sommerfeld points, each of which is expressed by the formal Fourier series. If they converge absolutely and uniformly on any compact set and form square integrable sections, this gives a complete orthogonal system of the space of holomorphic $L^2$-sections of \smash{$\bigl(L,\nabla^L\bigr)^{\otimes N}\to (M,N\omega ,J)$}. We also give a sufficient condition for their convergence and square integrability. Even if $J$ is not integrable, by considering an approximation of $D$, we can obtain an orthogonal family of $L^2$-sections of $L^{\otimes N}$ indexed by the Bohr--Sommerfeld points $B_{\rm BS}$.

The limit used in this paper is slightly different from the adiabatic limit in Riemannian geometry. When a fiber bundle $\pi\colon M\to B$ and a Riemannian metric $g$ on $M$ are given, we can consider the decomposition $(TM, g)=(V, g_V)\oplus (H, g_H)$, where $V$ is the tangent bundle along the fiber with fiber metric $g_V:=g|_V$ and $H$ is the orthogonal complement of $V$ with respect to $g$ with fiber metric $g_H:=g|_H$. For each $t>0$, we deform $g$ by
$g^t:=g_V\oplus tg_H$. Then, in Riemannian geometry, the adiabatic limit is the procedure for taking the limit of geometric objects associated with
$g^t$ as $t\to \infty$. But, since such a deformation of Riemannian metrics does not fit into our symplectic context, we modify the deformation. Namely, in this paper, we use a one-parameter family \smash{$\big\{ J^t\big\}_{t>0}$} of compatible almost complex structures on $(M,\omega)$ such that the corresponding one-parameter family of Riemannian metrics is \smash{$\big\{ g^t=\frac{1}{t}g_V\oplus tg_H\big\}_{t>0}$}, and investigate the behavior of $\vartheta^t_b$ \big(resp.\ \smash{${\tilde \vartheta}^t_b$}\big) as $t$ goes to $\infty$.

The paper is organized as follows. In Section~\ref{Developing Lagrangian fibration}, we first briefly review some well-known facts about integral affine geometry and Lagrangian fibrations. Then, by using these, we prove Corollary~\ref{complete base}. In Section~\ref{equiv quantization}, we discuss the $\pi_1(B)$-equivariant Spin${}^c$ quantization of ${(M_0, \omega_0)\to \R^n}$
with standard prequantum line bundle $\bigl( L_0, \nabla^{L_0}\bigr)$ and give a statement of Proposition~\ref{prop commutativity2}. In Section~\ref{the integrable case}, we prove Theorem~\ref{main integrable case} step by step, and explain the relation between $\vartheta^t_b$ and Jacobi's classical theta function. Finally, in Section~\ref{approximation}, we prove Theorem~\ref{main non-integrable case}. The requirements for Fourier series are explained in Appendix~\ref{appendixA}. A proof of Proposition~\ref{prop commutativity2} is given in Appendix~\ref{appendixB} because it is done by a very long direct computation.

\subsection{Notations}
For $x={}^t(x_1, \ldots ,x_n)$ and $y={}^t(y_1, \ldots ,y_n)\in \R^n$, let us denote the standard inner product~${\sum_{i=1}^nx_iy_i}$ by $x\cdot y$. We use the notation $\partial_{x_i}$ for $\frac{\partial}{\partial x_i}$. We also use the following notations:
\begin{align*}
&T^n:=(\R/\Z)^n, \qquad
(M_0, \omega_0):=\left( \R^n\times T^n, \sum_{i=1}^n{\rm d}x_i\wedge {\rm d}y_i\right),\\
&\bigl( L_0, \nabla^{L_0}\bigr):=\left(\R^n\times T^n\times \C , d-2\pi\sqrt{-1}\sum_{i=1}^nx_i {\rm d}y_i\right),
\end{align*}
where $(x_1,\ldots ,x_n)$ and $(y_1,\ldots ,y_n)$ are the coordinates of $\R^n$ and $T^n$, respectively. 
In this paper, all manifolds and maps are supposed to be smooth unless otherwise stated. When we use the term \lq\lq group action", we mean \lq\lq left group action" unless otherwise specified.

\section{Unfolding Lagrangian fibrations}\label{Developing Lagrangian fibration}
\subsection{Integral affine structures}
Let $B$ be a manifold.
\begin{Definition}\label{def integral affine structure}
An \emph{integral affine atlas} of $B$ is an atlas $\{ (U_\alpha,\phi_\alpha)\}$ of $B$ whose coordinate transformation $\phi_\alpha\circ \phi_\beta^{-1}$ on each non-empty overlap $U_{\alpha\beta}:=U_\alpha \cap U_\beta$ is an integral affine transformation. Namely, on each non-empty overlap $U_{\alpha\beta}:=U_\alpha \cap U_\beta$, there exist locally constant maps $A_{\alpha\beta}\colon U_{\alpha\beta}\to \GL_n(\Z)$ and $c_{\alpha\beta}\colon U_{\alpha\beta}\to \R^n$ such that $\phi_\alpha\circ \phi_\beta^{-1}$ is of the form \smash{$\phi_\alpha\circ \phi_\beta^{-1} (x)=A_{\alpha\beta}x+c_{\alpha\beta}$}.
Two integral affine atlases $\{ (U_\alpha, \phi_\alpha)\}$ and $\{ (U'_\beta ,\phi'_\beta)\}$ of $B$ are said to be \emph{equivalent} if on each non-empty overlap $U_\alpha\cap U'_\beta$, \smash{$\phi_\alpha\circ (\phi'_{\beta})^{-1}$} is an integral affine transformation.
An \emph{integral affine structure} on $B$ is an equivalence class of integral affine atlases of~$B$. A manifold equipped with integral affine structure is called an \emph{integral affine manifold}.
\end{Definition}

\begin{Example}\label{affine Euclidean}
An $n$-dimensional Euclidean space $\R^n$ is equipped with a natural integral affine structure.
\end{Example}

Let us give examples of integral affine manifolds obtained from integral affine actions on $\R^n$.
\begin{Example}\label{affine torus}\qquad
\begin{itemize}\itemsep=0pt
\item[(1)] Let $v_1,\ldots ,v_n\in \R^n$ be a linear basis of $\R^n$ and $C=(v_1\ \cdots\ v_n)\in \GL_n(\R)$ the matrix whose $i$th column vector is $v_i$ for $i=1,\ldots , n$. $\Z^n$ acts on $\R^n$ by $\rho_\gamma (x):=x+C\gamma$
for~${\gamma\in \Z^n}$ and $x\in \R^n$. Since the action preserves the natural integral affine structure on $\R^n$, the quotient space, which is topologically $T^n$, is equipped with an integral affine structure.
\item[(2)] Let $\lambda\in \N$ be a positive integer and $a, b\in \R_{>0}$ positive real numbers. Define the $\Z^2$-action on $\R^2$ as follows. First, for the standard basis $e_1$, $e_2$ of $\Z^2$, let us define the integral affine transforms $\rho_{e_1}$, $\rho_{e_2}$ by
\[
\rho_{e_1}(x):=x+
\begin{pmatrix}
a \\
0
\end{pmatrix},\qquad
\rho_{e_2}(x):=
\begin{pmatrix}
1 & \lambda \\
0 & 1
\end{pmatrix}
x+
\begin{pmatrix}
0 \\
b
\end{pmatrix}
\]
for $x\in \R^2$. Since $\rho_{e_1}$ and $\rho_{e_2}$ are commutative, they form the $\Z^2$-action on $\R^2$ by
\[
\rho_\gamma (x):=\rho_{e_1}^{\gamma_1}\circ \rho_{e_2}^{\gamma_2}(x)
\]
for each $\gamma={}^t(\gamma_1,\gamma_2)\in \Z^2$. By the same manner as in (1),
the quotient space is equipped with an integral affine structure. It is shown in~\cite[Theorem A]{Mis} that the quotient space is topologically $T^2$, but the induced integral affine structure is not isomorphic to that obtained in (1) for $n=2$ and there are only these two integral affine structures on $T^2$ up to isomorphism.
\end{itemize}
\end{Example}

\begin{Example}
For $\gamma={}^t(\gamma_1, \gamma_2, \gamma_3)$, $\gamma'={}^t(\gamma'_1, \gamma'_2, \gamma'_3)\in \Z^3$, define the product $\gamma\circ \gamma' \in \Z^3$ by
\[
\gamma\circ \gamma':=
\begin{pmatrix}
1 & \hphantom{-}0 & \hphantom{-}0 \\
0 & \hphantom{-}0 & -1 \\
0 & -1 & \hphantom{-}0
\end{pmatrix}^{\gamma_1}
\gamma'+\gamma.
\]
Then, $\Z^3$ with product $\circ$ is a non abelian group $\bigl(\Z^3,\circ \bigr)$. $\bigl(\Z^3,\circ \bigr)$ acts on $\R^3$ by
\[
\rho_\gamma(x):=
\begin{pmatrix}
1 & \hphantom{-}0 & \hphantom{-}0 \\
0 & \hphantom{-}0 & -1 \\
0 & -1 & \hphantom{-}0
\end{pmatrix}^{\gamma_1}
x+\gamma.
\]
The action preserves the natural integral affine structure on $\R^3$. Therefore, the quotient space~${\R^3/\bigl(\Z^3,\circ \bigr)}$ is equipped with the integral affine structure induced from that of $\R^3$.
\end{Example}

\begin{Example}
Let $n\ge 2$. For $\gamma={}^t(\gamma_1, \ldots, \gamma_n)$, $\gamma'={}^t(\gamma'_1, \ldots, \gamma'_n)\in \Z^n$, define the product~${\gamma\circ \gamma' \in \Z^n}$ by
\[
\gamma\circ \gamma':=
\begin{pmatrix}
1 & & & \\
 & (-1)^{\gamma_1} & & \\
 & & \ddots & \\
 & & & (-1)^{\gamma_{n-1}}
\end{pmatrix}
\gamma'+\gamma.
\]
Then, $\Z^n$ with product $\circ$ is a non abelian group $(\Z^n,\circ )$. $(\Z^n,\circ )$ acts on $\R^n$ by
\[
\rho_\gamma(x):=
\begin{pmatrix}
1 & & & \\
 & (-1)^{\gamma_1} & & \\
 & & \ddots & \\
 & & & (-1)^{\gamma_{n-1}}
\end{pmatrix}
x+\gamma.
\]
The action preserves the natural integral affine structure on $\R^n$. Therefore, the quotient space~${\R^n/(\Z^n,\circ )}$ is equipped with the integral affine structure induced from that of $\R^n$. For~${n\!=\!2}$, the quotient space is topologically a Klein bottle.
\end{Example}

\begin{Example}\label{Jordan block}
Let $n\ge 2$ and $\lambda_1,\ldots, \lambda_{n-1}\in \Z$. For $\gamma={}^t(\gamma_1, \ldots, \gamma_n)$, $\gamma'={}^t(\gamma'_1, \ldots, \gamma'_n)\in \Z^n$, define the product $\gamma\circ \gamma' \in \Z^n$ by
\[
\gamma\circ \gamma' :=
\begin{pmatrix}
1 & \lambda_1 & & & \\
 & 1 & \lambda_2 & & \\
 & & \ddots & \ddots & \\
 & & & 1 & \lambda_{n-1}\\
 & & & & 1
\end{pmatrix}^{\gamma_n}
\gamma'+\gamma.
\]
$\Z^n$ with product $\circ$ is a group $(\Z^n,\circ )$, which is non abelian for $n\ge 3$. $(\Z^n,\circ )$ acts on $\R^n$ by
\[
\rho_\gamma (x):=
\begin{pmatrix}
1 & \lambda_1 & & & \\
 & 1 & \lambda_2 & & \\
 & & \ddots & \ddots & \\
 & & & 1 & \lambda_{n-1}\\
 & & & & 1
\end{pmatrix}^{\gamma_n}
x+\gamma.
\]
Then, the quotient space $\R^n/(\Z^n,\circ )$ is equipped with the integral affine structure induced from that of $\R^n$. In the case where $n=2$ and $\lambda_1>0$, it coincides with the one given in Example~\ref{affine torus}\,(2) with $a=b=1$.
\end{Example}

\begin{Example}\label{incomplete}
Let $\Z /4\Z\cong \left\{
\pm \left(\begin{smallmatrix}
1 & 0 \\
0 & 1
\end{smallmatrix}\right),
\pm \left(\begin{smallmatrix}
0 & -1 \\
1 & \hphantom{-}0
\end{smallmatrix}\right)
\right\}$ act on $\bigl(\R^2\bigr)^n\setminus \{ 0\}$ naturally. Then, the quotient space is a non-compact manifold and equipped with the integral affine structure induced from that of $\bigl(\R^2\bigr)^n\setminus \{0\}$.
\end{Example}

Let $B$ be an $n$-dimensional connected integral affine manifold, $p\colon \tilde B\to B$ the universal covering of $B$. It is clear that $\tilde B$ is also equipped with the integral affine structure so that $p$ is an integral affine map. We set $\Gamma:=\pi_1(B)$. $\Gamma$ acts on $\tilde B$ from the right as a deck transformation. For each $\gamma \in \Gamma$ we denote by $\sigma_\gamma$ the inverse of the deck transformation corresponding to $\gamma$. Then, $\sigma\colon \gamma\mapsto \sigma_\gamma$ defines a left action \smash{$\sigma\in \Hom\bigl(\Gamma ,\Aut\bigl(\tilde B\bigr)\bigr)$}.
The following proposition is well known in affine geometry. See \cite[p.~641]{GH} for a proof.
\begin{Proposition}\label{developing map}
There exists an integral affine immersion $\dev\colon \tilde B\to \R^n$ and a homomorphism~${\rho\colon \Gamma\to \GL_n(\Z)\ltimes \R^n}$ such that the image of $\dev$ is an open set of $\R^n$ and $\dev$ is equivariant with respect to $\sigma$ and $\rho$. Such an integral affine immersion is unique up to the composition of an integral affine transformation on $\R^n$.
\end{Proposition}
We will prove a version of this proposition (see Proposition~\ref{Lagrangian fibration on simply connected base}) when B is equipped with a~Lagrangian fibration on it in Section~\ref{complete bases}.

\begin{Proposition}
Let $B$, $p\colon \tilde B\to B$, $\dev\colon \tilde B\to \R^n$, and $\rho\colon \Gamma\to \GL_n(\Z)\ltimes \R^n$ be as in Proposition~{\rm\ref{developing map}}. Suppose that $B$ is compact and the $\Gamma$-action $\rho$ on $\R^n$ is properly discontinuous. Then, $\dev$ is surjective.
\end{Proposition}
\begin{proof}
We denote the image of $\dev$ by $O$. By Proposition~\ref{developing map}, $O$ is an open set in $\R^n$. So, it is sufficient to show that $O$ is closed in $\R^n$. Since the $\Gamma$-action $\rho$ on $\R^n$ is properly discontinuous, the quotient space $\R^n /\Gamma$ becomes a Hausdorff space and the natural projection $q\colon \R^n\to \R^n/\Gamma$ is continuous. $O$ is preserved by the $\Gamma$-action $\rho$ on $\R^n$ since $\dev$ is $\Gamma$-equivariant. Then, $\dev$ induces a continuous surjective map \smash{${\overline \dev}\colon B={\tilde B}/\Gamma \to O/\Gamma$}. Since $B$ is compact, $O/\Gamma$ is a compact subset in the Hausdorff space $\R^n/\Gamma$. In particular, it is also closed. Hence, $O=q^{-1}\left( O/\Gamma\right)$ is also closed in $\R^n$.
\end{proof}

\begin{Corollary}
Let $B$, $p\colon \tilde B\to B$, $\dev\colon \tilde B\to \R^n$, and $\rho\colon \Gamma\to \GL_n(\Z)\ltimes \R^n$ be as in Proposition~{\rm\ref{developing map}} and assume that $B$ compact. If the image of $\rho$ is included in $(\GL_n(\Z)\cap {\rm O}(n))\ltimes \R^n$ and the subgroup $\rho(\Gamma )$ of $(\GL_n(\Z)\cap {\rm O}(n))\ltimes \R^n$ is discrete, then $\dev$ is surjective.
\end{Corollary}
\begin{proof}
It follows from \cite[Theorem~3.1.3]{Wolf}.
\end{proof}

\begin{Definition}
The integral affine immersion $\dev$ is called a \emph{developing map}. $B$ is said to be \emph{complete} if $\dev$ is bijective. $B$ is said to be \emph{incomplete} if $B$ is not complete.
\end{Definition}
\begin{Example}
All of the above examples are complete other than Example~\ref{incomplete} for $n\ge 2$.
\end{Example}

\begin{Example}
Let $B$ be an $n$-dimensional compact integral affine manifold $B$ with integral affine atlas $\{ (U_\alpha,\phi_\alpha)\}$ as in Definition~\ref{def integral affine structure}. If on each non-empty overlap $U_{\alpha\beta}$, the Jacobi matrix of $\phi_\alpha\circ \phi_\beta^{-1}$ is contained in $\GL_n(\Z)\cap {\rm O}(n)$, then $B$ has a flat Riemannian metric. Hence, by Bieberbach's theorem~\cite{Bie1, Bie2}, $B$ is finitely covered by $T^n$. In particular, $B$ is complete. For flat Riemannian manifolds, see \cite[Chapter~3]{Wolf}.
\end{Example}

\subsection{Lagrangian fibrations}\label{Lagrangian fibration}
In this section, let us recall Lagrangian fibrations and explain how integral affine structures are associated with Lagrangian fibrations. After that, let us recall their classification by Duistermaat. For more details, see \cite{Du,MR2587398,Y4}.

Let $(M,\omega)$ be a symplectic manifold.
\begin{Definition}
A map $\pi\colon (M,\omega)\to B$ from $(M,\omega)$ to a manifold $B$ is called a {\it Lagrangian fibration} if $\pi$ is a fiber bundle whose fiber is a Lagrangian submanifold of $(M,\omega)$.
\end{Definition}

\begin{Example}\label{local model for Lagrangian fibration}
Let $T^n=(\R/\Z)^n$ be an $n$-dimensional torus.
$\R^n\times T^n$ admits a standard symplectic structure $\omega_0=\sum_i{\rm d}x_i\wedge {\rm d}y_i$, where $(x_1,\ldots ,x_n)$ and $(y_1,\ldots ,y_n)$ are the coordinates of $\R^n$ and $T^n$, respectively. Then, the projection $\pi_0\colon (\R^n\times T^n,\omega_0)\to \R^n$ to $\R^n$ is a Lagrangian fibration.
\end{Example}

The following theorem shows that Example~\ref{local model for Lagrangian fibration} is the local model of Lagrangian fibrations.
\begin{Theorem}\label{Arnold-Liouville}
Let $\pi\colon (M,\omega)\to B$ be a Lagrangian fibration with compact, path-connected fibers. Then, for each $b\in B$, there exists a chart $(U,\phi)$ of $B$ containing $b$ and a symplectomorphism $\varphi\colon \bigl( \pi^{-1}(U), \omega|_{\pi^{-1}(U)}\bigr)\to ( \phi (U)\times T^n, \omega_0)$ such that the following diagram commutes:
\[
\xymatrix{
\bigl(\pi^{-1}(U),\omega|_{\pi^{-1}(U)}\bigr)\ar[d]^{\pi}\ar[r]^{\varphi} & (\phi (U)\times T^n, \omega_0)\ar[d]^{\pi_0} \\
U\ar[r]^{\phi} & \phi(U).
}
\]
\end{Theorem}
\begin{proof}
In \cite[Section~2]{Du}, Duistermaat showed that any Lagrangian fibration with compact, path-connected fibers is locally identified with a regular part of a completely integrable Hamiltonian system. Theorem~\ref{Arnold-Liouville} follows from this fact together with Arnold--Liouville's theorem which claims that a regular part of a completely integrable Hamiltonian system is locally identified with $\pi_0\colon (\R^n\times T^n,\omega_0)\to \R^n$. For Arnold--Liouville's theorem, see~\cite{Ar,Du,MR2587398,S}.
\end{proof}

In particular, Theorem~\ref{Arnold-Liouville} says that any Lagrangian fibration with compact, path-connected fibers has a torus as its fiber.

In this paper, we consider only Lagrangian fibrations with compact, path-connected fibers. In the rest of this paper, \lq\lq Lagrangian fibration" refers to a Lagrangian fibration having compact, path-connected fibers.

Now we investigate automorphisms of the local model. By the direct computation shows the following lemma. See also \cite[Lemma~2.5]{S}.

\begin{Lemma}\label{fiber-preserving symplectomorphism}
Let $\varphi\colon (\R^n\times T^n,\omega_0)\to (\R^n\times T^n,\omega_0)$ be a fiber-preserving symplectomorphism of $\pi_0\colon (\R^n\times T^n,\omega_0)\to \R^n$ which covers a map $\phi\colon \R^n\to \R^n$. Then, there exist a~matrix~${A\in \GL_n(\Z)}$, a constant $c\in \R^n$, and a map $u\colon \R^n\to T^n$ with ${}^tA \J u$ symmetric
such that~$\varphi$ is written as
\[
\varphi (x,y)=\bigl( Ax+c,{}^tA^{-1}y+u(x)\bigr)
\]
for any $(x,y)\in \R^n\times T^n$, where $\J u$ is the Jacobi matrix of $u$.
\end{Lemma}

By Theorem~\ref{Arnold-Liouville} and Lemma~\ref{fiber-preserving symplectomorphism}, we can obtain the following proposition.
\begin{Proposition}\label{integral affine structure}
Let $\pi\colon (M,\omega)\to B$ be a Lagrangian fibration. Then, there exists an atlas~${\{ (U_\alpha,\phi_\alpha)\}_{\alpha\in A}}$ of $B$ and for each $\alpha\in A$ there exists a symplectomorphism
\[
\varphi_\alpha\colon\ \bigl(\pi^{-1}(U_\alpha),\omega|_{\pi^{-1}(U_\alpha)}\bigr)\to (\phi_\alpha(U_\alpha)\times T^n,\omega_0)
\]
 such that the following diagram commutes:
\[
\xymatrix{
\bigl(\pi^{-1}(U_\alpha),\omega|_{\pi^{-1}(U_\alpha)}\bigr)\ar[d]^{\pi}\ar[r]^{\varphi_\alpha} & (\phi_\alpha (U_\alpha)\times T^n, \omega_0)\ar[d]^{\pi_0} \\
U_\alpha\ar[r]^{\phi_\alpha} & \phi_\alpha(U_\alpha).
}
\]
Moreover, on each non-empty overlap $U_{\alpha\beta}$ there exist locally constant maps ${A_{\alpha\beta}}\colon { U_{\alpha\beta}}\!\to\! \GL_n(\Z)$, $c_{\alpha\beta}\colon U_{\alpha\beta}\to \R^n$, and a map $u_{\alpha\beta}\colon U_{\alpha\beta}\to T^n$ with ${}^tA_{\alpha\beta}\J \bigl(u_{\alpha\beta}\circ \phi_\beta^{-1}\bigr)$ symmetric, such that the overlap map is written as
\begin{equation}\label{transition}
\varphi_\alpha\circ\varphi_\beta^{-1}(x,y)=\bigl(A_{\alpha\beta}x+c_{\alpha\beta},{}^tA^{-1}_{\alpha\beta}y+u_{\alpha\beta}\circ \phi_\beta^{-1}(x)\bigr)
\end{equation}
for any $(x,y)\in \phi_\beta(U_{\alpha\beta})\times T^n$.
\end{Proposition}

Proposition~\ref{integral affine structure} implies that the base manifold of a Lagrangian fibration has an integral affine structure.
Conversely, suppose that a manifold $B$ admits an integral affine structure and let~$\{ (U_\alpha,\phi_\alpha)\}_{\alpha\in A}$ be an integral affine atlas of $B$. Then, we can construct a Lagrangian fibration on $B$ in the following way. For each $\alpha \in A$, let $\overline \phi_\alpha\colon T^*B|_{U_\alpha}\to \phi_\alpha(U_\alpha) \times \R^n$ be the local trivialization of the cotangent bundle $T^*B$ induced from $(U_\alpha,\phi_\alpha)$. On each nonempty overlap~$U_{\alpha\beta}$, suppose that $\phi_{\alpha}\circ \phi_{\beta}^{-1}$ is written by $\phi_{\alpha}\circ \phi_{\beta}^{-1}(x)=A_{\alpha\beta}x+c_{\alpha\beta}$ as in Definition~\ref{def integral affine structure}. Then, the overlap map is written as
\begin{equation}\label{overlap map of T*B}
\overline\phi_\alpha\circ \bigl(\overline\phi_\beta\bigr)^{-1}(x,y)=\bigl(A_{\alpha\beta}x+c_{\alpha\beta},{}^tA^{-1}_{\alpha\beta}y\bigr).
\end{equation}
Since $A_{\alpha\beta}$ is in $\GL_n(\Z)$, \eqref{overlap map of T*B} preserves the integer lattice $\Z^n$ in each fiber $\R^n$. Hence, it induces the fiber-preserving symplectomorphism from $\pi_0\colon(\phi_\beta(U_{\alpha\beta})\times T^n,\omega_0)\to \phi_\beta(U_{\alpha\beta})$ to $\pi_0\colon (\phi_\alpha(U_{\alpha\beta})\times T^n,\omega_0)\to \phi_\alpha(U_{\alpha\beta})$ which covers $\phi_{\alpha}\circ \phi_{\beta}^{-1}$. By patching $\{\pi_0\colon (\phi_{\alpha}(U_\alpha)\times T^n, \omega_0)$ $\to \phi_\alpha(U_\alpha)\}_{\alpha\in A}$ together by these symplectomorphisms, we obtain a new Lagrangian fibration $\pi_{\rm can}\colon (M_{\rm can},\omega_{\rm can})\to B$, namely,
\[
(M_{\rm can},\omega_{\rm can}):=\coprod_{\alpha\in A}(\phi_\alpha(U_\alpha)\times T^n, \omega_0)/{\sim}
\]
and
\[
\pi_{\rm can}([x_\alpha,y_\alpha]):=\phi_\alpha^{-1}(x_\alpha)
\]
for $(x_\alpha,y_\alpha)\in \phi_\alpha(U_\alpha)\times T^n$. This construction does not depend on the choice of an equivalent integral affine atlas and depends only on the integral affine structure on $B$.
We call $\pi_{\rm can}\colon (M_{\rm can},\omega_{\rm can})\to B$ the {\it canonical model}.
%
We summarize the above argument to the following proposition.
\begin{Proposition}\label{necessary and sufficient condition}
Let $B$ be a manifold. $B$ is a base space of a Lagrangian fibration if and only if $B$ admits an integral affine structure.
\end{Proposition}

Let us give a classification theorem of Lagrangian fibrations in the required form in this paper.
Let $\pi\colon (M,\omega)\to B$ be a Lagrangian fibration. Then, associated with $\pi\colon (M,\omega)\to B$, $B$ has an integral affine atlas $\{ (U_\alpha,\phi_\alpha)\}_{\alpha\in A}$ as in Proposition~\ref{integral affine structure}. Let $\pi_{\rm can}\colon (M_{\rm can},\omega_{\rm can})\to B$ be the canonical model associated with the integral affine structure on $B$. On each $U_\alpha$, let
\[
\varphi_\alpha\colon \ \bigl(\pi^{-1}(U_\alpha),\omega|_{\pi^{-1}(U_\alpha)}\bigr)\to (\phi_\alpha(U_\alpha)\times T^n,\omega_0)
\]
 be a local trivialization of $\pi\colon (M,\omega)\to B$ as in Proposition~\ref{integral affine structure}, and
 \[
 \overline\phi_\alpha\colon \ \bigl(\pi_{\rm can}^{-1}(U_\alpha), \omega_{\rm can}\bigr)\to (\phi_\alpha(U_\alpha)\times T^n,\omega_0)
 \]
 be the local trivialization of $\pi_{\rm can}\colon (M_{\rm can},\omega_{\rm can})\to B$ naturally induced from $(U_\alpha,\phi_\alpha)$ as explained above.\footnote{Here we use the same notation as the local trivialization of $T^*B$ because we have no confusion.} Then their composition
\[
h_\alpha:={\overline\phi_\alpha}^{-1}\circ \varphi_\alpha\colon\ \bigl(\pi^{-1}(U_\alpha),\omega|_{\pi^{-1}(U_\alpha)}\bigr)\to \bigl(\pi_{\rm can}^{-1}(U_\alpha), \omega_{\rm can}\bigr)
\]
gives a local identification between them.
On each $U_\alpha\cap U_\beta$, suppose that $\varphi_\alpha\circ\varphi_\beta^{-1}$ is written as in \eqref{transition}. Then, \smash{$h_\alpha\circ h_\beta^{-1}$} can be written as
\[
h_\alpha\circ h_\beta^{-1}(p)={\overline\phi_\alpha}^{-1}\bigl(A_{\alpha\beta}x+c_{\alpha\beta},{}^tA_{\alpha\beta}^{-1}y+u_{\alpha\beta}(\pi(p))\bigr),
\]
where $\overline\phi_\beta(p)=(x,y)$. $u_{\alpha\beta}$ induces the local section $\tilde u_{\alpha\beta}$ of $\pi_{\rm can}\colon (M_{\rm can},\omega_{\rm can})\to B$ on $U_{\alpha\beta}$ by
\[
\tilde u_{\alpha\beta}(b):=[\phi_\alpha(b),u_{\alpha\beta}(b)]
\]
for $b\in U_{\alpha\beta}$. It is easy to see that $\tilde u_{\alpha\beta}$ satisfies $\tilde u_{\alpha\beta}^*\omega_{\rm can}=0$. A section with this condition is said to be \emph{Lagrangian}.

Let $\SCS$ be the sheaf of germs of Lagrangian section of $\pi_{\rm can}\colon (M_{\rm can},\omega_{\rm can})\to B$. $\SCS$ is the sheaf of Abelian groups since the fiber of $\pi_{\rm can}\colon (M_{\rm can},\omega_{\rm can})\to B$ has the structure of an Abelian group by construction. By definition $\{ \tilde u_{\alpha\beta}\}$ forms a \v{C}ech one-cocycle on $B$ with coefficients in $\SCS$. The cohomology class determined by $\{ \tilde u_{\alpha\beta}\}$ does not depend on the choice of a specific integral affine atlas and depends only on $\pi\colon (M,\omega)\to B$. We denote the cohomology class by~${u \in H^1(B;\SCS)}$.
$u$ is called the \emph{Chern class} of $\pi\colon (M,\omega)\to B$ in~\cite{Du}.

Lagrangian fibrations on the same integral affine manifold are classified with the Chern classes.
\begin{Theorem}[\cite{Du}]\label{classification}
For two Lagrangian fibrations $\pi_1\colon (M_1,\omega_1)\to B$ and ${\pi_2\colon (M_2,\omega_2)\to B}$ on the same integral affine manifold $B$, there exists a fiber-preserving symplectomorphism between them which covers the identity if and only if their Chern classes $u_1$ and $u_2$ agree with each other. Moreover, if an integral affine manifold $B$ and the cohomology class $u\in H^1(B;\SCS)$ are given, then there exists a Lagrangian fibration $\pi\colon (M,\omega)\to B$ that realizes them.
\end{Theorem}
\begin{Remark}
By the construction of $u$, there exists a fiber-preserving symplectomorphism between $\pi \colon (M,\omega)\to B$ and $\pi_{\rm can}\colon (M_{\rm can},\omega_{\rm can})\to B$ that covers the identity of $B$ if and only if $u$ vanishes. In particular, if $u$ vanishes, $\pi \colon (M,\omega)\to B$ possesses a Lagrangian section since~${\pi_{\rm can}\colon (M_{\rm can},\omega_{\rm can})\to B}$ has the zero section which is Lagrangian. Conversely, it can be shown that any Lagrangian fibration with Lagrangian section is identified with the canonical model.
\end{Remark}

\subsection{Lagrangian fibrations with complete bases}\label{complete bases}
Assume that $\pi \colon (M,\omega)\to B$ is a Lagrangian fibration with $n$-dimensional connected base mani\-fold~$B$, $p\colon \tilde B\to B$ the universal covering of $B$. We denote by \smash{$\tilde \pi\colon \bigl(\tilde M, \tilde \omega\bigr) \to \tilde B$} the pullback of $\pi\colon (M,\omega)\to B$ to $\tilde B$. Let $\Gamma$ be the fundamental group of $B$ and $\sigma \in \Hom\bigl( \Gamma ,\Aut\bigl(\tilde B\bigr)\bigr)$ the action of $\Gamma$ defined as the inverse of the deck transformation as in Proposition~\ref{developing map}. By definition, \smash{$\tilde M$} admits a natural lift of $\sigma$ which preserves $\tilde \omega$. The $\Gamma$-action on \smash{$\bigl(\tilde M, \tilde \omega\bigr)$} is denoted by $\tilde \sigma$. By Proposition~\ref{developing map}, we have a developing map $\dev\colon \tilde B\to \R^n$ and the homomorphism $\rho\colon \Gamma\to \GL_n(\Z)\ltimes \R^n$. We denote the image of $\dev$ by $O$. Note that the $\Gamma$-action $\rho$ on $\R^n$ preserves $O$ since $\dev$ is $\Gamma$-equivariant.
\begin{Proposition}\label{Lagrangian fibration on simply connected base}
There exists a Lagrangian fibration $\pi'\colon (M',\omega')\to O$, a fiber-preserving symplectic immersion \smash{$\tilde \dev \colon \bigl(\tilde M,\tilde \omega\bigr)\to (M',\omega')$} which covers $\dev$, and a lift $\tilde \rho$ of the $\Gamma$-action $\rho$ on $O$ to $(M',\omega')$ such that \smash{$\tilde \dev$} is $\Gamma$-equivariant with respect to $\tilde \sigma$ and $\tilde \rho$.
\end{Proposition}
\begin{proof}
By Proposition~\ref{necessary and sufficient condition}, $B$ admits an integral affine structure determined by $\pi$, and it also induces the integral affine structure on $\tilde B$. Let $\{ (U_\alpha,\phi_\alpha'')\}$ be the integral affine atlas of~$\tilde B$ and~${\big\{ \bigl(\tilde \pi^{-1}(U_\alpha),\omega|_{\tilde \pi^{-1}(U_\alpha)}, \varphi''_\alpha\bigr) \big\}}$ the local trivializations of \smash{$\tilde \pi\colon \bigl(\tilde M,\tilde \omega\bigr)\to \tilde B$} as in Proposition~\ref{integral affine structure} so that on each non-empty overlap $U_{\alpha\beta}$, there exist locally constant maps ${A_{\alpha\beta}\colon U_{\alpha\beta}\to \GL_n(\Z)}$ and \smash{$c'_{\alpha\beta}\colon U_{\alpha\beta}\to \R^n$}, and a map \smash{$u'_{\alpha\beta}\colon U_{\alpha\beta}\to T^n$} with \smash{${}^tA_{\alpha\beta}\J \bigl(u'_{\alpha\beta}\circ (\phi''_\beta)^{-1}\bigr)$} symmetric such that~\smash{${\varphi''_\alpha\circ(\varphi''_\beta)^{-1}}$} is written as in \eqref{transition}. Then, $A_{\alpha\beta}$'s form a \v{C}ech one-cocycle $\{ A_{\alpha\beta}\}\in C^1(\{U_\alpha\};\allowbreak\GL_n(\Z))$ and defines a cohomology class \smash{$[\{ A_{\alpha\beta}\}]\in H^1\bigl(\tilde B; \GL_n(\Z)\bigr)$}. It is well known, for example, see \cite[Appendix A]{LM}, that \smash{$H^1\bigl(\tilde B; \GL_n(\Z)\bigr)$} is identified with the moduli space of homo\-morphisms from $\pi_1\bigl(\tilde B\bigr)$ to $\GL_n(\Z)$. Since \smash{$\pi_1\bigl(\tilde B\bigr)$} is trivial, there exists a \v{C}ech zero-cocycle $\{ A_\alpha\}\in C^0(\{U_\alpha\};\GL_n(\Z))$ such that $A_{\alpha\beta}=A_\alpha A_\beta^{-1}$ on each $U_{\alpha\beta}$. By using the cocycle we modify the local trivializations \smash{$\big\{ \bigl(\tilde \pi^{-1}(U_\alpha),\omega|_{\tilde \pi^{-1}(U_\alpha)}, \varphi''_\alpha\bigr) \big\}$} and the integral affine atlas $\{(U_\alpha, \phi_\alpha'')\}$ by replacing $\varphi''_\alpha$, $\phi''_\alpha$ by
\[
\varphi'_\alpha(\tilde p):=\bigl(A_\alpha^{-1}\times {}^tA_\alpha\bigr)\circ \phi_\alpha''(\tilde p),\qquad \phi'_\alpha:=A_\alpha^{-1}\phi''_\alpha
\]
for each $\alpha \in A$, respectively. Then, on each $U_{\alpha\beta}$, $\varphi'_\alpha\circ (\varphi'_\beta)^{-1}$ is written as
\[
\varphi'_\alpha\circ (\varphi'_\beta)^{-1}(\tilde x,y)=\bigl(\tilde x+c_{\alpha\beta},y+u_{\alpha\beta}\circ (\phi'_\beta)^{-1}(\tilde x)\bigr),
\]
where we set $c_{\alpha\beta}:=A_\alpha^{-1}c'_{\alpha\beta}$ and $u_{\alpha\beta}:={}^tA_\alpha u'_{\alpha\beta}$. Then, $c_{\alpha\beta}$'s form a \v{C}ech one-cocycle $\{ c_{\alpha\beta}\}\in C^1(\{U_\alpha\};\R^n)$ and defines a cohomology class \smash{$[\{ c_{\alpha\beta}\}]\in H^1\bigl(\tilde B; \R^n\bigr)$}. By the universal coefficients theorem, $H^1\bigl(\tilde B; \R^n\bigr)$ is identified with $\Hom \bigl(H_1\bigl(\tilde B;\Z\bigr), \R^n\bigr)$, which is trivial. So there exists a~\v{C}ech zero-cocycle $\{ c_\alpha\}\in C^0(\{U_\alpha\};\R^n)$ such that $c_{\alpha\beta}=c_\alpha -c_\beta$ on each $U_{\alpha\beta}$. By using the cocycle, we again modify $\big\{ \bigl(\tilde \pi^{-1}(U_\alpha),\omega|_{\tilde \pi^{-1}(U_\alpha)}, \varphi'_\alpha\bigr) \big\}$ and $\{(U_\alpha, \phi'_\alpha)\}$ by replacing $\varphi'_\alpha$, $\phi_\alpha'$ by
\[
\varphi_\alpha (\tilde p):=\varphi'_\alpha (\tilde p)-(c_\alpha ,0),\qquad
\phi_\alpha\bigl(\tilde b\bigr):=\phi_\alpha'\bigl(\tilde b\bigr)-c_\alpha ,
\]
respectively for each $\alpha \in A$. Then, on each $U_{\alpha\beta}$, $\phi_\alpha$ coincides with $\phi_\beta$ and $\varphi_\alpha\circ \varphi_\beta^{-1}$ is written~as
\[
\varphi_\alpha\circ \varphi_\beta^{-1}(\tilde x, y)=\bigl(\tilde x, y+u_{\alpha\beta}\circ \phi_\beta^{-1}(\tilde x)\bigr).
\]
Now we define the map $\dev\colon \tilde B\to \R^n$ by
\[
\dev \bigl(\tilde b\bigr):=\phi_\alpha \bigl(\tilde b\bigr)
\]
if $\tilde b$ is in $U_\alpha$. It is well defined, and by construction, it is an integral affine immersion whose image is $\displaystyle \cup_{\alpha\in A}\phi_\alpha(U_\alpha)$.
$(M',\omega')$ is defined by
\[
(M',\omega'):=\coprod_{\alpha \in A} (\phi_\alpha(U_\alpha)\times T^n,\omega_0 )/{\sim} ,
\]
where $(x_\alpha, y_\alpha)\in \phi_\alpha(U_\alpha)\times T^n$ and $(x_\beta, y_\beta)\in \phi_\beta(U_\beta)\times T^n$ are in the relation $(x_\alpha, y_\alpha)\sim (x_\beta, y_\beta)$ if they satisfy $(x_\alpha, y_\alpha)=\varphi_\alpha\circ \varphi_\beta^{-1}(x_\beta, y_\beta)$, and $\pi'\colon (M',\omega')\to O$ is defined to be the first projection. \smash{$\tilde \dev\colon \bigl(\tilde M,\tilde \omega\bigr)\to (M',\omega')$} is defined by
\[
\tilde \dev(\tilde p):=[\varphi_\alpha (\tilde p)]
\]
if $\tilde p$ is in $\tilde \pi^{-1}(U_\alpha)$.

Without loss of generality, we can assume that each $U_\alpha$ is connected,and for each $\gamma\in \Gamma$ and~${\alpha\in A}$ there uniquely exists $\alpha'\in A$ such that the deck transformation $\sigma_\gamma$ maps $U_\alpha$ onto $U_{\alpha'}$. Then, its lift $\tilde \sigma_\gamma$ to \smash{$\bigl(\tilde M,\tilde \omega\bigr)$} maps $\tilde \pi^{-1}(U_\alpha)$ to $\tilde \pi^{-1}(U_{\alpha'})$. By Lemma~\ref{fiber-preserving symplectomorphism}, $\phi_{\alpha'}\circ \sigma_\gamma \circ \phi_\alpha^{-1}$
can be written as
\[
\phi_{\alpha'}\circ \sigma_\gamma \circ \phi_\alpha^{-1}(\tilde x)=A_\gamma^{\alpha'\alpha}\tilde x+c_\gamma^{\alpha'\alpha}
\]
for some $A_\gamma^{\alpha'\alpha}\in \GL_n(\Z )$, $c_\gamma^{\alpha'\alpha}\in \R^n$.
Since $\phi_\alpha$ coincides with $\phi_\beta$ on each overlap $U_{\alpha\beta}$, $\phi_{\alpha'}\circ \phi_\gamma \circ \phi_\alpha(\tilde x)=\smash{A_\gamma^{\alpha'\alpha}\tilde x+c_\gamma^{\alpha'\alpha}}$ also agrees with \smash{$\phi_{\beta'}\circ \phi_\gamma \circ \phi_\beta(\tilde x)=A_\gamma^{\beta'\beta}\tilde x+c_\gamma^{\beta'\beta}$} on the overlap $\phi_\alpha(U_{\alpha\beta})=\phi_\beta(U_{\alpha\beta})$. This implies \smash{$A_\gamma^{\alpha'\alpha}$}'s and \smash{$c_\gamma^{\alpha'\alpha}$}'s do not depend on $\alpha$ and depends only on $\gamma$. In fact, for each $\gamma\in \Gamma$ and $\alpha_0\in A$, we set
\[
A_0:=\big\{\alpha\in A \mid A_\gamma^{\alpha'_0\alpha_0}=A_\gamma^{\alpha'\alpha}\text{ and }c_\gamma^{\alpha'_0\alpha_0}=c_\gamma^{\alpha'\alpha} \big\} .
\]
$A_0$ contains all $\beta \in A$ with $U_{\alpha_0\beta}\not=\varnothing$. In particular, $A_0$ is not empty since $\alpha_0\in A_0$. Then, we have
\[
\biggl(\bigcup_{\alpha\in A_0}U_\alpha\biggr)\cup \biggl(\bigcup_{\alpha\in A\setminus A_0}U_\alpha\biggr) =\tilde B,\qquad \biggl(\bigcup_{\alpha\in A_0}U_\alpha\biggr)\cap \biggl(\bigcup_{\alpha\in A\setminus A_0}U_\alpha\biggr) =\varnothing .
\]
If the complement $A\setminus A_0$ is not empty, this contradicts to the connectedness of $\tilde B$. So we denote them by $A_\gamma$ and $c_\gamma$, respectively. Thus, we define the homomorphism $\rho\colon \Gamma \to \GL_n(\Z)\ltimes \R^n$ by
\[
\rho_\gamma:=(A_\gamma, c_\gamma).
\]
$\Gamma$ acts on $\R^n$ by $\rho_\gamma(x)=A_\gamma x+c_\gamma$ for $\gamma\in \Gamma$ and $x\in \R^n$. The lift $\tilde \rho_\gamma$ of $\rho_\gamma$ to $(M',\omega')$ is defined~by
\[
\tilde \rho_\gamma([x_\alpha, y_\alpha]):=\big[\varphi_{\alpha'}\circ \tilde \sigma_\gamma \circ \varphi_\alpha^{-1}(x_\alpha,y_\alpha)\big]
\]
if $(x_\alpha,y_\alpha)$ is in $\phi_\alpha(U_\alpha)\times T^n$. By construction, $\tilde\rho$ is a lift of $\rho$, and $\tilde \rho$ and $\rho$ satisfy $\smash{\tilde \dev (\tilde \sigma_\gamma (\tilde p))}=\smash{\tilde \rho_\gamma \bigl(\tilde \dev (\tilde p)\bigr)}$ and \smash{$\dev \bigl(\sigma_\gamma \bigl(\tilde b\bigr)\bigr)=\rho_\gamma \bigl(\dev \bigl(\tilde b\bigr)\bigr)$}, respectively.
\end{proof}

\begin{Remark}\label{remarks on M'}\quad
\begin{itemize}\itemsep=0pt
\item[(1)] By construction, the $n$-dimensional torus $T^n$ acts freely on $M'$ preserving $\omega'$ from the right so that $\pi'\colon M'\to O$ is a principal $T^n$-bundle.
\item[(2)] If $\pi\colon (M,\omega)\to B$ admits a Lagrangian section, the restriction of $\pi_0\colon {(\R^n\times T^n,\omega_0)\to \R^n}$ to $O$ can be taken as $\pi'\colon (M',\omega')\to O$. In fact, in this case, since $\pi\colon (M,\omega)\to B$ is identified with the canonical model, we can take a system of local trivializations $\big\{\bigl(\pi^{-1}(U_\alpha),\varphi_\alpha\bigr)\big\}$ with $u_{\alpha\beta}=0$ on each overlaps $U_{\alpha\beta}$. By applying the construction of $\pi'\colon (M',\omega')\to O$ given in the proof of Proposition~\ref{Lagrangian fibration on simply connected base} to such a $\big\{\bigl(\pi^{-1}(U_\alpha),\varphi_\alpha\bigr)\big\}$, we can show the claim.
\end{itemize}
\end{Remark}

Suppose that $(M,\omega)$ is prequantizable and let $\bigl(L,\nabla^L\bigr)\to (M,\omega)$ be a prequantum line bundle. We denote by \smash{$\bigl(\tilde L,\nabla^{\tilde L}\bigr)\to \bigl(\tilde M,\tilde \omega\bigr)$} the pullback of $\bigl(L,\nabla^L\bigr)\to (M,\omega)$ to \smash{$\bigl(\tilde M,\tilde \omega\bigr)$}. By definition, \smash{$\tilde L$}~admits a natural lift of the $\Gamma$-action $\tilde \sigma$ on \smash{$\bigl(\tilde M,\tilde\omega\bigr)$} which preserves \smash{$\nabla^{\tilde L}$}. The $\Gamma$-action on \smash{$\bigl(\tilde L, \nabla^{\tilde L}\bigr)$} is denoted by $\vardbtilde \sigma$. Then, we have the following prequantum version of Proposition~\ref{Lagrangian fibration on simply connected base}.
\begin{Proposition}
There exists a prequantum line bundle \smash{$\bigl(L',\nabla^{L'}\bigr)\to (M',\omega')$}, a bundle immersion \smash{$\vardbtilde \dev\colon \bigl(\tilde L,\nabla^{\tilde L}\bigr)\to \bigl(L',\nabla^{L'}\bigr)$} which covers \smash{$\tilde \dev$}, and a lift \smash{$\vardbtilde \rho$} of the $\Gamma$-action $\tilde \rho$ on $(M',\omega')$ to \smash{$\bigl(L',\nabla^{L'}\bigr)$} such that \smash{$\vardbtilde \dev$} is equivariant with respect to $\vardbtilde \sigma$ and $\vardbtilde \rho$.
\end{Proposition}
\begin{proof}
Let $\{(U_\alpha,\phi_\alpha)\}$ and $\big\{ \bigl(\tilde \pi^{-1}(U_\alpha),\omega|_{\tilde \pi^{-1}(U_\alpha)}, \varphi_\alpha\bigr) \big\}$ be the integral affine atlas of $\tilde B$ and the local trivializations of \smash{$\tilde \pi\colon \bigl(\tilde M,\tilde \omega\bigr)\to \tilde B$} obtained in the proof of Proposition~\ref{Lagrangian fibration on simply connected base}, respectively. Then, for each $\alpha\in A$ there exists a prequantum line bundle
\[
\bigl(\phi_\alpha(U_\alpha)\times T^n\times \C ,\nabla^{\tilde L_\alpha}\bigr)\to (\phi_\alpha(U_\alpha)\times T^n,\omega_0)
\]
 and a bundle isomorphism \smash{$\psi_\alpha\colon \bigl(\tilde L,\nabla^{\tilde L}\bigr)|_{\tilde \pi^{-1}(U_\alpha)}\!\to\! \bigl(\phi_\alpha(U_\alpha)\times T^n\times \C ,\nabla^{\tilde L_\alpha}\bigr)$} which covers $\varphi$. Now we define $\bigl(L',\nabla^{L'}\bigr)$ by
\[
\bigl(L',\nabla^{L'}\bigr):=\coprod_{\alpha \in A}\bigl(\phi_\alpha(U_\alpha)\times T^n\times \C ,\nabla^{\tilde L_\alpha}\bigr)/{\sim} ,
\]
where $(x_\alpha,y_\alpha,z_\alpha)\in \phi_\alpha(U_\alpha)\times T^n\times \C$ and $(x_\beta,y_\beta,z_\beta)\in \phi_\beta(U_\beta)\times T^n\times \C$ are in the relation $(x_\alpha,y_\alpha,z_\alpha)\sim (x_\beta,y_\beta,z_\beta)$ if they satisfy $(x_\beta,y_\beta,z_\beta)=\psi_\alpha\circ \psi_\beta^{-1}(x_\beta,y_\beta,z_\beta)$.
$\smash{\vardbtilde \dev\colon \bigl(\tilde L,\nabla^{\tilde L}\bigr)}\to \smash{\bigl(L',\nabla^{L'}\bigr)}$ is defined by
\[
\vardbtilde \dev (\tilde v):=[\psi_\alpha (\tilde v)]
\]
if $\tilde v$ is in \smash{$\bigl(\tilde L,\nabla^{\tilde L}\bigr)|_{\tilde \pi^{-1}(U_\alpha)}$}.

Suppose that for each $\gamma\in \Gamma$ the deck transformation $\sigma_\gamma$ maps each $U_\alpha$ to some $U_{\alpha'}$ as before. Then, \smash{$\vardbtilde \sigma_\gamma$} maps \smash{$\tilde L_{\tilde \pi^{-1}(U_\alpha)}$} to \smash{$\tilde L_{\tilde \pi^{-1}(U_{\alpha'})}$}. Then, the $\Gamma$-action $\vardbtilde \rho$ is defined by
\[
\vardbtilde \rho_\gamma (x_\alpha,y_\alpha, z_\alpha):=\big[\psi_{\alpha'}\circ \vardbtilde \sigma_\gamma \circ \psi_\alpha^{-1}(x_\alpha,y_\alpha, z_\alpha)\big]
\]
if $(x_\alpha,y_\alpha, z_\alpha)$ is in $\phi_\alpha(U_\alpha)\times T^n\times \C$.
\end{proof}

In the case where $B$ is complete, we obtain the following corollary.
\begin{Corollary}\label{complete base}
Let $\pi\colon (M,\omega)\to B$ be a Lagrangian fibration with connected $n$-dimensional base $B$ and $\bigl(L,\nabla^L\bigr)\to (M,\omega)$ a prequantum line bundle on $(M,\omega)$. Let $p\colon \tilde B\to B$ be the universal covering of $B$. Let us denote by \smash{$\bigl(\tilde M,\tilde \omega\bigr)$} the pullback of $(M,\omega)$ to \smash{$\tilde B$} and denote by~\smash{$\bigl(\tilde L,\nabla^{\tilde L}\bigr)$} the pullback of $\bigl(L,\nabla^L\bigr)$ to \smash{$\bigl(\tilde M,\tilde \omega\bigr)$}.
If $B$ is complete, there exist an integral affine isomorphism \smash{$\dev\colon \tilde B\to\R^n$}, a fiber-preserving symplectomorphism \smash{$\tilde \dev\colon \bigl(\tilde M,\tilde \omega\bigr)\to(\R^n\times T^n,\omega_0)$}, and a bundle isomorphism \smash{$\vardbtilde \dev\colon \bigl(\tilde L,\nabla^{\tilde L}\bigr)\to \bigl(\R^n\times T^n\times \C , d-2\pi\sqrt{-1}x\cdot {\rm d}y\bigr)$} such that \smash{$\tilde \dev$} covers $\dev$ and \smash{$\vardbtilde \dev$} covers \smash{$\tilde \dev$}, respectively. Here $x\cdot {\rm d}y$ denotes \smash{$\sum_{i=1}^nx_i{\rm d}y_i$}.
Moreover, let $\sigma$ be the $\Gamma$-action on $\tilde B$ defined as the inverse of deck transformations, $\tilde \sigma$ the natural lift of $\sigma$ to \smash{$\bigl(\tilde M,\tilde \omega\bigr)$}, and \smash{$\vardbtilde \sigma$} the natural lift of $\tilde \sigma$ to \smash{$\bigl(\tilde L,\nabla^{\tilde L}\bigr)$}, respectively. Then, there exist an integral affine $\Gamma$-action~${\rho\colon \Gamma \to \GL_n(\Z)\ltimes \R^n}$ on $\R^n$, its lifts $\tilde \rho$ and \smash{$\vardbtilde \rho$} to $(\R^n\times T^n,\omega_0)$ and $\bigl(\R^n\times T^n\times \C , d-2\pi\sqrt{-1}x\cdot {\rm d}y\bigr)$, respectively such that $\dev$, \smash{$\tilde \dev$}, and \smash{$\vardbtilde \dev$} are $\Gamma$-equivariant.
\end{Corollary}

\begin{proof}
By construction of \smash{$\tilde \dev$} given in the proof of Proposition~\ref{Lagrangian fibration on simply connected base}, if $\dev$ is bijective, so is \smash{$\tilde \dev$}. The argument in \cite[p.\ 696]{Du} and Theorem~\ref{classification} also show that $\pi_0\colon (\R^n\times T^n,\omega_0)\to \R^n$ is the unique Lagrangian fibration on $\R^n$ up to fiber-preserving symplectomorphism covering the identity. In particular, $\pi'\colon (M',\omega')\to \R^n$ is identified with $\pi_0\colon (\R^n\times T^n,\omega_0)\to \R^n$.

Concerning the prequantum line bundle, it is sufficient to show that $(\R^n\times T^n,\omega_0)$ has a~unique prequantum line bundle $\bigl(\R^n\times T^n\times \C , d-2\pi\sqrt{-1}x\cdot {\rm d}y\bigr)$ up to bundle isomorphism. Since~$\omega_0$ is exact, any prequantum line bundle on $(\R^n\times T^n,\omega_0)$ is trivial as a complex line bundle. Let~${\bigl(\R^n\times T^n\times \C ,d-2\pi\sqrt{-1}\alpha\bigr)}$ be a prequantum line bundle on $(\R^n\times T^n,\omega_0)$ with connection~${d-2\pi\sqrt{-1}\alpha}$. Then, $\alpha-x\cdot {\rm d}y$ defines a de Rham cohomology class in $H^1(\R^n\times T^n; \R)$. Since $H^1(\R^n\times T^n; \R)$ is isomorphic to $H^1(T^n;\R)$, in terms of the generators ${\rm d}y_i$'s of $H^1(T^n;\R)$, $\alpha-x\cdot {\rm d}y$ can be described as
\[
\alpha-x\cdot {\rm d}y=\sum_{i=1}^n\tau_i{\rm d}y_i+{\rm d}f
\]
for some $\tau_1,\ldots ,\tau_n\in \R$ and $f\in C^\infty(\R^n\times T^n)$. Now we define the bundle isomorphism $\psi\colon \R^n\times T^n\times \C \to \R^n\times T^n\times \C$ by
\[
\psi(x,y,z):=\bigl(x+(\tau_i),y,{\rm e}^{-2\pi\sqrt{-1}f(x,y)}z\bigr) .
\]
Then, $\psi$ satisfies $\psi^*\bigl(d-2\pi\sqrt{-1}x\cdot {\rm d}y\bigr)=d-2\pi\sqrt{-1}\alpha$.
\end{proof}

In the rest of this paper, we use the notations $(M_0, \omega_0):=\bigl( \R^n\times T^n, \sum_{i=1}^n{\rm d}x_i\wedge {\rm d}y_i\bigr)$ and~${\bigl( L_0, \nabla^{L_0}\bigr):=\bigl(\R^n\times T^n\times \C , d-2\pi\sqrt{-1}x\cdot {\rm d}y\bigr)}$ for simplicity.

\begin{Remark}[Hermitian metric on $\bigl( L_0, \nabla^{L_0}\bigr)$]\label{Hermitian metric}
By Corollary~\ref{complete base}, any Lagrangian fibration~${\pi\colon (M,\omega)\to B}$ on a complete $B$ with prequantum line bundle $\bigl(L,\nabla^L\bigr)\to (M,\omega)$ is obtained as the quotient space of the $\Gamma$-action on $\pi_0\colon (M_0,\omega_0)\to \R^n$ with prequantum line bundle $\bigl( L_0, \nabla^{L_0}\bigr)
\to (M_0,\omega_0)$. By definition, the prequantum line bundle $\bigl(L,\nabla^L\bigr)\to (M,\omega)$ is equipped with a Hermitian metric $ \< \cdot , \, \cdot \>_L$ compatible with $\nabla^L$.\footnote{A Hermitian metric $ \< \cdot , \, \cdot \>_L$ on $L$ is compatible with $\nabla^L$ if it satisfies $d (\< s_1, s_2\>_L)=\<\nabla^L s_1, s_2\>_L +\< s_1, \nabla^L s_2\>_L$ for all $s_1, s_2\in \Gamma ( L)$. } The pull-back of $ \< \cdot , \, \cdot \>_L$ to~${\bigl( L_0, \nabla^{L_0}\bigr)\to (M_0,\omega_0)}$
coincides with the one induced from the standard Hermitian inner product on $\C$ up to constant. In fact, it is easy to see that, up to constant, it is the unique Hermitian metric on $\bigl( L_0, \nabla^{L_0}\bigr)\to (M_0,\omega_0)$
compatible with $\nabla^{L_0}$.
In the rest of this paper, we assume that $\bigl( L_0, \nabla^{L_0}\bigr)\to (M_0,\omega_0)$
is always equipped with the Hermitian metric though we do not specify it.
\end{Remark}

\subsection[The lifting problem of fiber-preserving symplectomorphisms to the prequantum line bundle]{The lifting problem of fiber-preserving symplectomorphisms \\ to the prequantum line bundle}

Let $\Gamma'$ be a group,
and suppose that $\Gamma'$ acts on $\pi_0\colon (M_0,\omega_0)\to \R^n$ as fiber-preserving symplectomorphisms. As in the previous section, we denote by $\rho\colon \Gamma' \to \GL_n(\Z)\ltimes \R^n$ the $\Gamma'$-action on $\R^n$ and also denote by $\tilde \rho$ its lift to $(M_0,\omega_0)$.
By Lemma~\ref{fiber-preserving symplectomorphism}, for each $\gamma\in \Gamma'$, there exist~${A_\gamma \in\GL_n(\Z)}$, $c_\gamma\in \R^n$, and a map $u_\gamma\colon \R^n\to T^n$ with ${}^tA_\gamma \J u_\gamma$ symmetric such that $\rho_\gamma$ and~$\tilde \rho_\gamma$ can be described as follows
\begin{equation}\label{Gamma-action on Rn times Tn}
\rho_\gamma (x)=A_\gamma x+c_\gamma ,\qquad {\tilde \rho}_\gamma (x,y)=\bigl(A_\gamma x+c_\gamma ,{}^tA_\gamma^{-1}y+u_\gamma (x)\bigr) .
\end{equation}
Note that since \eqref{Gamma-action on Rn times Tn} is a $\Gamma'$-action, $A_\gamma$, $c_\gamma$, and $u_\gamma$ satisfy the following conditions:
\begin{gather}
A_{\gamma_1\gamma_2}=A_{\gamma_1}A_{\gamma_2},\qquad
c_{\gamma_1\gamma_2}=A_{\gamma_1}c_{\gamma_2}+c_{\gamma_1},\qquad
u_{\gamma_1\gamma_2}(x)={}^tA_{\gamma_1}^{-1}u_{\gamma_2}(x)+u_{\gamma_1}(\rho_{\gamma_2}(x))\label{condition for action}
\end{gather}
for $\gamma_1$, $\gamma_2\in \Gamma'$, and $x\in \R^n$.
Let ${\tilde u}_\gamma={}^t\bigl({\tilde u}_\gamma^1,\ldots ,{\tilde u}_\gamma^n\bigr)\colon \R^n\to \R^n$ be a lift of $u_\gamma$. For ${\tilde u}_\gamma$ and $i=1,\ldots ,n$, we put
\[
\int_0^{x_i}{\tilde u}_\gamma (x){\rm d}x_i:=
{}^t\left(\int_0^{x_i}{\tilde u}_\gamma^1 (x){\rm d}x_i, \ldots , \int_0^{x_i}{\tilde u}_\gamma^n (x){\rm d}x_i\right)
\]
and
\[
F_\gamma^i(x):=\left({}^tA_\gamma \int_0^{x_i}{\tilde u}_\gamma (x){\rm d}x_i\right)_i=\sum_{j=1}^n\bigl({}^tA_\gamma\bigr)_{ij}\int_0^{x_i}{\tilde u}_\gamma^j (x){\rm d}x_i.
\]
Let $N\in \N$ be a positive integer. Each $\tilde \rho_\gamma$ preserves $N\omega_0$, hence, $\Gamma'$ also acts on $\pi_0\colon (M_0,N\omega_0)\to \R^n$ as fiber-preserving symplectomorphisms. Then, we examine in detail the conditions for the $\Gamma'$-action to have a lift to \smash{$\bigl( L_0, \nabla^{L_0}\bigr)^{\otimes N}\to (M_0,N\omega_0)$}. The purpose of this subsection is to show the following lemma which gives the necessary and sufficient condition on the existence of a lift of the $\Gamma'$-action, and which also gives the explicit formula for the lift when this condition is satisfied.

\begin{Lemma}\label{condition for lift of Gamma-action to tilde_L}\quad
\begin{itemize}\itemsep=0pt
\item[$(1)$] For each $\gamma\in \Gamma'$, there exists a bundle automorphism \smash{$\vardbtilde \rho_\gamma$} of \smash{$\bigl( L_0, \nabla^{L_0}\bigr)^{\otimes N}$} preserving the Hermitian metric and the connection such that \smash{$\vardbtilde \rho_\gamma$} covers $\tilde \rho_\gamma$
if and only if $c_\gamma$ is contained in $\frac{1}{N}\Z^n$. Moreover, in this case, \smash{$\vardbtilde \rho_\gamma$} can be described as follows
\begin{equation}\label{vardbtilde rho}
{\vardbtilde \rho}_\gamma (x,y,z)=\bigl( {\tilde \rho}_\gamma (x,y),\ g_\gamma {\rm e}^{2\pi\sqrt{-1}N\{{\tilde g}_\gamma(x)+c_\gamma\cdot ({}^tA_\gamma^{-1}y)\}}z\bigr)
\end{equation}
 for $(x,y,z)\in L_0^{\otimes N}\cong \R^n\times T^n\times \C$, where $g_\gamma$ is an arbitrary element in ${\rm U}(1)$ and
\[
{\tilde g}_\gamma(x):=\rho_\gamma (x)\cdot {\tilde u}_\gamma(x)-c_\gamma\cdot {\tilde u}_\gamma (0)-\sum_{i=1}^nF_\gamma^i(0,\ldots ,0,x_i,\ldots ,x_n).
\]
The formula \eqref{vardbtilde rho} does not depend on the choice of ${\tilde u}_\gamma$.\footnote{In the rest of this paper, we often use the notation $u_\gamma$ instead of ${\tilde u}_\gamma$.}
\item[$(2)$] Under the condition given in $(1)$, the map \smash{$\vardbtilde \rho\colon \Gamma' \to \Aut \bigl(\bigl( L_0, \nabla^{L_0}\bigr)^{\otimes N}\bigr)$} defined by \eqref{vardbtilde rho} is a homomorphism if and only if
the map $g\colon \Gamma'\ni \gamma\mapsto g_\gamma \in {\rm U}(1)$ is a homomorphism and for all $\gamma_1, \gamma_2 \in \Gamma'$ and $x\in \R^n$, the following condition holds:
\begin{gather*}
\bigl\{-c_{\gamma_1}\cdot u_{\gamma_1}(0)+c_{\gamma_1}\cdot {}^tA_{\gamma_1}^{-1}u_{\gamma_2}(0)+\rho_{\gamma_1}(c_{\gamma_2})\cdot u_{\gamma_1}(\rho_{\gamma_2}(0))\bigr\}\\
\qquad-\sum_{i=1}^n\left({}^tA_{\gamma_1}\int_0^{(\rho_{\gamma_2}(x))_i}u_{\gamma_1}(0,\ldots ,0,\tau_i,(\rho_{\gamma_2}(x))_{i+1},\ldots ,(\rho_{\gamma_2}(x))_n){\rm d}\tau_i\right)_i \\
\qquad+\sum_{i=1}^n\left({}^tA_{\gamma_2}{}^tA_{\gamma_1}\int_0^{x_i}u_{\gamma_1}(\rho_{\gamma_2}(0,\ldots 0,\tau_i,x_{i+1},\ldots ,x_n)){\rm d}\tau_i\right)_i \in\frac{1}{N}\Z .
\end{gather*}
\end{itemize}
\end{Lemma}
\begin{proof}
For each $\gamma\in \Gamma'$, we put
\[
{\vardbtilde \rho}_\gamma (x,y,z)=\bigl( {\tilde \rho}_\gamma (x,y),{\rm e}^{2\pi\{{\tilde g}^R_\gamma (x,y)+\sqrt{-1}{\tilde g}^I_\gamma(x,y)\}}z\bigr) ,
\]
where \smash{${\tilde g}^R_\gamma$} and \smash{${\tilde g}^I_\gamma$} are real valued functions on $M_0$. By the direct computation, it is easy to see
that \smash{${\vardbtilde \rho}_\gamma$} preserves \smash{$\nabla^{L_0^{\otimes N}}=d-2\pi\sqrt{-1}Nx\cdot {\rm d}y$} if and only if \smash{${\tilde g}^R_\gamma$} is constant and \smash{${\tilde g}^I_\gamma$} satisfies the following conditions:
\begin{gather}
\partial_{x_i}{\tilde g}^I_\gamma=N(A_\gamma x+c_\gamma)\cdot \partial_{x_i} {\tilde u}_\gamma,\label{partial x}\\
\partial_{y_i} {\tilde g}^I_\gamma=N\bigl(A_\gamma^{-1}c_\gamma\bigr)_i \label{partial y}
\end{gather}
for $i=1, \ldots , n$.
The conditions for the complete integrability of the system of partial differential equations \eqref{partial x} and \eqref{partial y} are as follows:
\begin{align}
\partial_{x_i}\partial_{x_j}{\tilde g}^I_\gamma&=\partial_{x_j}\partial_{x_i}{\tilde g}^I_\gamma , \label{xx}\\
\partial_{x_i}\partial_{y_j}{\tilde g}^I_\gamma&=\partial_{y_j}\partial_{x_i}{\tilde g}^I_\gamma , \label{xy}\\
\partial_{y_i}\partial_{y_j}{\tilde g}^I_\gamma&=\partial_{y_j}\partial_{y_i}{\tilde g}^I_\gamma \label{yy}
\end{align}
for $i, j=1, \ldots , n$. From \eqref{partial x} and \eqref{partial y}, \eqref{xy} and \eqref{yy} are true because both sides of each vanish. From \eqref{partial x}, \eqref{xx} can be expressed as
\begin{equation}\label{xx2}
\bigl( {}^tA_\gamma\partial_{x_i}u_\gamma (x)\bigr)_j=\bigl( {}^tA_\gamma\partial_{x_j}u_\gamma (x)\bigr)_i
\end{equation}
for $i, j=1, \ldots , n$. But, since ${}^tA_\gamma \J u_\gamma$ is symmetric, \eqref{xx2} is also valid. Therefore, we know that there exists ${\tilde g}^I_\gamma$ that satisfies \eqref{partial x} and \eqref{partial y}. In fact, such a ${\tilde g}^I_\gamma$ is given by
 \begin{align}
 {\tilde g}^I_\gamma (x,y)={}&{\tilde g}^I_\gamma (0,0)+N\biggl\{\rho_\gamma (x)\cdot {\tilde u}_\gamma(x)-c_\gamma \cdot {\tilde u}_\gamma(0) \nonumber\\
 & -\sum_{i=1}^nF_\gamma^i(0,\ldots ,0,x_i,\ldots ,x_n)
 +c_\gamma \cdot {}^tA_\gamma^{-1} y\biggr\}.\label{g(x,y)}
 \end{align}
 Since $y\in T^n$, ${\tilde g}^I_\gamma$ should satisfies \smash{${\rm e}^{2\pi\sqrt{-1}{\tilde g}^I_\gamma (0,e_i)}={\rm e}^{2\pi\sqrt{-1}{\tilde g}^I_\gamma (0,0)}$} for all $i=1,\ldots ,n$ and $\gamma\in \Gamma'$. This holds if and only if $A_\gamma^{-1}Nc_\gamma \cdot e_i\in \Z$ for all $i=1,\ldots ,n$ and $\gamma\in \Gamma'$. Since $A_\gamma \in \GL_n(\Z)$, this is equivalent to the condition $Nc_\gamma \in \Z^n$. In this case, we put \smash{$g_\gamma:={\rm e}^{2\pi({\tilde g}^R_\gamma (0,0)+\sqrt{-1}{\tilde g}^I_\gamma (0,0))}$}. Since \smash{${\vardbtilde \rho}_\gamma$} preserves the Hermitian metric on $\bigl( L_0, \nabla^{L_0}\bigr)\to (M_0,\omega_0)$,
$g_\gamma$ is contained in ${\rm U}(1)$. The formula \eqref{vardbtilde rho} does not depend on the choice of ${\tilde u}_\gamma$ since the difference of two lifts of $u_\gamma$ is in $\Z^n$. This proves $(1)$.

The map \smash{${\vardbtilde \rho}$} defined in $(2)$ is a homomorphism if and only if \smash{${\tilde g}^I_\gamma (x,y)-{\tilde g}^I_\gamma (0,0)$} defined by~\eqref{g(x,y)} satisfies the cocycle condition. By a direct computation using \eqref{condition for action}, it is equivalent to the ones given in $(2)$.
 \end{proof}

\begin{Example}\label{B times T}
Let $B$ be the $n$-dimensional integral affine torus given in Example~\ref{affine torus}\,(1) for a~linear basis $v_1,\ldots ,v_n\in \R^n$. The product $B\times T^n$ admits a symplectic structure $\omega$ so that the trivial torus bundle $\pi\colon (B\times T^n,\omega)\to B$ becomes a Lagrangian fibration. This is obtained as the quotient space of the action of $\Gamma':=\Z^n$ on $\pi_0\colon (M_0,\omega_0)\to \R^n$ which is defined by
\[
{\tilde \rho}_\gamma (x,y)=(x+C\gamma, y)
\]
for $\gamma\in \Gamma'$ and $(x,y)\in M_0$, where $C=(v_1\cdots v_n)\in \GL_n(\R)$.
Let $N\in \N$ be a positive number. The $\Gamma'$-action ${\tilde \rho}$ on $(M_0, N\omega_0)$ has a lift to the prequantum line bundle $\smash{\bigl( L_0, \nabla^{L_0}\bigr)^{\otimes N}}
\to (M_0,N\omega_0)$ if and only if all $v_i$'s lie in $\frac{1}{N}\Z^n$, and in this case ${\vardbtilde \rho}$ is given by
\[
{\vardbtilde \rho}_\gamma (x,y,z)=\bigl({\tilde \rho}_\gamma (x,y), g_\gamma {\rm e}^{2\pi\sqrt{-1}NC\gamma\cdot y}z\bigr)
\]
for $\gamma\in \Gamma'$ and $(x,y,z)\in L_0^{\otimes N}\cong \R^n\times T^n\times \C$, where $g\colon \Gamma' \ni \gamma\mapsto g_\gamma\in {\rm U}(1)$ is an arbitrary homomorphism.
\end{Example}

\begin{Example}[the Kodaira--Thurston manifold]\label{Kodaira--Thurston}
Let $\Gamma'$ be $\Z^2$. Let us consider the $\Gamma'$-action on $\pi_0\colon \bigl(\R^2\times T^2,\omega_0\bigr)\to \R^2$ which is defined by
\[
\rho_\gamma(x):=x+\gamma,\qquad \tilde \rho_\gamma (x,y):=(\rho_\gamma(x), y+u_\gamma (x))
\]
for $\gamma \in \Gamma'$ and $(x,y)\in \R^2\times T^2$, where $u_\gamma (x)={}^t(0,\gamma_1x_2)$. The Lagrangian fibration given by the quotient of this action is denoted by $\pi\colon (M,\omega)\to B$. $M$ was first observed by Kodaira in~\cite{Kodaira} and Thurston pointed out in~\cite{Thurston} that $(M,\omega)$ does not admits any K\"ahler structure. $M$ is nowadays called the Kodaira--Thurston manifold.
Let $N\in \N$ be a positive number. The $\Gamma'$-action~${\tilde \rho}$ on $\bigl(\R^2\times T^2, N\omega_0\bigr)$ has a lift to the prequantum line bundle $\bigl(\R^2\times T^2\times \C, d-2\pi\sqrt{-1}Nx\cdot {\rm d}y\bigr)\to \bigl(\R^2\times T^2,N\omega_0\bigr)$ if and only if $N$ is even, and in this case the lift ${\vardbtilde \rho}$ is given by
\[
{\vardbtilde \rho}_\gamma (x,y,z)=\bigl({\tilde \rho}_\gamma (x,y), g_\gamma {\rm e}^{2\pi\sqrt{-1}N\{\frac{1}{2}\gamma_1x_2^2+\gamma_1\gamma_2x_2+\gamma\cdot y\}}z\bigr)
\]
for $\gamma\in \Gamma'$ and $(x,y,z)\in \R^2\times T^2\times \C$, where $g\colon \Gamma' \ni \gamma\mapsto g_\gamma\in {\rm U}(1)$ is an arbitrary homomorphism.
\end{Example}

\begin{Example}\label{Ex4.8}
Let $B$ be the $n$-dimensional integral affine torus given in Example~\ref{affine torus}\,(1) for a linear basis $v_1,\ldots ,v_n\in \R^n$. When all $v_i$'s are integer vectors, i.e., $v_1,\ldots ,v_n\in \Z^n$, we can generalize Examples~\ref{B times T} and~\ref{Kodaira--Thurston} in the following way. For $i, j=1,\ldots , n$, choose $u_{ij}\in \Z^n$ satisfying $u_{ij}=u_{ji}$.
For each $\gamma \in \Gamma':=\Z^n$, we define the map $u_\gamma\colon \R^n\to T^n$ by
\[
u_\gamma (x):=
\begin{pmatrix}
u_{11}\cdot \gamma & \cdots & u_{1n}\cdot \gamma \\
\vdots & & \vdots \\
u_{n1}\cdot \gamma & \cdots & u_{nn}\cdot \gamma
\end{pmatrix}
x ,
\]
and we also define the action of $\Gamma'$ on $\pi_0\colon (M_0,\omega_0)\to \R^n$ by
\begin{equation}\label{Gamma-action 1}
{\tilde \rho}_\gamma (x,y)=\left(x+C\gamma,\ y+u_\gamma (x)\right)
\end{equation}
for $\gamma\in \Gamma'$ and $(x,y)\in M_0$, where $C=(v_1\ \cdots\ v_n)$. Then, the quotient $\pi\colon (M,\omega)\to B$ obtained as the $\Gamma'$-action~\eqref{Gamma-action 1} is a Lagrangian fibration on $B$.
Let $N\in \N$ be a positive number. The $\Gamma'$-action ${\tilde \rho}$ on $(M_0, N\omega_0)$ has a lift to the prequantum line bundle \smash{$\bigl( L_0, \nabla^{L_0}\bigr)^{\otimes N}
\to (M_0,N\omega_0)$} if and only if $\frac{N}{2}v_i\cdot U_jv_i\in \Z$ for all $i, j=1,\ldots , n$, where
\[
U_j:=
\begin{pmatrix}
(u_{11})_j & \cdots & (u_{1n})_j \\
\vdots & & \vdots \\
(u_{n1})_j & \cdots & (u_{nn})_j
\end{pmatrix}.
\]
And in this case, the lift \smash{${\vardbtilde \rho}$} is given by
\[
{\vardbtilde \rho}_\gamma (x,y,z)=\bigl({\tilde \rho}_\gamma (x,y), g_\gamma {\rm e}^{2\pi\sqrt{-1}N[\frac{1}{2}\{\rho_\gamma (x)\cdot u_\gamma (\rho_\gamma (x))-\rho_\gamma(0)\cdot u_\gamma (\rho_\gamma (0))\}+\rho_\gamma(0)\cdot y]}z\bigr)
\]
for $\gamma\in \Gamma'$ and \smash{$(x,y,z)\in L_0^{\otimes N}\cong \R^n\times T^n\times \C$}, where $g\colon \Gamma' \ni \gamma\mapsto g_\gamma\in {\rm U}(1)$ is an arbitrary homomorphism.
\end{Example}

\begin{Example}\label{Jordan block2}
Let $n\ge 2$ and $\lambda_1,\ldots, \lambda_{n-1}\in \Z$. Let $\Gamma'$ be the group $(\Z^n,\circ )$ given in Example~\ref{Jordan block}. For each $\gamma\in \Gamma'$, let $A_\gamma$ be the matrix
\[
A_\gamma:=
\begin{pmatrix}
1 & \lambda_1 & & & \\
 & 1 & \lambda_2 & & \\
 & & \ddots & \ddots & \\
 & & & 1 & \lambda_{n-1}\\
 & & & & 1
\end{pmatrix}^{\gamma_n}
\]
and $u_\gamma\colon \R^n\to T^n$ the map defined by
\[
u_\gamma (x):=
\begin{pmatrix}
0 \\
\vdots \\
0 \\
\gamma_nx_n
\end{pmatrix}.
\]
Let us consider the $\Gamma'$-action ${\tilde\rho}$ on $\pi_0\colon ( M_0,\omega)\to \R^n$ which is defined by
\begin{equation}\label{Gamma-action 2}
{\tilde \rho}_\gamma (x,y):=\bigl( A_\gamma x+\gamma , {}^tA_\gamma^{-1}y+u_\gamma (x)\bigr)
\end{equation}
for $\gamma\in \Gamma'$ and $(x,y)\in M_0$. Then, the quotient $\pi\colon (M,\omega)\to B$ obtained as the $\Gamma'$-action~\eqref{Gamma-action 2} is a Lagrangian fibration on the integral affine manifold $B$ obtained in Example~\ref{Jordan block}.
Let ${N\in \N}$ be a positive number. The $\Gamma'$-action ${\tilde \rho}$ on $(M_0, N\omega_0)$ has a lift to the prequantum line bundle~\smash{$\bigl( L_0, \nabla^{L_0}\bigr)^{\otimes N}
\to (M_0,N\omega_0)$} if and only if $N$ is even, and in this case the lift ${\vardbtilde \rho}$ is given by
\[
{\vardbtilde \rho}_\gamma (x,y,z)=\bigl({\tilde \rho}_\gamma (x,y), g_\gamma {\rm e}^{2\pi\sqrt{-1}N\{ \gamma_n x_n(\frac{1}{2}x_n+\gamma_n)+\gamma\cdot ({}^tA_\gamma^{-1}y)\}}z\bigr)
\]
for $\gamma\in \Gamma'$ and \smash{$(x,y,z)\in L_0^{\otimes N}\cong \R^n\times T^n\times \C$}, where $g\colon \Gamma' \ni \gamma\mapsto g_\gamma\in {\rm U}(1)$ is an arbitrary homomorphism.
\end{Example}

\section[Degree-zero harmonic spinors and integrability of almost complex structures]{Degree-zero harmonic spinors and integrability\\ of almost complex structures}\label{equiv quantization}
Let $N\in \N$ be a positive integer. For a compatible almost complex structure $J$ on the total space of the Lagrangian fibration $\pi_0\colon (M_0,N\omega_0)\to \R^n$, let $D$ be the associated Spin${}^c$ Dirac operator with coefficients in the prequantum line bundle \smash{$\bigl( L_0, \nabla^{L_0}\bigr)^{\otimes N}
\to (M_0,N\omega_0)$}. An element in the kernel $\ker D$ of $D$ is called a harmonic spinor. In this section, for $J$ which is invariant along the fiber in the sense of Lemma~\ref{invariant J}, we investigate the condition on the existence of nontrivial degree-zero harmonic spinors, i.e., nontrivial sections which is contained in $\ker D$. For the construction and properties of the Spin${}^c$ Dirac operator, see \cite{MR1365745, LM}.

\subsection{Bohr--Sommerfeld points}
Let $\pi\colon (M,\omega)\to B$ be a Lagrangian fibration with prequantum line bundle $\bigl(L,\nabla^L\bigr)\to (M,\omega)$. We recall the definition of Bohr--Sommerfeld points.
\begin{Definition}
A point $b\in B$ is said to be \emph{Bohr--Sommerfeld} if $\bigl(L,\nabla^L\bigr)|_{\pi^{-1}(b)}$ admits a~nontrivial covariant constant section. We denote the set of Bohr--Sommerfeld points by $B_{\rm BS}$.
\end{Definition}

Let us detect Bohr--Sommerfeld points for $\pi_0\colon (M_0,N\omega_0)\to \R^n$ with prequantum line bundle~\smash{$\bigl( L_0, \nabla^{L_0}\bigr)^{\otimes N}\to (M_0,N\omega_0)$}.
\begin{Proposition}\label{BS}
A point $x\in \R^n$ is Bohr--Sommerfeld if and only if $x$ is contained in $\frac{1}{N}\Z^n$, i.e., $\R^n_{\rm BS}=\frac{1}{N}\Z^n$. Moreover, for a Bohr--Sommerfeld point $x\in \frac{1}{N}\Z^n$, a covariant constant section $s$ of \smash{$\bigl( L_0, \nabla^{L_0}\bigr)^{\otimes N}\big|_{\pi_0^{-1}(x)}$} is of the form \smash{$s(y)=s(0){\rm e}^{2\pi\sqrt{-1}Nx\cdot y}$}.
\end{Proposition}
\begin{proof}
For a fixed $x\in \R^n$, a section $s$ of \smash{$\bigl( L_0, \nabla^{L_0}\bigr)^{\otimes N}\big|_{\pi_0^{-1}(x)}\to \pi_0^{-1}(x)$} is covariant constant if and only if $s$ satisfies
\[
0=\nabla^{L_0^{\otimes N}}_{\partial_{y_i}}s=\partial_{y_i}s-2\pi\sqrt{-1}Nx_is
\]
for $i=1,\ldots , n$. Hence, any covariant constant section $s$ should be of the form $s(y)=\smash{s(0){\rm e}^{2\pi\sqrt{-1}Nx\cdot y}}$. Since $\pi_0^{-1}(x)$ is a torus, $s$ is periodic with respect to $y_i$'s. In particular, $s$ satisfies~\smash{$s(0)=s(e_i)=s(0){\rm e}^{2\pi\sqrt{-1}Nx_i}$} for $i=1,\ldots , n$. This implies that \smash{$\bigl( L_0, \nabla^{L_0}\bigr)^{\otimes N}\!\big|_{\pi_0^{-1}(x)}\!\!\to \pi_0^{-1}(x)$} admits a nontrivial covariant constant section if and only if $Nx_i\in \Z$ for $i=1,\ldots , n$.
\end{proof}

\begin{Remark}\label{action preserves BS}
Suppose that $\pi_0\colon (M_0,N\omega_0)\to \R^n$ is equipped with an action of a group $\Gamma$ which preserves all the data, and its lift \smash{$\vardbtilde \rho$} to \smash{$\bigl( L_0, \nabla^{L_0}\bigr)^{\otimes N}$} is given by \eqref{vardbtilde rho}. Then, by Lemma~\ref{condition for lift of Gamma-action to tilde_L}\,(1), the $\Gamma$-action $\rho$ on $\R^n$ preserves $\R^n_{\rm BS}$. When the $\Gamma$-action $\rho$ on $\R^n$ is properly discontinuous and free, let $F\subset \R^n$ be a fundamental domain of the $\Gamma$-action $\rho$ on $\R^n$. Then, the map
\begin{equation}\label{Gamma times F}
\Gamma\times \left(F\cap \frac{1}{N}\Z^n\right)\ni \left(\gamma ,\dfrac{m}{N}\right)\mapsto N\rho_\gamma \left(\dfrac{m}{N}\right)\in \Z^n
\end{equation}
can be defined and is bijective. In particular, if a Lagrangian fibration $\pi\colon (M,N\omega)\to B$ with prequantum line bundle \smash{$\bigl(L,\nabla^L\bigr)^{\otimes N}\to (M,N\omega)$} is obtained as the quotient space of the $\Gamma$-action, then $F\cap \frac{1}{N}\Z^n$ is identified with $B_{\rm BS}$.
\end{Remark}

\subsection{Almost complex structures}
Let ${\mathcal S}_n$ be the Siegel upper half space, namely, the space of $n\times n$ symmetric complex matrices whose imaginary parts are positive definite
\[
{\mathcal S}_n:=\big\{ Z=X+\sqrt{-1}Y\in M_n(\C)\mid X, Y\in M_n(\R), {}^tZ=Z, \text{and }Y\text{ is positive definite} \big\} .
\]
It is well known that ${\mathcal S}_n$ is identified with the space of compatible complex structures on the~$2n$-dimensional standard symplectic vector space. See \cite[Chapter II, Section~4]{Mumford1}.

For a tangent vector $u=\sum_{i=1}^n \{(u_x)_i\partial_{x_i}+(u_y)_i\partial_{y_i} \}\in T_{(x,y)}M_0$ at a point $(x,y)\in M_0$, we use the following notation:
\[
u
= (\partial_{x_1},\ldots, \partial_{x_n}, \partial_{y_1},\ldots, \partial_{y_n} )
\begin{pmatrix}
(u_x)_1 \\
\vdots \\
(u_x)_n \\
(u_y)_1\\
\vdots \\
(u_y)_n
\end{pmatrix}
=(\partial_x, \partial_y)
\begin{pmatrix}
u_x \\
u_y
\end{pmatrix},
\]
where
\[
\partial_x=(\partial_{x_1},\ldots, \partial_{x_n}),\qquad \partial_y=(\partial_{y_1},\ldots, \partial_{y_n}),\qquad
u_x=
\begin{pmatrix}
(u_x)_1\\
\vdots \\
(u_x)_n
\end{pmatrix},\qquad
u_y=
\begin{pmatrix}
(u_y)_1\\
\vdots\\
(u_y)_n
\end{pmatrix}.
\]
In terms of the notations of tangent vectors
$u=(\partial_x, \partial_y)
\left(\begin{smallmatrix}
u_x \\
u_y
\end{smallmatrix}\right)$
and
$v=(\partial_x, \partial_y)
\left(\begin{smallmatrix}
v_x \\
v_y
\end{smallmatrix}\right)\in T_{(x,y)}M_0$, $\omega_0$ can be described by
\begin{equation*}\label{omega}
\omega_0(u,v)=\left({}^tu_x, {}^tu_y\right)
\begin{pmatrix}
0 & I \\
-I & 0
\end{pmatrix}
\begin{pmatrix}
v_x\\
v_y
\end{pmatrix}.
\end{equation*}
Since the tangent bundle $TM_0$ is trivial, the space of compatible almost complex structures on $(M_0,\omega_0)$ is identified with the space of $C^\infty$ maps from $M_0$ to ${\mathcal S}_n$. For $Z=X+\sqrt{-1}Y\in C^\infty ( M_0,{\mathcal S}_n)$, the corresponding almost complex structure $J_Z$ is given as follows:
\begin{equation}\label{J}
J_Zu:=(\partial_x, \partial_y)
\begin{pmatrix}
XY^{-1} & -Y-XY^{-1}X \\
Y^{-1} & -Y^{-1}X
\end{pmatrix}_{(x,y)}
\begin{pmatrix}
u_x \\
u_y
\end{pmatrix}
\end{equation}
for $u=(\partial_x, \partial_y)
\left(\begin{smallmatrix}
u_x \\
u_y
\end{smallmatrix}\right)\in T_{(x,y)}M_0$.\footnote{
$\left(\begin{smallmatrix}
XY^{-1} & -Y-XY^{-1}X \\
Y^{-1} & -Y^{-1}X
\end{smallmatrix}\right)_{(x,y)}$, $\bigl(XY^{-1}\bigr)_{(x,y)}$ etc.\
are the values of the maps
$\left(\begin{smallmatrix}
XY^{-1} & -Y-XY^{-1}X \\
Y^{-1} & -Y^{-1}X
\end{smallmatrix}\right)$, $XY^{-1}$ etc.\
at $(x,y)$. We will often omit the subscript \lq\lq${}_{(x,y)}$" for simplicity unless it causes confusion.} Then, the Riemannian metric $g$ determined by $\omega_0$ and $J_Z$ can be described by
\begin{align}
g(u, v):&=\omega_0(u,Jv)=\bigl({}^tu_x, {}^tu_y\bigr)
\begin{pmatrix}
0 & I \\
-I & 0
\end{pmatrix}
\begin{pmatrix}
XY^{-1} & -Y-XY^{-1}X \\
Y^{-1} & -Y^{-1}X
\end{pmatrix}
\begin{pmatrix}
v_x\\
v_y
\end{pmatrix}\nonumber\\
&=\bigl({}^tu_x, {}^tu_y\bigr)
\begin{pmatrix}
Y^{-1} & -Y^{-1}X \\
-XY^{-1} & Y+XY^{-1}X
\end{pmatrix}
\begin{pmatrix}
v_x\\
v_y
\end{pmatrix}.\label{g}
\end{align}

Let $J=J_Z$ be the almost complex structure on $(M_0,\omega_0)$ corresponding to a given $Z=X+\sqrt{-1}Y\in C^\infty ( M_0,{\mathcal S}_n)$. Then, $(-J\partial_y, \partial_y)= (-J\partial_{y_1},\ldots,-J\partial_{y_n}, \partial_{y_1},\ldots ,\partial_{y_n} )$ is also a basis of the tangent space of $(M_0,\omega_0)$. For each tangent vector $u\in T_{(x,y)}M_0$, by using this basis, $u$ is expressed as follows:
\[
u=\sum_i\{ (u_H)_i(-J\partial_{y_i})+(u_V)_i\partial_{y_i}\}=(-J\partial_y, \partial_y)
\begin{pmatrix}
u_H\\
u_V
\end{pmatrix}.
\]
Then, we have the following transition formula between $ (\partial_x, \partial_y )$ and $(-J\partial_y, \partial_y)$:
\begin{equation*}
u=(-J\partial_y, \partial_y)
\begin{pmatrix}
u_H\\
u_V
\end{pmatrix}
=(\partial_x, \partial_y)\left(
\begin{pmatrix}
-XY^{-1} & Y+XY^{-1}X \\
-Y^{-1} & Y^{-1}X
\end{pmatrix}
\begin{pmatrix}
0\\
u_H
\end{pmatrix}
+
\begin{pmatrix}
0\\
u_V
\end{pmatrix}
\right) .
\end{equation*}
By this formula, we obtain the following lemma.
\begin{Lemma}\label{g2}
In terms of this notation, the Riemannian metric $g$ defined by \eqref{g} can be described by
\begin{align*}
g(u,v)={}&\bigl(0,{}^tu_H\bigr)
\begin{pmatrix}
Y^{-1} & -Y^{-1}X \\
-XY^{-1} & Y+XY^{-1}X
\end{pmatrix}
\begin{pmatrix}
0\\
v_H
\end{pmatrix}\\
&
+\bigl(0,{}^tu_V\bigr)
\begin{pmatrix}
Y^{-1} & -Y^{-1}X \\
-XY^{-1} & Y+XY^{-1}X
\end{pmatrix}
\begin{pmatrix}
0\\
v_V
\end{pmatrix}.
\end{align*}
\end{Lemma}

Suppose that a group $\Gamma$ acts on $\pi_0\colon (M_0,\omega_0)\to \R^n$ and the $\Gamma$-actions $\rho$ on $\R^n$ and $\tilde \rho$ on~$(M_0,\omega_0)$ are written as in \eqref{Gamma-action on Rn times Tn}. Then, it is easy to see the following lemma.
\begin{Lemma}
The $\Gamma$-action ${\tilde \rho}$ on $(M_0,\omega_0)$ preserves the almost complex structure $J=J_Z$ on~$(M_0,\omega_0)$ corresponding to $Z=X+\sqrt{-1}Y\in C^\infty ( M_0,{\mathcal S}_n)$ if and only if the following conditions hold:
\begin{align}
& A_\gamma\bigl(XY^{-1}\bigr)_{(x,y)}=\bigl(XY^{-1}\bigr)_{{\tilde \rho}_\gamma (x,y)}A_\gamma-\bigl(Y+XY^{-1}X\bigr)_{{\tilde \rho}_\gamma (x,y)}(\J u_\gamma)_x , \label{Gamma-invariance=J1}\\
& A_\gamma\bigl(Y+XY^{-1}X\bigr)_{(x,y)}=\bigl(Y+XY^{-1}X\bigr)_{{\tilde \rho}_\gamma (x,y)} {}^tA_\gamma^{-1},\label{Gamma-invariance=J2} \\
& (\J u_\gamma)_x\bigl(XY^{-1}\bigr)_{(x,y)}+{}^tA_\gamma^{-1}Y^{-1}_{(x,y)}=Y^{-1}_{{\tilde \rho}_\gamma (x,y)}A_\gamma -\bigl(Y^{-1}X\bigr)_{{\tilde \rho}_\gamma (x,y)}(\J u_\gamma)_x.\nonumber
\end{align}
\end{Lemma}
\begin{proof}
For all $\gamma \in \Gamma$ and $(x,y)\in (M_0,\omega_0)$, the condition
\[
({\rm d}{\tilde \rho}_\gamma)_{(x,y)}\circ J_{(x,y)}=J_{{\tilde \rho}_\gamma (x,y)}\circ ({\rm d}{\tilde \rho}_\gamma)_{(x,y)}
\]
 implies above three equalities together with the following equality:
\[
(\J u_\gamma)_x\bigl(Y+XY^{-1}X\bigr)_{(x,y)}+{}^tA_\gamma^{-1}\bigl(Y^{-1}X\bigr)_{(x,y)}=\bigl(Y^{-1}X\bigr)_{{\tilde \rho}_\gamma (x,y)}{}^tA_\gamma^{-1}.
\]
But, this can be obtained from \eqref{Gamma-invariance=J1}, \eqref{Gamma-invariance=J2}, and ${}^t\bigl({}^tA_\gamma (\J u_\gamma)_x\bigr)={}^tA_\gamma (\J u_\gamma)_x$.
\end{proof}

Let $\pi\colon (M,\omega)\to B$ be a Lagrangian fibration with complete $n$-dimensional base $B$ and $p\colon \tilde B\to B$ the universal covering of $B$. By Corollary~\ref{complete base}, the pullback of $\pi\colon (M,\omega)\to B$ to~$\tilde B$ is identified with $\pi_0\colon (M_0,\omega_0)\to \R^n$ and $\pi\colon (M,\omega)\to B$ can be obtained as the quotient of the $\Gamma=\pi_1(B)$-action on $\pi_0\colon (M_0,\omega_0)\to \R^n$. In particular, for each compatible almost complex structure $J$ on $(M,\omega)$, there exists a map $Z_J=X+\sqrt{-1}Y\in C^\infty\left(M_0,{\mathcal S}_n\right)$ such that the pullback $p^*J$ of $J$ to $p^*(M,\omega)$ coincides with $J_{Z_J}$. Then, we have the following lemma.
\begin{Lemma}[\textup{\cite[Corollary~9.15]{FFY2}}]\label{invariant J}
For any Lagrangian fibration $\pi\colon (M,\omega)\to B$, there exists a~compatible almost complex structure $J$ on $(M,\omega)$ such that the corresponding map $Z_J$ does not depend on $y_1,\ldots ,y_n$. We say such $J$ to be invariant along the fiber.
\end{Lemma}
\begin{proof}
Take a Riemannian metric $g'$ on $(M,\omega)$. Then, the pullback $p^*g'$ is $\pi_1(B)$-invariant. Moreover, $p^*(M,\omega)$ admits a free $T^n$-action, and this $T^n$-action together with the $\pi_1(B)$-action forms an action of the semi-direct product $\pi_1(B)\ltimes T^n$ of $T^n$ and $\pi_1(B)$. By averaging $p^*g'$ over~$T^n$, we obtain a Riemannian metric on $p^*M$ invariant under the $\pi_1(B)\ltimes T^n$-action. It is easy to see that $p^*\omega$ is also $\pi_1(B)\ltimes T^n$-invariant, so by the standard method using the $\pi_1(B)\ltimes T^n$-invariant Riemannian metric and $p^*\omega$, we can obtain a $\pi_1(B)\ltimes T^n$-invariant compatible almost complex structure on $p^*(M,\omega)$. In particular, since the almost complex structure is still invariant under $\pi_1(B)$-action, it descends to $(M,\omega)$. This is the required almost complex structure.
\end{proof}

\subsection[A condition on the existence of nontrivial harmonic spinors of degree-zero]{A condition on the existence of nontrivial harmonic spinors\\ of degree-zero}
\label{Spin-c Dirac operator}
For a map $Z=X+\sqrt{-1}Y\in C^\infty ( M_0,{\mathcal S}_n)$, we set
\begin{equation}\label{Omega}
\Omega:=\bigl(Y+XY^{-1}X\bigr)^{-1}ZY^{-1}.
\end{equation}

\begin{Lemma}\label{property of Omega} $\Omega$ has the following properties:
\begin{itemize}\itemsep=0pt
\item[$(1)$] $\Omega={\overline Z}^{-1}$, where ${\overline Z}=X-\sqrt{-1}Y$.
\item[$(2)$] $\Omega$ is symmetric, i.e., ${}^t\Omega=\Omega$.
\end{itemize}
\end{Lemma}
\begin{proof}
A direct computation shows that $\Omega{\overline Z}=I$. This proves (1). (2) follows from (1) since~$Z$ is symmetric.
\end{proof}

Let $N\in \N$ be a positive integer. Let $J=J_Z$ be the compatible almost complex structure on $(M_0,N\omega_0)$ corresponding to a given $Z=X+\sqrt{-1}Y\in C^\infty ( M_0,{\mathcal S}_n)$. Then, the Riemannian metric $Ng:=N\omega_0 (\cdot , J\cdot )$ defines an isomorphism $f\colon T^*M_0\cong TM_0$ by $\tau=Ng\left( f(\tau),\ \cdot \right)$ for~${\tau \in T^*M_0}$. For $i=1,\ldots ,n$, let $\Omega_i$ denote the $i$th column vector of $\Omega$, and $\re\Omega_i$ and $\im\Omega_i$ be the real and imaginary parts of $\Omega_i$, respectively. Then, we can show the following lemma.
\begin{Lemma}\label{dual}
For $i=1,\ldots ,n$,
\[
f({\rm d}x_i)=-\frac{1}{N}J\partial_{y_i},\qquad f({\rm d}y_i)=(-J\partial_y, \partial_y)
\begin{pmatrix}
\dfrac{1}{N}\re\Omega_i \vspace{1mm}\\
\dfrac{1}{N}\im\Omega_i
\end{pmatrix}.
\]
\end{Lemma}
\begin{proof}
We prove the latter. The former can be proved by the same way.
Put $f({\rm d}y_i)=\smash{(-J\partial_y, \partial_y)
\left(\begin{smallmatrix}
Y^i_H \\
Y^i_V
\end{smallmatrix}\right)}$. By definition, for each $i, j=1,\ldots ,n$, we have
\begin{gather}
{\rm d}y_i(-J\partial_{y_j})=Ng\left((-J\partial_y, \partial_y)
\begin{pmatrix}
Y^i_H \\
Y^i_V
\end{pmatrix},
(-J\partial_y, \partial_y)
\begin{pmatrix}
e_j \\
0
\end{pmatrix}
\right), \label{-J partial_y}\\
{\rm d}y_i(\partial_{y_j})=Ng\left((-J\partial_y, \partial_y)
\begin{pmatrix}
Y^i_H \\
Y^i_V
\end{pmatrix},
(-J\partial_y, \partial_y)
\begin{pmatrix}
0 \\
e_j
\end{pmatrix}
\right). \label{partial_y}
\end{gather}
Since $-J\partial_{y_j}$ is written as
\[
-J\partial_{y_j}=(\partial_x, \partial_y)
\begin{pmatrix}
-XY^{-1} & Y+XY^{-1}X \\
-Y^{-1} & Y^{-1}X
\end{pmatrix}
\begin{pmatrix}
0 \\
e_j
\end{pmatrix}
\]
by \eqref{J}, the left-hand side of \eqref{-J partial_y} is \smash{$\bigl(Y^{-1}X\bigr)_{ij}$}. On the other hand, by Lemma~\ref{g2}, the right-hand side of \eqref{-J partial_y} can be described as $NY^i_H\cdot \bigl(Y+XY^{-1}X\bigr)e_j$. This implies $Y^{-1}X=N{}^t\bigl(Y^1_H\cdots Y^n_H\bigr)\bigl(Y+XY^{-1}X\bigr)$. Since $Y$ is positive definite, so is $Y+XY^{-1}X$. In particular, $N\bigl(Y+XY^{-1}X\bigr)$ is invertible. By using ${}^tX=X$, ${}^tY=Y$ together with this fact, we can obtain \smash{$\bigl(Y^1_H\cdots Y^n_H\bigr)=\frac{1}{N}\bigl(Y+XY^{-1}X\bigr)^{-1}XY^{-1}$}. By the same way, from \eqref{partial_y}, we obtain $I=N{}^t\bigl(Y^1_V\cdots Y^n_V\bigr)\bigl(Y+XY^{-1}X\bigr)$, i.e., \smash{$\bigl(Y^1_V\cdots Y^n_V\bigr)=\frac{1}{N}\bigl(Y+XY^{-1}X\bigr)^{-1}$}. Hence, $\frac{1}{N}\Omega=\bigl(Y^1_H\cdots Y^n_H\bigr)+\sqrt{-1}\bigl(Y^1_V\cdots Y^n_V\bigr)$.
\end{proof}

Define the Hermitian metric on $(M_0, N\omega_0, Ng, J)$ by
\begin{equation}\label{Hermitian metric on M}
h(u,v):=Ng(u,v)+\sqrt{-1}Ng(u,Jv)
\end{equation}
for $u, v\in T_{(x,y)}M_0$. Let $(W,c)$ be the Clifford module bundle associated with $(Ng,J)$, i.e., as a~complex vector bundle, $W$ is defined by
\[
W:=\wedge^\bullet (TM_0, J)\otimes_{\C}\bigl(L_0^{\otimes N}\bigr).
\]
$W$ is equipped with the Hermitian metric induced from $h$ and that on $L_0$, and also equipped with the Hermitian connection, which is denoted by $\nabla^W$, induced from the Levi-Civita connection~$\nabla^{LC}$ of $(M_0,Ng)$ and $\nabla^{L_0}$. $c$ is the Clifford multiplication $c\colon TM_0\to \End_{\C}(W)$ which is defined by
$
c(u)(\tau):=u\wedge \tau -u \llcorner_h \tau
$
for $u\in TM_0$ and $\tau \in W$, where $\llcorner_h $ is the contraction with respect to the Hermitian metric~$h$ on~${(M_0, N\omega_0, Ng, J)}$. It is well known that $W$ is identified with $\wedge^\bullet(T^*M_0)^{0,1}\otimes_{\C}\bigl({L_0}^{\otimes N}\bigr)$ as a Clifford module bundle since $h$ induces the isomorphism from $(TM_0, J)$ to $(T^*M_0)^{0,1}$ as Hermitian vector bundles. See \cite[pp.\ 12--13]{MR1365745} for more details.

Now let us define the Spin${}^c$ Dirac operator $D \colon \Gamma (W)\to \Gamma (W)$ by the composition of the following maps:
\[
\xymatrix{
D\colon\ \Gamma(W)\ar[r]^{\nabla^W}& \Gamma (T^*M_0\otimes W)\ar[r]^{f\otimes \id_W} & \Gamma (TM_0\otimes W)\ar[r]^{c} & \Gamma (W).
}
\]
We compute the action of $D$ on degree zero elements in $\Gamma (W)$. We identify a section of ${L_0}$ with a complex valued function on $M_0$. By using Lemma~\ref{dual}, for a section $s$ of ${L_0}^{\otimes N}$, $Ds$ can be computed as
\begin{align*}
Ds&=c\circ (f\otimes \id_W)\circ \nabla^Ws=c\circ (f\otimes \id_W)\bigl({\rm d}s-2\pi\sqrt{-1}Nx\cdot {\rm d}ys\bigr)\\
&=\sum_{i=1}^n\bigl\{c(f({\rm d}x_i))(\partial_{x_i}s)+c(f({\rm d}y_i))\bigl(\partial_{y_i}s-2\pi\sqrt{-1}Nx_is\bigr)\bigr\}\\
&=-\frac{\sqrt{-1}}{N}\sum_{i=1}^n\partial_{y_i}\otimes_{\C}\biggl\{\partial_{x_i}s+\sum_{j=1}^n\Omega_{ij}\bigl(\partial_{y_j}s-2\pi\sqrt{-1}Nx_js\bigr)\biggr\}.
\end{align*}
In particular, the equality $Ds=0$ is equivalent to
\begin{equation}\label{Dirac2}
0=
\begin{pmatrix}
\partial_{x_1}s\\
\vdots \\
\partial_{x_n}s
\end{pmatrix}
+\Omega
\begin{pmatrix}
\partial_{y_1}s-2\pi\sqrt{-1}Nx_1s\\
\vdots \\
\partial_{y_n}s-2\pi\sqrt{-1}Nx_ns
\end{pmatrix}.
\end{equation}

For a section $s$ of ${L_0}^{\otimes N}$, let us consider the Fourier series expansion of $s$ with respect to $y_i$'s. For each $x\in \R^n$, as a function of $y_i$'s, $s(x, \cdot )$ can be expressed as the Fourier series
\begin{equation}\label{Fourier}
s(x,y)=\sum_{m\in \Z^n}a_m(x){\rm e}^{2\pi\sqrt{-1}m\cdot y},
\end{equation}
where
$
a_m(x):=\int_{T^n}s(x,y){\rm e}^{-2\pi\sqrt{-1}m\cdot y}{\rm d}y
$
for $m\in \Z^n$. Suppose that $Z$ does not depend on $y_1, \ldots , y_n$ as in Lemma~\ref{invariant J}. Then, by using the Fourier series \eqref{Fourier}, the equation $Ds=0$ can be reduced to the following system of differential equations for $a_m$'s with variables $x_1,\ldots ,x_n$.
\begin{Lemma}
$s$ satisfies $Ds=0$ if and only if $a_m$'s satisfy
\begin{equation}\label{Dirac3}
0=
\begin{pmatrix}
\partial_{x_1}a_m\\
\vdots \\
\partial_{x_n}a_m
\end{pmatrix}
+2\pi\sqrt{-1}a_m\Omega
(m-Nx)
\end{equation}
for all $m\in \Z^n$.
\end{Lemma}
\begin{proof}
By Lemma~\ref{Fourier for derivatives}, the partial derivatives $\partial_{x_j}s$ and $\partial_{y_j}s$ have the following Fourier series with respect to $y_i$'s:
\begin{align}
\partial_{x_j}s(x,y)&=\sum_{m\in \Z^n}\partial_{x_j}a_m(x){\rm e}^{2\pi\sqrt{-1}m\cdot y},\label{Fourier for sx} \\
\partial_{y_j}s(x,y)&=\sum_{m\in \Z^n}2\pi\sqrt{-1}m_ja_m(x){\rm e}^{2\pi\sqrt{-1}m\cdot y} \qquad\text{for $j=1,\ldots , n$.}\label{Fourier for sy}
\end{align}
Then, substituting \eqref{Fourier}, \eqref{Fourier for sx} and \eqref{Fourier for sy} into \eqref{Dirac2}, we can obtain~\eqref{Dirac3}.
\end{proof}

We investigate the equation \eqref{Dirac3}.
\begin{Lemma}
Let $a_m$ be a solution of \eqref{Dirac3} for some $m\in \Z^n$. If there exists $p\in \R^n$ such that~${a_m(p)=0}$, then $a_m(x)=0$ for all $x\in \R^n$.
\end{Lemma}
\begin{proof}
First, fix the variables $x_2, \ldots, x_n$ to equal $p_2, \ldots ,p_n$, respectively. Then, the first entry of \eqref{Dirac3}, i.e., $0=\partial_{x_1}a_m+2\pi\sqrt{-1}a_m(\Omega(m-Nx))_1$ can be thought of as an ordinary differential equation on $x_1$, and $a_m(x_1, p_2,\ldots ,p_n)$ is its solution with initial condition $a_m(p)=0$. On the other hand, the trivial solution also has the same initial condition. By the uniqueness of the solution of the ordinary differential equation, $a_m(x_1, p_2,\ldots ,p_n)=0$ for any $x_1$. Next, by fixing variables $x_3, \ldots, x_n$ with $p_3, \ldots ,p_n$ and fixing $x_1$ with arbitrary value, $a_m(x_1, x_2, p_3,\ldots ,p_n)$ is a~solution of $0=\partial_{x_2}a_m+2\pi\sqrt{-1}a_m(\Omega(m-Nx))_2$ with initial condition $a_m(x_1, p_2,\ldots ,p_n)=0$. Then, $a_m(x_1, x_2,p_3, \ldots ,p_n)=0$ for any $x_1$, $x_2$. By repeating the process for $x_3, \ldots ,x_n$, we can show that $a_m(x)=0$.
\end{proof}

\begin{Lemma}\label{basis of kerD}
If $a_m$ is a nontrivial smooth solution of \eqref{Dirac3} for some $m\in \Z^n$, then the condition
\begin{equation}\label{commutativity}
\left(\left(\partial_{x_i}\Omega\right)_x(m-Nx)\right)_j=\left(\left(\partial_{x_j}\Omega\right)_x(m-Nx)\right)_i
\end{equation}
holds  for all $i, j=1,\ldots , n$ and $x\in \R^n$.
Conversely, if there exists $m\in \Z^n$ such that \eqref{commutativity} holds for all $i, j=1,\ldots , n$ and $x\in \R^n$,
then \eqref{Dirac3} has a unique nontrivial solution up to constant. Moreover, in this case, each solution $a_m$ of \eqref{Dirac3} has the following form:
\begin{equation}\label{a(x)}
a_m(x)=a_m\left(\frac{m}{N}\right){\rm e}^{-2\pi\sqrt{-1}\sum_{i=1}^nG_m^i(\frac{m_1}{N},\ldots ,\frac{m_{i-1}}{N},x_i,\ldots ,x_n)},
\end{equation}
where $a_m(\frac{m}{N})$ can be taken as an arbitrary constant in $\C$ and
$
G_m^i(x):=\bigl(\int_{\frac{m_i}{N}}^{x_i}\Omega(m-Nx){\rm d}x_i\bigr)_i $.
\end{Lemma}
\begin{proof}
Since $a_m$ is smooth, $a_m$ satisfies $\partial_{x_i}\partial_{x_j}a_m=\partial_{x_j}\partial_{x_i}a_m$ for all $i, j=1,\ldots ,n$. By differentiating \eqref{Dirac3}, we have
\begin{align*}
\partial_{x_i}\partial_{x_j}a_m={}&-2\pi\sqrt{-1}a_m\left\{-2\pi\sqrt{-1}\sum_{k=1}^n\Omega_{ik}(m_k-Nx_k)\sum_{l=1}^n\Omega_{jl}(m_l-Nx_l)\right.\\
&\left.
+\sum_{l=1}^n\left(\partial_{x_i}\Omega_{jl}\right)(m_l-Nx_l)-N\Omega_{ji}\right\}
\end{align*}
for $i, j=1, \ldots ,n$ and $x\in \R^n$. The condition \eqref{commutativity} is obtained from this equation.

Conversely, suppose there exists $m\in \Z^n$ such that \eqref{commutativity} holds for all $i, j=1,\ldots , n$ and $x\in \R^n$. By solving the differential equation appeared as the $i$th component of \eqref{Dirac3} for $i=1,\ldots ,n$, we have
\begin{equation}\label{a(xi-0)}
a_m(x)=a_m\left(x_1,\ldots ,x_{i-1},\frac{m_i}{N},x_{i+1},\ldots ,x_n\right){\rm e}^{-2\pi\sqrt{-1}G_m^i(x)}.
\end{equation}
Using \eqref{a(xi-0)} repeatedly, we obtain the formula~\eqref{a(x)}.
By using \eqref{commutativity}, we can show that \eqref{a(x)} does not depend on the order of applying \eqref{a(xi-0)} to $x_i$'s as in the proof of Lemma~\ref{condition for lift of Gamma-action to tilde_L}. Hence, \eqref{a(x)} is well defined.
\end{proof}

\begin{Definition}
We say $m\in \Z^n$ to be \emph{integrable} if~\eqref{commutativity} holds for all $i, j=1,\ldots , n$ and $x\in \R^n$.
\end{Definition}

For each $m\in \Z^n$ which is integrable, define the section $s_m\in \Gamma\bigl({L_0}^{\otimes N}\bigr)$ by
\begin{equation}\label{s_m}
s_m(x,y):={\rm e}^{2\pi\sqrt{-1}\{-\sum_{i=1}^nG_m^i(\frac{m_1}{N},\ldots ,\frac{m_{i-1}}{N},x_i,\ldots ,x_n)+m\cdot y\}} .
\end{equation}
By the elliptic regularity of $D$ and Lemma~\ref{basis of kerD}, we can obtain the following.
\begin{Proposition}\label{existence condition}
If
\[
s\!=\!\sum_{m\in \Z^n}a_m(x){\rm e}^{2\pi\sqrt{-1}m\cdot y}\in \Gamma\bigl({L_0}^{\otimes N}\bigr)
\]
 is a nontrivial solution of~${0\!=\!Ds}$, then all $m\in \Z^n$ with $a_m\not=0$ are integrable. Conversely, suppose that there exists $m\in \Z^n$ such that $m$ is integrable. Then, the section $s_m$ defined by~\eqref{s_m}
satisfies $0=Ds_m$.
\end{Proposition}
The following proposition gives a geometric interpretation of the condition \eqref{commutativity}.
\begin{Proposition}\label{prop commutativity2}
The following conditions are equivalent:
\begin{itemize}\itemsep=0pt
\item[$(1$)] All $m\in \Z^n$ are integrable.
\item[$(2)$] $\partial_{x_i}\Omega_{jk}=\partial_{x_j}\Omega_{ik}$ for all $i, j, k=1,\ldots , n$.
\item[$(3)$] $\nabla^{LC} J=0$, where $\nabla^{LC}$ is the Levi-Civita connection with respect to $g$.
\end{itemize}
\end{Proposition}
A proof of Proposition~\ref{prop commutativity2} is given in Appendix~\ref{appendixB}.

\begin{Remark}\label{holomorphic structure}
When one of (hence, all) the conditions in Proposition~\ref{prop commutativity2} holds, $\left(M_0, \omega_0, J, g\right)$ is a K\"ahler manifold and $J$ induces a holomorphic structure on ${L_0}$ such that $\nabla^{L_0}$ is the canonical connection.
\end{Remark}

\subsection[The Gamma-equivariant case]{The $\boldsymbol{\Gamma}$-equivariant case}
Suppose that $\pi_0\colon (M_0,N\omega_0,J)\to \R^n$ with prequantum line bundle $\bigl({L_0}, \nabla^{L_0}\bigr)^{\otimes N}{\hspace*{-1mm}\to} (M_0,N\omega_0,J)$ is equipped with an action of a group $\Gamma$ which preserves all the data, and the $\Gamma$-actions are described by \eqref{Gamma-action on Rn times Tn} and \eqref{vardbtilde rho} as before. We assume that the $\Gamma$-action $\rho$ on $\R^n$ is properly discontinuous and free.
Since the $\Gamma$-action preserves all the data, the Spin${}^c$ Dirac operator $D$ is $\Gamma$-equivariant. In particular, $\Gamma$ acts on $\Gamma \bigl({L_0}^{\otimes N}\bigr)\cap \ker D$.
\begin{Lemma}\label{Gamma-invariance-lemma}
Let $s$ be a section of ${L_0}^{\otimes N}$ with the Fourier series of the form \eqref{Fourier}.
Then, $s$ is $\Gamma$-equivariant, i.e., ${\vardbtilde \rho}_\gamma\circ s=s\circ {\tilde \rho}_\gamma$ for all $\gamma\in \Gamma$ if and only if $a_m$ satisfies the following condition:
\begin{gather}\label{Gamma-invariance}
a_{N\rho_\gamma(\frac{m}{N})}(\rho_\gamma (x))=g_\gamma a_m(x){\rm e}^{2\pi\sqrt{-1}N\{ {\tilde g}_\gamma(x)
-\rho_\gamma(\frac{m}{N})\cdot u_\gamma (x)\}}
\end{gather}
for all $\gamma \in \Gamma$, $m\in \Z^n$, and $x\in \R^n$. In particular, any $\Gamma$-equivariant section of ${L_0}^{\otimes N}$ can be written as follows:
\begin{align}
s(x,y)={}&\sum_{(\gamma ,\frac{m}{N})\in \Gamma \times (F\cap \frac{1}{N}\Z^n)}g_{\gamma}a_m(\rho_{\gamma^{-1}}(x))\nonumber\\
&\times{\rm e}^{2\pi\sqrt{-1}N\{{\tilde g}_\gamma (\rho_{\gamma^{-1}}(x))-\rho_{\gamma}(\frac{m}{N})\cdot u_{\gamma}(\rho_{\gamma^{-1}}(x))\}}{\rm e}^{2\pi\sqrt{-1}N\rho_{\gamma}(\frac{m}{N})\cdot y}.\label{Gamma-equiv section}
\end{align}
\end{Lemma}
\begin{proof}
By computing the both sides separately, we have
\begin{align*}
{\vardbtilde \rho}_\gamma\circ s (x,y)&=g_\gamma {\rm e}^{2\pi\sqrt{-1}N \{{\tilde g}_\gamma (x)+c_\gamma\cdot {}^tA_\gamma^{-1}y \}}\sum_{m\in \Z^n}a_m(x){\rm e}^{2\pi\sqrt{-1}m\cdot y} \\
&=g_\gamma\sum_{m\in \Z^n}a_m(x){\rm e}^{2\pi\sqrt{-1}N{\tilde g}_\gamma (x)}{\rm e}^{2\pi\sqrt{-1}N\rho_{\gamma} (\frac{m}{N} )\cdot {}^tA_\gamma^{-1}y},\\
s\circ {\tilde \rho}_\gamma(x,y)&=\sum_{l\in \Z^n}a_l (\rho_\gamma(x) ){\rm e}^{2\pi\sqrt{-1}l\cdot ({}^tA_\gamma^{-1}y+u_\gamma(x) )} \\
&=\sum_{m\in \Z^n}a_{N\rho_{\gamma} (\frac{m}{N} )} (\rho_\gamma(x) ){\rm e}^{2\pi\sqrt{-1}N\rho_{\gamma} (\frac{m}{N} )\cdot u_\gamma(x)}{\rm e}^{2\pi\sqrt{-1}N\rho_{\gamma} (\frac{m}{N} )\cdot {}^tA_\gamma^{-1}y}.
\end{align*}
Here, in the last equality, we replace $l$ by $N\rho_{\gamma}(\frac{m}{N})$. Note that the map
\[
\Z^n\ni m\mapsto N\rho_{\gamma}\left(\frac{m}{N}\right)\in \Z^n
\]
 is bijective. Then, \smash{${\vardbtilde \rho}_\gamma\circ s=s\circ {\tilde \rho}_\gamma$} for all $\gamma\in \Gamma$ implies
\[
g_\gamma a_m(x){\rm e}^{2\pi\sqrt{-1}N{\tilde g}_\gamma (x)}=a_{N\rho_{\gamma}(\frac{m}{N})}(\rho_\gamma(x)){\rm e}^{2\pi\sqrt{-1}N\rho_{\gamma}(\frac{m}{N})\cdot u_\gamma(x)}
\]
for all $m\in \Z^n$. In particular, by \eqref{Gamma times F} and~\eqref{Gamma-invariance}, $s$ can be rewritten as follows:
\begin{align*}
s(x,y)={}&\sum_{l\in \Z^n}a_l(x){\rm e}^{2\pi\sqrt{-1}l\cdot y} \stackrel{\eqref{Gamma times F}}{=}\sum_{(\gamma ,\frac{m}{N})\in \Gamma \times (F\cap \frac{1}{N}\Z^n)}a_{N\rho_{\gamma}(\frac{m}{N})}(x){\rm e}^{2\pi\sqrt{-1}N\rho_{\gamma}(\frac{m}{N})\cdot y} \\
\stackrel{\eqref{Gamma-invariance}}{=}{}&\sum_{ (\gamma ,\frac{m}{N} )\in \Gamma \times (F\cap \frac{1}{N}\Z^n )}g_{\gamma}a_m (\rho_{\gamma^{-1}}(x) )\\
&\times {\rm e}^{2\pi\sqrt{-1}N \{{\tilde g}_\gamma (\rho_{\gamma^{-1}}(x) )-\rho_{\gamma}(\frac{m}{N})\cdot u_{\gamma} (\rho_{\gamma^{-1}}(x) ) \}}{\rm e}^{2\pi\sqrt{-1}N\rho_{\gamma}(\frac{m}{N})\cdot y}.\tag*{\qed}
\end{align*} \renewcommand{\qed}{}
\end{proof}

In the $\Gamma$-equivariant case, the condition~\eqref{commutativity} has a symmetry in the following sense.
\begin{Lemma}\label{equiv commutativity}
There exists $m_0\in \Z^n$  with $\frac{m_0}{N}\in F$ such that $m_0$ is integrable if and only if for any~${\gamma\in \Gamma}$ $m=N\rho_\gamma (\frac{m_0}{N})$
is integrable. Moreover, let $a_{m_0}$ be a nontrivial solution of \eqref{Dirac3} for $m_0$. For each $\gamma \in \Gamma$, we define \smash{$a_{N\rho_\gamma (\frac{m_0}{N})}$} in such a way that it satisfies \eqref{Gamma-invariance}. Then, \smash{$a_{N\rho_\gamma (\frac{m_0}{N})}$} is a nontrivial solution of \eqref{Dirac3} for $m=N\rho_\gamma (\frac{m_0}{N})$.
\end{Lemma}
\begin{proof}
Suppose that there exists $m_0\in \Z^n$ with $\frac{m_0}{N}\in F$ such that $m_0$ is integrable. By Lemma~\ref{basis of kerD}, \eqref{Dirac3} for $m_0$ has a nontrivial solution $a_{m_0}$. Then, for each $\gamma \in \Gamma$, define \smash{$a_{N\rho_\gamma (\frac{m_0}{N})}$} by \eqref{Gamma-invariance}. By Lemma~\ref{basis of kerD} again, in order to show this lemma, it is sufficient to prove \smash{$a_{N\rho_\gamma (\frac{m_0}{N})}$} is a solution of \eqref{Dirac3} for $m=N\rho_\gamma (\frac{m_0}{N})$. Let us compute the Jacobi matrix of the both sides of~\eqref{Gamma-invariance}. The left-hand side is
\begin{align}
\J (a_{N\rho_\gamma (\frac{m_0}{N})}\circ \rho_\gamma )_x
&= (\J a_{N\rho_\gamma (\frac{m_0}{N})} )_{\rho_\gamma(x)} (\J \rho_\gamma )_x\nonumber \\
&= (\partial_{x_1}a_{N\rho_\gamma (\frac{m_0}{N})}, \ldots , \partial_{x_n}a_{N\rho_\gamma (\frac{m_0}{N})} )_{\rho_\gamma (x)}A_\gamma . \label{LHS}
\end{align}
The right-hand side is
\begin{gather}
\J \bigl(g_\gamma a_m(x){\rm e}^{2\pi\sqrt{-1}N\{ {\tilde g}_\gamma(x)-\rho_\gamma(\frac{m}{N})\cdot u_\gamma (x)\}}\bigr)_x \nonumber\\
\qquad=g_\gamma {\rm e}^{2\pi\sqrt{-1}N\{ {\tilde g}_\gamma(x)-\rho_\gamma(\frac{m}{N})\cdot u_\gamma (x)\}}(\J a_m)_x +g_\gamma a_m(x)\J \bigl( {\rm e}^{2\pi\sqrt{-1}N\{ {\tilde g}_\gamma(x)-\rho_\gamma(\frac{m}{N})\cdot u_\gamma (x)\}}\bigr) \nonumber \\
\qquad\stackrel{\eqref{Dirac3}}{=}-2\pi\sqrt{-1}g_\gamma {\rm e}^{2\pi\sqrt{-1}N\{ {\tilde g}_\gamma(x)-\rho_\gamma(\frac{m}{N})\cdot u_\gamma (x)\}}a_m(x) {}^t(\Omega_x(m-Nx))\nonumber \\
\phantom{\qquad=}{}+2\pi\sqrt{-1}Ng_\gamma a_m(x){\rm e}^{2\pi\sqrt{-1}N\{ {\tilde g}_\gamma(x)-\rho_\gamma(\frac{m}{N})\cdot u_\gamma (x)\}} \J \left( {\tilde g}_\gamma(x)-\rho_\gamma\left(\frac{m}{N}\right)\cdot u_\gamma (x)\right) \nonumber \\
\qquad\stackrel{\eqref{Gamma-invariance}}{=}-2\pi\sqrt{-1}a_{N\rho_\gamma (\frac{m_0}{N})}(\rho_\gamma(x)){}^t\left(\Omega_xA_\gamma^{-1}\left(N\rho_\gamma \left(\frac{m}{N}\right)-N\rho_\gamma (x)\right)\right)\nonumber \\
\phantom{\qquad=}{}+2\pi\sqrt{-1}Na_{N\rho_\gamma (\frac{m_0}{N})}(\rho_\gamma(x))\J \left( {\tilde g}_\gamma(x)-\rho_\gamma\left(\frac{m}{N}\right)\cdot u_\gamma (x)\right) .\label{RHS}
\end{gather}
For each $i=1,\ldots ,n$, the direct computation shows
\begin{gather*}
\partial_{x_i}\left( {\tilde g}_\gamma(x)-\rho_\gamma\left(\frac{m}{N}\right)\cdot u_\gamma (x)\right)\\
\qquad=(\partial_{x_i}u_\gamma)_x\cdot \left(\rho_\gamma (x)-\rho_\gamma \left(\frac{m}{N}\right)\right)+\left({}^tA_\gamma u_\gamma(x)\right)_i-\left({}^tA_\gamma u_\gamma(0,\ldots ,0,x_i,\ldots ,x_n)\right)_i\\
\phantom{\qquad=}{} -\sum_{j<i}\int_0^{x_j}\bigl({}^tA_\gamma \J u_\gamma\bigr)_{ji}(0,\ldots ,0,x_j,\ldots ,x_n){\rm d}x_j \\
\qquad=(\partial_{x_i}u_\gamma)_x\cdot \left(\rho_\gamma (x)-\rho_\gamma \left(\frac{m}{N}\right)\right)
+\sum_{j<i}\int_0^{x_j}\partial_{x_j}\left({}^tA_\gamma u_\gamma (0,\ldots ,0,x_j,\ldots ,x_n)\right)_i{\rm d}x_j \\
\phantom{\qquad=}{} -\sum_{j<i}\int_0^{x_j}\left({}^tA_\gamma \J u_\gamma\right)_{ji}(0,\ldots ,0,x_j,\ldots ,x_n){\rm d}x_j \\
\qquad=(\partial_{x_i}u_\gamma)_x\cdot \left(\rho_\gamma (x)-\rho_\gamma \left(\frac{m}{N}\right)\right)
+\sum_{j<i}\int_0^{x_j}\left({}^tA_\gamma \J u_\gamma \right)_{ij}(0,\ldots ,0,x_j,\ldots ,x_n){\rm d}x_j\\
\phantom{\qquad=}{} -\sum_{j<i}\int_0^{x_j}\left({}^tA_\gamma \J u_\gamma\right)_{ji}(0,\ldots ,0,x_j,\ldots ,x_n){\rm d}x_j \\
\qquad=-(\partial_{x_i}u_\gamma)_x\cdot \left(\rho_\gamma \left(\frac{m}{N}\right)-\rho_\gamma (x)\right) .
\end{gather*}
In the last equality, we used ${}^t\bigl({}^tA_\gamma \J u_\gamma\bigr)={}^tA_\gamma \J u_\gamma$. Hence, we have
\begin{equation}\label{RHS2}
\J \left( {\tilde g}_\gamma(x)-\rho_\gamma\left(\frac{m}{N}\right)\cdot u_\gamma (x)\right)=-{}^t\left(\rho_\gamma \left(\frac{m}{N}\right)-\rho_\gamma (x)\right)\left(\J u_\gamma\right)_x.
\end{equation}
By \eqref{LHS}, \eqref{RHS} and \eqref{RHS2}, we obtain
\begin{gather*}
{}^tA_\gamma {}^t( \J a_{N\rho_\gamma (\frac{m_0}{N})})_{\rho_\gamma (x)}\\
\qquad=-2\pi\sqrt{-1}a_{N\rho_\gamma (\frac{m_0}{N})}(\rho_\gamma(x))\bigl(\Omega_xA_\gamma^{-1}+{}^t(\J u_\gamma)_x\bigr)\left(N\rho_\gamma \left(\frac{m}{N}\right)-N\rho_\gamma (x)\right) .
\end{gather*}
On the other hand, by \eqref{Gamma-invariance=J1} and \eqref{Gamma-invariance=J2}, we have
\begin{equation}\label{Gamma-relation of Omega}
{}^tA_\gamma \Omega_{\rho_\gamma (x)}=\Omega_xA_\gamma^{-1}+{}^t\left(\J u_\gamma\right)_x.
\end{equation}
This proves the lemma.
\end{proof}

\begin{Remark}\label{equiv commutativity2}
By Remark~\ref{action preserves BS} and Lemma~\ref{equiv commutativity}, all $\frac{m}{N}\in F\cap \frac{1}{N}\Z^n$ are integrable if and only if the condition~(1), hence all conditions in Proposition~\ref{prop commutativity2} holds.
\end{Remark}

\section{The integrable case}\label{the integrable case}
In this section, we investigate the case where the almost complex structure is integrable in details. We use the setting and the notations introduced in the previous section.

\subsection[Definition and properties of vartheta\_m/N]{Definition and properties of $\boldsymbol{\vartheta_{\frac{m}{N}}}$}
Let $\frac{m}{N}\in F\cap \frac{1}{N}\Z^n$ be the point which is integrable, and $a_m$ the nontrivial solution of~\eqref{Dirac3} of the form~\eqref{a(x)} with $a_m(\frac{m}{N})=1$. For each $\gamma \in \Gamma$, define $a_{N\rho_\gamma (\frac{m}{N})}$ in such a~way that it satisfies~\eqref{Gamma-invariance}. As we showed in Lemma~\ref{equiv commutativity}, $a_{N\rho_\gamma (\frac{m}{N})}$ is a nontrivial solution of~\eqref{Dirac3} for~$N\rho_\gamma (\frac{m}{N})$. Then, we define the formal Fourier series $\vartheta_{\frac{m}{N}}$ by
\begin{equation}\label{theta}
\vartheta_{\frac{m}{N}}(x,y):=\sum_{\gamma \in \Gamma}a_{N\rho_\gamma (\frac{m}{N})}(x){\rm e}^{2\pi\sqrt{-1}N\rho_\gamma (\frac{m}{N})\cdot y}.
\end{equation}

\begin{Proposition}\label{expression of theta}\quad
\begin{itemize}\itemsep=0pt
\item[$(1)$] $\vartheta_{\frac{m}{N}}$ has the following expression:
\begin{align*}
\vartheta_{\frac{m}{N}}(x,y)
=\sum_{\gamma \in \Gamma}g_\gamma {\rm e}^{2\pi\sqrt{-1}\{\Theta (\frac{m}{N}, \gamma, x)+N\rho_\gamma (\frac{m}{N})\cdot y\}},
\end{align*}
where
\begin{align*}
\Theta \left(\frac{m}{N}, \gamma, x\right):={}&-\sum_{i=1}^nG_m^i\left(\frac{m_1}{N},\ldots ,\frac{m_{i-1}}{N},(\rho_{\gamma^{-1}}(x))_i,\ldots ,(\rho_{\gamma^{-1}}(x))_n\right)\\
&+N\left\{{\tilde g}_\gamma (\rho_{\gamma^{-1}}(x))-\rho_\gamma \left(\frac{m}{N}\right)\cdot u_\gamma (\rho_{\gamma^{-1}}(x))\right\}.
\end{align*}
\item[$(2)$] $\vartheta_{\frac{m}{N}}$ can be described as \smash{$\vartheta_{\frac{m}{N}}\!=\!\sum_{\gamma \in \Gamma}{\vardbtilde \rho}_\gamma \circ s_m\circ {\tilde \rho}_{\gamma^{-1}}$}, where $s_m$ is the section defined by~\eqref{s_m}.
\item[$(3)$] If $Y+XY^{-1}X$ is constant, then $\vartheta_{\frac{m}{N}}$ converges absolutely and uniformly on any compact set.
\end{itemize}
\end{Proposition}
\begin{proof}
(1) and (2) are obtained by \eqref{Gamma-invariance}, \eqref{a(x)}, \eqref{vardbtilde rho} and \eqref{s_m}.
Let us prove (3). By \eqref{condition for action} and \eqref{Gamma-invariance=J2}, we obtain
\[
{}^tA_{\gamma^{-1}}\bigl(Y+XY^{-1}X\bigr)^{-1}A_{\gamma^{-1}}=\bigl(Y+XY^{-1}X\bigr)^{-1}.
\]
By using this formula together with the assumption, the expression in (1) can be rewritten as
\begin{equation*}
\vartheta_{\frac{m}{N}}(x,y)=\sum_{\gamma \in \Gamma}g_\gamma {\rm e}^{2\pi\sqrt{-1} [\frac{\sqrt{-1}N}{2} (x-\rho_{\gamma} ( \frac{m}{N} ) )\cdot (Y+XY^{-1}X)^{-1} (x-\rho_{\gamma}(\frac{m}{N}) ) +\text{real part} ]}.
\end{equation*}
Since $\bigl(Y+XY^{-1}X\bigr)^{-1}$ is positive definite, there exists a positive constant $c>0$ such that \smash{$\bigl(Y+XY^{-1}X\bigr)^{-1}\ge cI$}. Then,
\begin{gather*}
 | g_\gamma {\rm e}^{2\pi\sqrt{-1} [\frac{\sqrt{-1}N}{2} (x-\rho_{\gamma} ( \frac{m}{N} ) )\cdot (Y+XY^{-1}X)^{-1} (x-\rho_{\gamma}(\frac{m}{N}) ) +\text{real part} ]} | \\
\qquad={\rm e}^{-N\pi (x-\rho_{\gamma} ( \frac{m}{N} ) )\cdot (Y+XY^{-1}X)^{-1} (x-\rho_{\gamma}(\frac{m}{N}) )}
\le {\rm e}^{-cN\pi \norm{ x-\rho_{\gamma} ( \frac{m}{N} )}^2}\\
\qquad= {\rm e}^{-cN\pi \norm{ x-\frac{l}{N}}^2}
\qquad \left(\text{put } l:=N\rho_{\gamma} \left( \frac{m}{N} \right) \right)\\
\qquad= \prod_{i=1}^n {\rm e}^{-cN\pi (x_i-\frac{l_i}{N} )^2}.
\end{gather*}
Hence, the series is dominated by \smash{$\prod_{i=1}^n\sum_{l_i\in \Z}{\rm e}^{-cN\pi (\frac{l_i}{N}-x_i)^2}$}. Any compact set is contained in a product of closed intervals $I_1\times \cdots \times I_n$, so it is sufficient to show that \smash{$\sum_{l\in \Z}{\rm e}^{-cN\pi (\frac{l}{N}-x)^2}$} converges uniformly on any closed interval $I$.
Suppose that $I$ is of the form $I:=[x_m, x_M]$. Set~${l_M:=\max \big\{ l\in \Z \mid \frac{l}{N}\in I\big\}}$ and $l_m:=\min \big\{ l\in \Z \mid \frac{l}{N}\in I\big\}$. On $I$, \smash{$\sum_{-k\le l\le k}{\rm e}^{-cN\pi (\frac{l}{N}-x)^2}$} can be estimated as
\begin{align*}
\sum_{-k\le l\le k}{\rm e}^{-cN\pi (\frac{l}{N}-x)^2}&=\bigg(\sum_{-k\le l<l_m}+\sum_{l_m\le l\le l_M}+\sum_{l_M\le l\le k}\bigg){\rm e}^{-cN\pi (\frac{l}{N}-x)^2}\\
&\le \sum_{-k\le l<l_m}{\rm e}^{\frac{-c\pi}{N} (l-Nx_m)^2}+(l_M-l_n+1)+\sum_{l_M< l\le k}{\rm e}^{\frac{-c\pi}{N} (l-Nx_M)^2}\\
&\le \int_{-k}^{l_m}{\rm e}^{\frac{-c\pi}{N}(\tau -Nx_m)^2}{\rm d}\tau +(l_M-l_n+1)+\int_{l_M}^k {\rm e}^{\frac{-c\pi}{N}(\tau -Nx_M)^2}{\rm d}\tau.
\end{align*}
It is well known that \smash{$\int_{-k}^{l_m}{\rm e}^{\frac{-c\pi}{N}(\tau -Nx_m)^2}{\rm d}\tau$} and \smash{$\int_{l_M}^k {\rm e}^{\frac{-c\pi}{N}(\tau -Nx_M)^2}{\rm d}\tau$} converge as ${k\to +\infty}$.\looseness=-1
\end{proof}

\begin{Lemma}\label{existence condition2}
Let $s$ be a section of ${L_0}^{\otimes N}$ with Fourier series of the form~\eqref{Fourier}. If $s$ is a~nontrivial $\Gamma$-equivariant solution of $0=Ds$, then there exists $\frac{m}{N}\in F\cap \frac{1}{N}\Z^n$ such that $m$ is integrable. Conversely, suppose that there exists ${\frac{m}{N}\in F\cap \frac{1}{N}\Z^n}$
such that $m$ is integrable
and $\vartheta_{\frac{m}{N}}$ converges absolutely and uniformly on any compact set. Then, $\vartheta_{\frac{m}{N}}$ is a nontrivial $\Gamma$-equivariant solution of $0=Ds$.
\end{Lemma}
\begin{proof}
Since $s=\sum_{l\in \Z^n}a_l(x){\rm e}^{2\pi\sqrt{-1}l\cdot y}$ is nontrivial solution of $0=Ds$, by Proposition~\ref{existence condition}, there exists $l\in \Z^n$ such that $a_l\not=0$. On the other hand, as is noticed in Remark~\ref{action preserves BS}, there exists $\left(\gamma ,\frac{m}{N}\right)\in \Gamma\times \left(F\cap \frac{1}{N}\Z^n\right)$ such that $l=N\rho_\gamma(\frac{m}{N})$. Since $s$ is $\Gamma$-equivariant, by \eqref{Gamma-invariance}, $0\not=a_l=a_{N\rho_\gamma(\frac{m}{N})}$ implies $a_m\not=0$.
The latter part follows from the definition of $\vartheta_{\frac{m}{N}}$.
\end{proof}

Let $\pi\colon (M,\omega)\to B$ be a Lagrangian fibration on a complete base $B$ with prequantum line bundle $\bigl(L,\nabla^L\bigr)\to (M,\omega)$. By Corollary~\ref{complete base}, they
are obtained as the quotient of an action of~${\Gamma:=\pi_1(B)}$ on $\bigl({L_0},\nabla^{L_0}\bigr)\to (M_0,\omega_0)$. Let $J$ be a compatible almost complex structure on~$(M,\omega)$ which is invariant along the fiber in the sense of Lemma~\ref{invariant J} and $D^M$ the associated Spin${}^c$ Dirac operator on $(M,N\omega)$ with coefficients in $L^{\otimes N}$. We denote by $D$ the Spin${}^c$ Dirac operator with coefficients in ${L_0}^{\otimes N}$ associated with the pull-back of $J$ to $M_0$. Since the $\Gamma$-action preserves all the data, $\Gamma \bigl(L^{\otimes N}\bigr)\cap \ker D^M$ is identified with \smash{$\bigl(\Gamma \bigl({L_0}^{\otimes N}\bigr)\cap \ker D\bigr)^\Gamma$}, the space of $\Gamma$-equivariant elements in $\Gamma \bigl({L_0}^{\otimes N}\bigr)\cap \ker D$.
If $J$ is integrable, so is the pull-back of $J$ to $M_0$. In this case, by Proposition~\ref{prop commutativity2}, all $\frac{m}{N}\in F\cap \frac{1}{N}\Z^n$ are integrable. So, one can consider $\vartheta_{\frac{m}{N}}$ for all $\frac{m}{N}\in F\cap \frac{1}{N}\Z^n$. By Lemma~\ref{existence condition2} and the above identification, if all $\vartheta_{\frac{m}{N}}$'s converge absolutely and uniformly on any compact set, then they can be thought of as elements of $\Gamma \bigl(L^{\otimes N}\bigr)\cap \ker D^M$, i.e., holomorphic sections of $L^{\otimes N}$ indexed by $B_{\rm BS}$. (As we noticed in Remark~\ref{action preserves BS}, $F\cap \frac{1}{N}\Z^n$ is identified with $B_{\rm BS}$.)

We choose the orientation on $M$ so that \smash{$(-1)^{\frac{n(n-1)}{2}}\frac{(N\omega)^n}{n!}$} is a positive volume form, and define the Hermitian inner product on $\Gamma \bigl(L^{\otimes N}\bigr)$ by
\[
(s, s')_{L^2(L^{\otimes N})}:=\int_M\langle s, s'\rangle_{L^{\otimes N}} (-1)^{\frac{n(n-1)}{2}}\dfrac{(N\omega)^n}{n!}
\]
for $s, s'\in \Gamma \bigl(L^{\otimes N}\bigr)$, where $\langle \cdot ,\ \cdot \rangle_{L^{\otimes N}}$ is the Hermitian metric of $L^{\otimes N}$. For $s\in \Gamma \bigl(L^{\otimes N}\bigr)$, we denote its $L^2$-norm by
\begin{equation*}
\norm{s}_{L^2 (L^{\otimes N} )}:= \{ (s, s )_{L^2 (L^{\otimes N} )} \}^{\frac{1}{2}}
\end{equation*}
and denote the space of $L^2$-sections of $L^{\otimes N}$ by $L^2\bigl(L^{\otimes N}\bigr)$. Then, we have the following theorem.
\begin{Theorem}\label{main2}
Let $\pi\colon (M,\omega)\to B$ be a Lagrangian fibration on a complete base $B$ and $\bigl(L,\nabla^L\bigr)\to (M,\omega)$ a prequantum line bundle. Let $J$ be a compatible integrable almost complex structure on $(M,\omega)$ which is invariant along the fiber in the sense of Lemma~{\rm\ref{invariant J}} and $D^M$ the associated Spin${}^c$ Dirac operator on $(M,N\omega)$ with coefficients in $L^{\otimes N}$ as above.
Assume that~\smash{$\vartheta_{\frac{m}{N}}$} converges absolutely and uniformly on any compact set and is square integrable as a section of~$L^{\otimes N}$ for each $\frac{m}{N}\in F\cap \frac{1}{N}\Z^n$.
Then, $L^2\bigl(L^{\otimes N}\bigr)\cap \ker D^M$ is a Hilbert space and \smash{$\{\vartheta_{\frac{m}{N}}\}_{\frac{m}{N}\in F\cap \frac{1}{N}\Z^n}$} is a complete orthogonal system of $L^2\bigl(L^{\otimes N}\bigr)\cap \ker D^M$ indexed by the Bohr--Sommerfeld points.
\end{Theorem}

\begin{proof}
By the definition of \smash{$\vartheta_{\frac{m}{N}}$} and the assumption of Theorem~\ref{main2}, \smash{$\{ \vartheta_{\frac{m}{N}}\}_{\frac{m}{N}\in F\cap \frac{1}{N}\Z^n}$} is an~orthogonal system of $L^2\bigl(L^{\otimes N}\bigr)$.
Suppose that \smash{$\mathrm{l.h.}(\{ \vartheta_{\frac{m}{N}}\}_{\frac{m}{N}\in F\cap \frac{1}{N}\Z^n})$} is the subspace of $L^2\bigl(L^{\otimes N}\bigr)$ generated by \smash{$\{ \vartheta_{\frac{m}{N}}\}_{\frac{m}{N}\in F\cap \frac{1}{N}\Z^n}$}, namely,
\[
\mathrm{l.h.}(\{ \vartheta_{\frac{m}{N}}\}_{\frac{m}{N}\in F\cap \frac{1}{N}\Z^n}):=\left\{ \sum_{i=1}^kc_i\vartheta_{\frac{m^i}{N}}\,\Big|\, k\in \N ,\, c_i\in \C ,\, \frac{m^i}{N}\in F\cap \frac{1}{N}\Z^n\right\},
\]
and we denote the closure of \smash{$\mathrm{l.h.}(\{ \vartheta_{\frac{m}{N}}\}_{\frac{m}{N}\in F\cap \frac{1}{N}\Z^n})$} in $L^2\bigl(L^{\otimes N}\bigr)$ by \smash{$\overline{\mathrm{l.h.}(\{ \vartheta_{\frac{m}{N}}\}_{\frac{m}{N}\in F\cap \frac{1}{N}\Z^n})}$}. Then, \smash{$\overline{\mathrm{l.h.}(\{ \vartheta_{\frac{m}{N}}\}_{\frac{m}{N}\in F\cap \frac{1}{N}\Z^n})}$} is described as
\begin{align*}
&\overline{\mathrm{l.h.}(\{ \vartheta_{\frac{m}{N}}\}_{\frac{m}{N}\in F\cap \frac{1}{N}\Z^n})} \\
&\qquad=\bigg\{ \sum_{\frac{m}{N}\in F\cap \frac{1}{N}\Z^n}c_{\frac{m}{N}}\vartheta_{\frac{m}{N}}\,\Big|\, c_{\frac{m}{N}}\in \C , \sum_{\frac{m}{N}\in F\cap \frac{1}{N}\Z^n}c_{\frac{m}{N}}\vartheta_{\frac{m}{N}}\text{ converges in }L^2\bigl(L^{\otimes N}\bigr)\bigg\} .
\end{align*}
In fact, any $\varphi=\sum_{\frac{m}{N}\in F\cap \frac{1}{N}\Z^n}c_{\frac{m}{N}}\vartheta_{\frac{m}{N}}$ in the right-hand side satisfies
\[
\lim_{k\to \infty}\Bigg\lVert \varphi -\sum_{\frac{m}{N}\in F\cap \frac{1}{N}\Z^n,\ \abs{m}\le k}c_{\frac{m}{N}}\vartheta_{\frac{m}{N}}\Bigg\rVert_{L^2(L^{\otimes N})}=0.
\]
This implies $\varphi$ is contained by the left-hand side. Conversely, since \smash{$\{ \vartheta_{\frac{m}{N}}\}_{\frac{m}{N}\in F\cap \frac{1}{N}\Z^n}$} is a completely orthogonal system of the subspace \smash{$\overline{\mathrm{l.h.}(\{ \vartheta_{\frac{m}{N}}\}_{\frac{m}{N}\in F\cap \frac{1}{N}\Z^n})}$}, then any $\varphi$ from the subspace \smash{$\overline{\mathrm{l.h.}(\{ \vartheta_{\frac{m}{N}}\}_{\frac{m}{N}\in F\cap \frac{1}{N}\Z^n})}$} satisfies
\[
\lim_{k\to \infty}\Bigg\lVert\varphi -\sum_{\frac{m}{N}\in F\cap \frac{1}{N}\Z^n,\abs{m}\le k}\dfrac{(\varphi ,\vartheta_{\frac{m}{N}})_{L^{\otimes N}}}{\norm{\vartheta_{\frac{m}{N}}}^2_{L^2(L^{\otimes N})}}\vartheta_{\frac{m}{N}}\Bigg\rVert_{L^2(L^{\otimes N})}=0.
\]
This implies $\varphi$ is contained by the right-hand side. We show that
\[
L^2\bigl(L^{\otimes N}\bigr)\cap \ker D^M=\smash{\overline{\mathrm{l.h.}(\{ \vartheta_{\frac{m}{N}}\}_{\frac{m}{N}\in F\cap \frac{1}{N}\Z^n})}}.
\]
 Let $s$ be an element of $L^2\bigl(L^{\otimes N}\bigr)\cap \ker D^M$. We think of $s$ as an element of \smash{$\bigl(\Gamma \bigl({L_0}^{\otimes N}\bigr)\cap \ker D\bigr)^\Gamma$}. By Lemma~\ref{Gamma-invariance-lemma}, $s$ can be written as in \eqref{Gamma-equiv section}. Then,
\begin{align*}
s(x,y)
&\stackrel{\eqref{Gamma-equiv section}}{=}\sum_{ (\gamma ,\frac{m}{N} )\in \Gamma \times (F\cap \frac{1}{N}\Z^n )}\hspace*{-2mm}g_{\gamma}a_m (\rho_{\gamma^{-1}}(x) ){\rm e}^{2\pi\sqrt{-1}N \{{\tilde g}_\gamma (\rho_{\gamma^{-1}}(x) )-\rho_\gamma (\frac{m}{N})\cdot u_\gamma (\rho_{\gamma^{-1}}(x) )+\rho_\gamma (\frac{m}{N})\cdot y \}}\\
&\stackrel{\eqref{a(x)}}{=}\sum_{ (\gamma ,\frac{m}{N} )\in \Gamma \times (F\cap \frac{1}{N}\Z^n )}g_\gamma a_m(\frac{m}{N}){\rm e}^{2\pi\sqrt{-1} \{\Theta (\frac{m}{N}, \gamma, x)+N\rho_\gamma (\frac{m}{N})\cdot y \}}\\
&=\sum_{\frac{m}{N}\in F\cap \frac{1}{N}\Z^n}a_m\left(\frac{m}{N}\right)\sum_{\gamma \in \Gamma}g_\gamma {\rm e}^{2\pi\sqrt{-1} \{\Theta (\frac{m}{N}, \gamma, x)+N\rho_\gamma (\frac{m}{N})\cdot y \}}\\
&=\sum_{\frac{m}{N}\in F\cap \frac{1}{N}\Z^n}a_m\left(\frac{m}{N}\right) \vartheta_{\frac{m}{N}}(x,y).
\end{align*}
(Note that it is well known that the Fourier series of $s$ pointwise converges absolutely. In particular, the order of terms of the Fourier series of $s$ in interchangeable.)
This implies $s$ is contained by \smash{$\overline{\mathrm{l.h.}(\{ \vartheta_{\frac{m}{N}}\}_{\frac{m}{N}\in F\cap \frac{1}{N}\Z^n})}$}.
Conversely, for any \smash{$s=\sum_{\frac{m'}{N}\in F\cap \frac{1}{N}\Z^n}c_{\frac{m'}{N}}\vartheta_{\frac{m'}{N}}$} in \smash{$\overline{\mathrm{l.h.}(\{ \vartheta_{\frac{m}{N}}\}_{\frac{m}{N}\in F\cap \frac{1}{N}\Z^n})}$}, let
\[
s=\sum_{(\gamma ,\frac{m}{N})\in \Gamma \times (F\cap \frac{1}{N}\Z^n)}b_{N\rho_\gamma (\frac{m}{N})}{\rm e}^{2\pi\sqrt{-1}N\rho_\gamma (\frac{m}{N})\cdot y}
\]
 be the Fourier series of $s$ with respect to $y_i$'s. Then, each $b_{N\rho_\gamma (\frac{m}{N})}$ is described by
\begin{align*}
b_{N\rho_\gamma (\frac{m}{N})}:={}&\int_{T^n}s(x,y){\rm e}^{-2\pi\sqrt{-1}N\rho_\gamma (\frac{m}{N})\cdot y}{\rm d}y\\
={}&\int_{T^n}\sum_{\frac{m'}{N}\in F\cap \frac{1}{N}\Z^n}c_{\frac{m'}{N}}\vartheta_{\frac{m'}{N}}(x,y){\rm e}^{-2\pi\sqrt{-1}N\rho_\gamma (\frac{m}{N})\cdot y}{\rm d}y\\
={}&\sum_{\frac{m'}{N}\in F\cap \frac{1}{N}\Z^n}c_{\frac{m'}{N}}\int_{T^n}\vartheta_{\frac{m'}{N}}(x,y){\rm e}^{-2\pi\sqrt{-1}N\rho_\gamma (\frac{m}{N})\cdot y}{\rm d}y\\
={}&c_{\frac{m}{N}}a_{N\rho_\gamma (\frac{m}{N})}(x).
\end{align*}
For each $\bigl(\gamma ,\frac{m}{N}\bigr)\in \Gamma \times \bigl(F\cap \frac{1}{N}\Z^n\bigr)$, by the definition of \smash{$\vartheta_{\frac{m}{N}}$}, \smash{$a_{N\rho_\gamma (\frac{m}{N})}(x)$} is a nontrivial solution of \eqref{Dirac3} for $N\rho_\gamma \bigl(\frac{m}{N}\bigr)$. Hence, so is $c_{\frac{m}{N}}a_{N\rho_\gamma (\frac{m}{N})}(x)$. This implies $s$ satisfies $D^Ms=0$.
\end{proof}

The condition in Proposition~\ref{expression of theta}\,(3) also gives a sufficient condition on the square integrability of $\vartheta_{\frac{m}{N}}$ as a section of $L^{\otimes N}$. A proof will be given later in more general case. See Lemma~\ref{norm of theta}.
\begin{Proposition}\label{square-integrability}
If $Y+XY^{-1}X$ is constant, then $\vartheta_{\frac{m}{N}}$ is square integrable as a section of~$L^{\otimes N}$.
\end{Proposition}

In particular, by Propositions~\ref{expression of theta} and~\ref{square-integrability}, and Theorem~\ref{main2}, we obtain the following corollary.
\begin{Corollary}\label{cor_main2}
Let $\pi\colon (M,\omega)\to B$ be a Lagrangian fibration on a complete base $B$ and $\bigl(L,\nabla^L\bigr)\to (M,\omega)$ a prequantum line bundle. Let $J$ be a compatible almost complex structure on $(M,\omega)$ which is invariant along the fiber in the sense of Lemma~{\rm\ref{invariant J}}. If $J$ is integrable and $Y+XY^{-1}X$ is constant, then \smash{$\{\vartheta_{\frac{m}{N}}\}_{\frac{m}{N}\in F\cap \frac{1}{N}\Z^n}$} is a complete orthogonal system of the space of square integrable holomorphic sections of $\bigl(L,\nabla^L\bigr)^{\otimes N}\to (M,N\omega, J)$ indexed by the Bohr--Sommerfeld points.
\end{Corollary}

Let us consider the special case where $\Gamma$ is trivial. In this case, $F=\R^n$, ${\bigl(L,\nabla^L\bigr)\!\to\! (M,\omega)\!\to\! B}$ is $\bigl({L_0},\nabla^{L_0}\bigr)\to (M_0,\omega_0)\to \R^n$, and $\vartheta_{\frac{m}{N}}$ is nothing but $s_m$ which is defined by \eqref{s_m} by Proposition~\ref{expression of theta}~(2).
Then, by Proposition~\ref{existence condition}, we have the following corollary.
\begin{Corollary}\label{main0}
Let $J$ be a compatible almost complex structure on $(M_0,\omega_0)$ which is invariant along the fiber in the sense of Lemma~{\rm\ref{invariant J}} and $D$ the associated Spin${}^c$ Dirac operator on~$(M_0,N\omega_0)$ with coefficients in ${L_0}^{\otimes N}$. Assume that $J$ is integrable and $s_m$ is in $L^2\bigl({L_0}^{\otimes N}\bigr)$ for all $m\in \Z^n$.
Then, $L^2\bigl({L_0}^{\otimes N}\bigr)\cap \ker D$ is a Hilbert space and $\left\{s_m\right\}_{m\in \Z^n}$ is a complete orthogonal system of $L^2\bigl({L_0}^{\otimes N}\bigr)\cap \ker D$. The latter assumption holds if $Y+XY^{-1}X$ is constant.
\end{Corollary}

\begin{Example}\label{EX4.8-2}
For Example~\ref{Ex4.8}, $Z=X+\sqrt{-1}Y$ can be chosen so that $Y+XY^{-1}X$ is a~constant map and $XY^{-1}$ and $Y^{-1}$ satisfy
\begin{gather*}
\bigl(XY^{-1}\bigr)_x=\bigl(Y+XY^{-1}X\bigr)
\begin{pmatrix}
u_{11}\cdot C^{-1}x & \cdots & u_{1n}\cdot C^{-1}x \\
\vdots & & \vdots \\
u_{n1}\cdot C^{-1}x & \cdots & u_{nn}\cdot C^{-1}x
\end{pmatrix}, \\
\bigl(Y^{-1}\bigr)_x=
\begin{pmatrix}
u_{11}\cdot C^{-1}x & \cdots & u_{1n}\cdot C^{-1}x \\
\vdots & & \vdots \\
u_{n1}\cdot C^{-1}x & \cdots & u_{nn}\cdot C^{-1}x
\end{pmatrix}
\bigl(Y+XY^{-1}X\bigr)\\
\phantom{\bigl(Y^{-1}\bigr)_x=}{}\times
\begin{pmatrix}
u_{11}\cdot C^{-1}x & \cdots & u_{1n}\cdot C^{-1}x \\
\vdots & & \vdots \\
u_{n1}\cdot C^{-1}x & \cdots & u_{nn}\cdot C^{-1}x
\end{pmatrix}+Y+XY^{-1}X.
\end{gather*}
In this case, $Y+XY^{-1}X$ is necessarily $I$ and $\Omega$ can be written as
\[
\Omega_x=
\begin{pmatrix}
u_{11}\cdot C^{-1}x & \cdots & u_{1n}\cdot C^{-1}x \\
\vdots & & \vdots \\
u_{n1}\cdot C^{-1}x & \cdots & u_{nn}\cdot C^{-1}x
\end{pmatrix}
+\sqrt{-1}\bigl(Y+XY^{-1}X\bigr)^{-1},
\]
and the condition~$(2)$ in Proposition~\ref{prop commutativity2} is equivalent to the following condition:
\begin{equation*}
\bigl({}^tC^{-1}u_{jk}\bigr)_i=\bigl({}^tC^{-1}u_{ik}\bigr)_j\qquad \text{for all }\quad i, j, k=1,\ldots , n.
\end{equation*}
Assume this condition as well as the condition $\frac{N}{2}v_i\cdot U_jv_i\in \Z$ for all $i, j=1,\ldots , n$. Then, for each $\frac{m}{N}\in F\cap \frac{1}{N}\Z^n$, $\vartheta_{\frac{m}{N}}$ is described by
\begin{align*}
\vartheta_{\frac{m}{N}}(x,y)=\sum_{\gamma \in \Gamma}g_\gamma {\rm e}^{2\pi\sqrt{-1}\{\Theta (\frac{m}{N}, \gamma, x)+N\rho_\gamma (\frac{m}{N})\cdot y\}},
\end{align*}
where
\begin{gather*}
\Theta \left(\frac{m}{N}, \gamma, x\right)
\\
\qquad=N\sum_{i=1}^n\sum_{j>i}\left( \rho_{\gamma^{-1}}(x)-\frac{m}{N}\right)_i \left( \rho_{\gamma^{-1}}(x)-\frac{m}{N}\right)_j\bigl( {}^tC^{-1}u_{ij}\bigr)\cdot
\begin{pmatrix}
\frac{m_1}{N} \\
\vdots \\
\frac{m_{i-1}}{N} \\
\frac{1}{2}\left( \rho_{\gamma^{-1}}(x)+\frac{m}{N}\right)_i \\
\left( \rho_{\gamma^{-1}}(x)\right)_{i+1} \\
\vdots \\
\left( \rho_{\gamma^{-1}}(x)\right)_n
\end{pmatrix}\\
\phantom{\qquad=}{} +\frac{N}{2}\sum_{i=1}^n\left( \rho_{\gamma^{-1}}(x)-\frac{m}{N}\right)_i^2\bigl( {}^tC^{-1}u_{ii}\bigr)\cdot
\begin{pmatrix}
\frac{m_1}{N} \\
\vdots \\
\frac{m_{i-1}}{N} \\
\frac{1}{3}\left(2\rho_{\gamma^{-1}}(x)+\frac{m}{N}\right)_i \\
\bigl( \rho_{\gamma^{-1}}(x)\bigr)_{i+1} \\
\vdots \\
\bigl( \rho_{\gamma^{-1}}(x)\bigr)_n
\end{pmatrix}\\
\phantom{\qquad=}{} + \frac{N}{2} \left\{\left( \rho_{\gamma^{-1}}(x)-\frac{m}{N}\right) \cdot \left(
\begin{pmatrix}
u_{11}\cdot \gamma & \cdots & u_{1n}\cdot \gamma \\
\vdots & & \vdots \\
u_{n1}\cdot \gamma & \cdots & u_{nn}\cdot \gamma
\end{pmatrix}
+ \sqrt{-1}\bigl(Y+XY^{-1}X\bigr)^{-1} \right)\right.\\
\left.\phantom{\qquad=}{}\times
\left( \rho_{\gamma^{-1}}(x)-\frac{m}{N}\right) - \frac{m}{N}\cdot
\begin{pmatrix}
u_{11}\cdot \gamma & \cdots & u_{1n}\cdot \gamma \\
\vdots & & \vdots \\
u_{n1}\cdot \gamma & \cdots & u_{nn}\cdot \gamma
\end{pmatrix}
\frac{m}{N}\right\}.
\end{gather*}
By Proposition~\ref{expression of theta}\,(3), \smash{$\vartheta_{\frac{m}{N}}$} converges absolutely and uniformly on any compact set.
\end{Example}

\subsection[The case when Z is constant]{The case when $\boldsymbol{Z}$ is constant}
Let $\pi\colon (M,\omega)\to B$ be a Lagrangian fibration on a complete $n$-dimensional $B$ with prequantum line bundle $\bigl(L,\nabla^L\bigr)\to (M,\omega)$. Then, it is obtained as the quotient of the $\Gamma:=\pi_1(B)$-action on $\pi_0\colon (M_0,\omega_0)\to \R^n$ with prequantum line bundle $\bigl({L_0}, \nabla^{L_0}\bigr)\to (M_0,\omega_0)$. Suppose that the $\Gamma$-actions are described by \eqref{Gamma-action on Rn times Tn} and \eqref{vardbtilde rho} as before. Let $J$ be a compatible almost complex structure on $(M,\omega)$ and $Z\in C^\infty (M_0, {\mathcal S}_n)$ be the map corresponding to the pull-back of $J$ to~$M_0$.
A situation in which $(2)$ in Proposition~\ref{prop commutativity2} holds occurs when $Z$ is a constant map. In this subsection, we discuss this case in detail. Note that in this case, $\J u_\gamma$ is a constant map for each $\gamma \in\Gamma$. It is obtained by \eqref{Gamma-invariance=J1}. Moreover, as a special case of the setting in the previous subsection, we can obtain the following theorem.
\begin{Theorem}\quad
\begin{itemize}\itemsep=0pt
\item[$(1)$] For each $\frac{m}{N}\in F\cap \frac{1}{N}\Z^n$, $\vartheta_{\frac{m}{N}}$ can be described as follows:
\begin{equation*}
\vartheta_{\frac{m}{N}}(x,y)=\sum_{\gamma \in \Gamma}g_\gamma {\rm e}^{2\pi\sqrt{-1}\{\Theta (\frac{m}{N}, \gamma, x)+N\rho_\gamma (\frac{m}{N})\cdot y\}},
\end{equation*}
where
\begin{align*}
\Theta \left(\frac{m}{N}, \gamma, x\right)={}&\frac{N}{2}\left\{\left(\rho_{\gamma^{-1}}(x)-\frac{m}{N}\right)\cdot \bigl(\Omega+{}^tA_\gamma \J u_\gamma\bigr)\left(\rho_{\gamma^{-1}}(x)-\frac{m}{N}\right)\right.\\
&\left.-\frac{m}{N}\cdot \bigl({}^tA_\gamma \J u_\gamma\bigr)\frac{m}{N}\right\}-N\rho_{\gamma}\left(\frac{m}{N}\right)\cdot u_\gamma (0) .
\end{align*}
\item[$(2)$] For each $\frac{m}{N}\in F\cap \frac{1}{N}\Z^n$, $\vartheta_{\frac{m}{N}}$ converges absolutely and uniformly on any compact set.
\item[$(3)$] $J$ is integrable and
\smash{$\{ \vartheta_{\frac{m}{N}}\}_{\frac{m}{N}\in F\cap \frac{1}{N}\Z^n}$} gives a complete orthogonal system of the space of square integrable holomorphic sections of $\bigl(L,\nabla^L\bigr)^{\otimes N}\to (M,N\omega, J)$.
\end{itemize}
\end{Theorem}
\begin{proof}
(1) is obtained from Proposition~\ref{expression of theta}\,(1). (2) is obtained by the assumption and Proposition~\ref{expression of theta}\,(3). The first half of (3) holds since $J$ is covariant constant with respect to the associated Levi-Civita connection. The other half is obtained by Lemma~\ref{norm of theta} later and Corollary~\ref{cor_main2}.
\end{proof}

When $Z$ is constant, the associated Riemannian metric of $M$ is flat. So, by Bieberbach's theorem, if $M$ is compact, then $M$ is finitely covered by the $2n$-dimensional torus $T^{2n}$. In particular, \smash{$\vartheta_{\frac{m}{N}}$}'s should be obtained from classical theta functions. So, let us see how \smash{$\vartheta_{\frac{m}{N}}$}'s relate with classical theta functions for Example~\ref{B times T} with $C=I$, in which $M$ itself is $T^{2n}$. First, let us briefly recall classical theta functions. For each $T\in {\mathcal S}_n$ and $a, b\in \Q^n$, the theta function with rational characteristics is a holomorphic 
section on the trivial holomorphic line bundle $\C^n\times \C\to \C^n$ which is defined by
\[
\vartheta
\begin{bmatrix}
a \\
b
\end{bmatrix}
( z, T):=\sum_{\gamma \in \Z^n}{\rm e}^{\pi\sqrt{-1}(\gamma+a)\cdot T ( \gamma+a)+2\pi\sqrt{-1}( \gamma+a)\cdot ( z+b)} .
\]
It is well known that $\vartheta
\begin{bmatrix}
a \\
b
\end{bmatrix}
( z, T)$ has the following quasi-periodicity:
\begin{gather*}
\vartheta
\begin{bmatrix}
a \\
b
\end{bmatrix}
( z+m, T)={\rm e}^{2\pi\sqrt{-1}a\cdot m}
\vartheta
\begin{bmatrix}
a \\
b
\end{bmatrix}
( z, T) ,\\
\vartheta
\begin{bmatrix}
a \\
b
\end{bmatrix}
( z+Tm, T)={\rm e}^{-2\pi\sqrt{-1}b\cdot m}{\rm e}^{-\pi\sqrt{-1}m\cdot Tm-2\pi\sqrt{-1}m\cdot z}
\vartheta
\begin{bmatrix}
a \\
b
\end{bmatrix}
( z, T)
\end{gather*}
for $m\in \Z^n$. For more details, see \cite[Chapter II, Section~1]{Mumford1} and \cite[Section~2]{Mumford3}.
Here we need the case where $T=N\Omega$, $a=\frac{m}{N}$, and $b=0$. In this case, define the $\Z^{2n}=\Z^n\times \Z^n$-action on~${\C^n\times \C\to \C^n}$ by
\begin{equation*}
(\gamma,\gamma')\cdot (z,w):=\bigl(z+N(-\Omega\gamma+\gamma'), {\rm e}^{-\pi\sqrt{-1}N\gamma \cdot \Omega\gamma +2\pi\sqrt{-1}\gamma \cdot z}w\bigr)
\end{equation*}
for $(\gamma,\gamma')\in \Z^{2n}$ and $(z,w)\in \C^n\times \C$. Also define the $\Z^{2n}$-action on the trivial complex line bundle $\R^{2n}\times \C\to \R^{2n}$ by
\begin{equation}\label{Z2n-action on R2n}
(\gamma,\gamma')\cdot (x,y,w):=\bigl(x+\gamma ,y+\gamma', {\rm e}^{2\pi\sqrt{-1}N\gamma\cdot y}w\bigr)
\end{equation}
for $(\gamma,\gamma')\in \Z^{2n}$ and $(x,y,w)\in \R^{2n}\times \C$. Note that by taking the quotient of the latter $\Z^n$-action of \eqref{Z2n-action on R2n}, we can recover Example~\ref{B times T} with $C=I$ and $g_\gamma=1$.
Let $F\colon \R^{2n}\to \C^n$ and~${{\tilde F}\colon \R^{2n}\times \C\to \C^n\times \C}$ be the $\R$-linear isomorphism and the bundle isomorphism covering~$F$ which are defined by
\begin{align*}
F(x,y)&:=N(-\Omega x+y) ,\qquad
{\tilde F}(x,y,w):=\bigl(N(-\Omega x+y) , {\rm e}^{-\pi\sqrt{-1}N x\cdot \Omega x}w\bigr) .
\end{align*}
Then, the direct computation shows the following theorem.
\begin{Theorem}\label{relation with theta}\quad
\begin{itemize}\itemsep=0pt
\item[$(1)$] $J_{\sqrt{-1}I}\circ F=F\circ (J_Z)$, i.e., $F$ is a $\C$-linear isomorphism from $\bigl(\R^{2n}, J_Z\bigr)$ to the standard complex vector space \smash{$(\C^n,J_{\sqrt{-1}I})$}.
\item[$(2)$] ${\tilde F}$ is equivariant with respect to the $\Z^{2n}$-actions defined above.
\item[$(3)$] $\vartheta_{\frac{m}{N}}$ satisfies \smash{${\tilde F}\circ \vartheta_{\frac{m}{N}}(x,y)=\vartheta
\left[
\begin{smallmatrix}
\frac{m}{N} \\
0
\end{smallmatrix}
\right]\left( F(x,y), N\Omega\right)$}, i.e.,
\[
\vartheta_{\frac{m}{N}}(x,y)={\rm e}^{\pi\sqrt{-1}N x\cdot \Omega x}\vartheta
\left[
\begin{matrix}
\frac{m}{N} \\
0
\end{matrix}
\right]\left( N(-\Omega x+y), N\Omega\right) .
\]
\end{itemize}
\end{Theorem}

\subsection{Adiabatic-type limit}\label{adiabatic limit}
In this subsection let us consider a one parameter family $\big\{\bigl(g^t, J^t\bigr)\big\}_{t>0}$ of the Riemannian metrics and the almost complex structures on a Lagrangian fibration so that the fiber shrinks as $t$ goes to $\infty$, and investigate the behavior of \smash{$\vartheta_{\frac{m}{N}}$} defined by \eqref{theta} as $t$ goes to $\infty$.

Let $Z=X+\sqrt{-1}Y\in C^\infty ( M_0,{\mathcal S}_n)$ be the map independent of $y_1, \ldots ,y_n$. Let $J=J_Z$ be the corresponding compatible almost complex structure on $(M_0,\omega_0)$. For each $t>0$, we define the almost complex structure $J^t$ by
\begin{align*}
J^tu:=( -J\partial_y, \partial_y)
\begin{pmatrix}
0 & \frac{-1}{t}\\
t & 0
\end{pmatrix}
\begin{pmatrix}
u_H\\
u_V
\end{pmatrix}
\end{align*}
for $u=(-J\partial_y, \partial_y)
\left(\begin{smallmatrix}
u_H\\
u_V
\end{smallmatrix}\right)\in T_{(x,y)}M_0$. It is easy to see the following lemma.
\begin{Lemma}\label{compatibility2}\quad
\begin{itemize}\itemsep=0pt
\item[$(1)$] For any $t>0$, $J^t$ is compatible with $\omega_0$. The map $Z^t\in C^\infty ( M_0,{\mathcal S}_n)$ corresponding to $J^t$ is described as
\begin{equation*}
Z^t=\left(\frac{1}{t}X+\sqrt{-1}Y\right)Y^{-1}\bigl(Y+XY^{-1}X\bigr)\left(tY+\frac{1}{t}XY^{-1}X\right)^{-1}Y.
\end{equation*}
$J^t$ can be also written as
\begin{gather*}
J^t\left((\partial_x, \partial_y)
\begin{pmatrix}
u_x \\
u_y
\end{pmatrix}\right)\\
\qquad=(\partial_x, \partial_y)
\frac{1}{t}
\begin{pmatrix}
XY^{-1} & -Y-XY^{-1}X \\
Y^{-1}\bigl(t^2Y+XY^{-1}X\bigr)\bigl(Y+XY^{-1}X\bigr)^{-1} & -Y^{-1}X
\end{pmatrix}
\begin{pmatrix}
u_x \\
u_y
\end{pmatrix}.
\end{gather*}
\item[$(2)$] For any $t>0$, let $g^t$ be the Riemannian metric corresponding to $\omega_0$ and $J^t$. Then, for $u=(-J\partial_y, \partial_y)
\left(\begin{smallmatrix}
u_H\\
u_V
\end{smallmatrix}\right)$, $v=(-J\partial_y, \partial_y)
\left(\begin{smallmatrix}
v_H\\
v_V
\end{smallmatrix}\right)
\in T_{(x,y)}M_0$, $g^t$ can be written by
\begin{align*}
g^t(u,v)={}&\omega_0\bigl( u,J^tv \bigr)
= t\bigl(0,{}^tu_H\bigr)
\begin{pmatrix}
Y^{-1} & -Y^{-1}X \\
-XY^{-1} & Y+XY^{-1}X
\end{pmatrix}
\begin{pmatrix}
0\\
v_H
\end{pmatrix}
\\
&
+\frac{1}{t}\bigl(0,{}^tu_V\bigr)
\begin{pmatrix}
Y^{-1} & -Y^{-1}X \\
-XY^{-1} & Y+XY^{-1}X
\end{pmatrix}
\begin{pmatrix}
0\\
v_V
\end{pmatrix}.
\end{align*}
\end{itemize}
\end{Lemma}

Suppose that a group $\Gamma$ acts on $\pi_0\colon (M_0,\omega_0)\to \R^n$ and the $\Gamma$-actions $\rho$ on $\R^n$ and $\tilde \rho$ on~${(M_0,\omega_0)}$ are written as in \eqref{Gamma-action on Rn times Tn}. 
\begin{Lemma}
The $\Gamma$-action $\tilde \rho$ preserves $J^t$ $($hence, $g^t)$ for all $t>0$ if and only if $\tilde \rho$ preserves~$J$.
\end{Lemma}

For $J^t$ and $g^t$ defined as above, the same arguments in Section~\ref{Spin-c Dirac operator} goes well, just by replacing~$J$,~$g$ by $J^t$, $g^t$. For each $t>0$, let \smash{$\vartheta^t_{\frac{m}{N}}$} be the one defined by \eqref{theta} for $J^t$ and $g^t$. Let us investigate the behavior of \smash{$\vartheta^t_{\frac{m}{N}}$} as $t$ goes to infinity.
For $t>0$, $\Omega^t$ defined by \eqref{Omega} for $Z^t$ can be described~as
\begin{equation}\label{Omega-t}
\Omega^t=\bigl(Y+XY^{-1}X\bigr)^{-1}\bigl(X+t\sqrt{-1}Y\bigr)Y^{-1}.
\end{equation}
Let $D^t$ be the corresponding Spin${}^c$ Dirac operator. Then, for a section $s$ of ${L_0}^{\otimes N}$, $D^ts$ can be described as
\[
D^ts=-\frac{\sqrt{-1}}{N}\sum_{i=1}^n\partial_{y_i}\otimes_{\C}\biggl\{\partial_{x_i}s+\sum_{j=1}^n\bigl(\Omega^t\bigr)_{ij}\bigl(\partial_{y_j}s -2\pi\sqrt{-1}Nx_js\bigr)\biggr\}.
\]
It is clear that
\begin{Lemma}
For any $t>0$, the condition $(2)$ in Proposition~{\rm\ref{prop commutativity2}} holds for $\Omega^t$ if and only if it holds for $\Omega=\Omega^1$. In particular, $J^t$ is integrable if and only if $J$ is integrable.
\end{Lemma}
Suppose that $\pi_0\colon (M_0,N\omega_0,J)\to \R^n$ with prequantum line bundle $\bigl({L_0}, \nabla^{L_0}\bigr)^{\otimes N}\to (M_0,$ $N\omega_0, J)$ is equipped with an action of a group $\Gamma$ which preserves all the data, and the $\Gamma$-actions are described by \eqref{Gamma-action on Rn times Tn} and \eqref{vardbtilde rho} as before. We assume that the $\Gamma$-action $\rho$ on $\R^n$ is properly discontinuous and free. Let $\pi \colon (M,N\omega)\to B$ and \smash{$\bigl(L,\nabla^L\bigr)^{\otimes N}\to (M,N\omega)$} be the Lagrangian fibration and the prequantum line bundle on it obtained by the quotient of the $\Gamma$-action.
On~$M$, we define the $L^p$-norm of a section $s$ of $L^{\otimes N}$ by
\[
\norm{s}_{L^p(L^{\otimes N})}:=\left(\int_M\langle s, s\rangle_{L^{\otimes N}}^{\frac{p}{2}} (-1)^{\frac{n(n-1)}{2}}\dfrac{(N\omega)^n}{n!}\right)^\frac{1}{p},
\]
where $\langle \cdot ,\ \cdot \rangle_{L^{\otimes N}}$ is the Hermitian metric of $L^{\otimes N}$ which is induced from the Hermitian metric~${\langle \cdot ,\ \cdot \rangle_{{L_0}^{\otimes N}}}$ of ${L_0}^{\otimes N}$. As noticed in Remark~\ref{Hermitian metric},
there exists a positive constant $C$ such that~${\langle \cdot ,\ \cdot \rangle_{{L_0}^{\otimes N}}}$ can be written as $\langle \cdot ,\ \cdot \rangle_{{L_0}^{\otimes N}} =C \langle \cdot ,\ \cdot \rangle_{\C}$, where $\langle \cdot ,\ \cdot \rangle_{\C}$ is the standard Hermitian inner product on $\C$.

For each $t>0$ and each point $\frac{m}{N}\in F\cap \frac{1}{N}\Z^n$ which is integrable, the corresponding \smash{$\vartheta^t_{\frac{m}{N}}$} is defined by \eqref{theta} for $\Omega^t$. We identify $F\cap \frac{1}{N}\Z^n$ with $B_{\rm BS}$ the set of Bohr--Sommerfeld points of $\pi \colon (M,N\omega)\to B$ with prequantum line bundle \smash{$\bigl(L,\nabla^L\bigr)^{\otimes N}\to (M,N\omega)$} and identify \smash{$\vartheta^t_{\frac{m}{N}}$} with the section of \smash{$\bigl(L,\nabla^L\bigr)^{\otimes N}\to (M,N\omega)$} which is induced from \smash{$\vartheta^t_{\frac{m}{N}}$}. Then, concerning the $L^p$-norm, we have the following lemma.
\begin{Lemma}\label{norm of theta}
Suppose that $Y+XY^{-1}X$ is constant. Then, as a section of $\bigl(L,\nabla^L\bigr)^{\otimes N}\to (M,N\omega)$, the $L^p$-norm of \smash{$\vartheta^t_{\frac{m}{N}}$} converges and it can be calculated as follows:
\begin{equation*}
\big\|\vartheta^t_{\frac{m}{N}}\big\|_{L^p(L^{\otimes N})}^p=C\sqrt{\det \bigl(Y+XY^{-1}X\bigr)}\left(\frac{N}{pt}\right)^\frac{n}{2}.
\end{equation*}
\end{Lemma}
\begin{proof}
Let $o(B)$ be the orientation bundle of $B$ which is defined as the quotient bundle of the trivial real line bundle $\R^n\times \R \to \R^n$ on the universal cover of $B$ by the $\Gamma$-action $\rho'_\gamma (x,r):=\left(\rho_\gamma (x), (\det A_\gamma) r\right)$ for $\gamma \in \Gamma$ and $(x,r)\in \R^n\times \R$. Then, we have a push-forward map~${\pi_*\colon \Omega^k(M)\to \Omega^{k-n}(B,o(B))}$, where $\Omega^\bullet (B,o(B))$ is the de Rham complex twisted by $o(B)$. $B$ has a natural density which we denote by $\abs{{\rm d}x}$. For densities, see \cite[Chapter I, Section~7]{BT}. Then,
\begin{align}
\big\|\vartheta^t_{\frac{m}{N}}\big\|_{L^p(L^{\otimes N})}^p&=\int_M\big\langle \vartheta^t_{\frac{m}{N}}, \vartheta^t_{\frac{m}{N}}\big\rangle_{L^{\otimes N}}^{\frac{p}{2}} (-1)^{\frac{n(n-1)}{2}}\dfrac{(N\omega)^n}{n!}\nonumber \\
&=\int_B\pi_*\left(\big\langle \vartheta^t_{\frac{m}{N}}, \vartheta^t_{\frac{m}{N}}\big\rangle_{L^{\otimes N}}^{\frac{p}{2}} (-1)^{\frac{n(n-1)}{2}}\dfrac{(N\omega)^n}{n!}\right) \nonumber\\
&=CN^n\sum_{\gamma \in \Gamma}\int_F {\rm e}^{-pN\pi t(\rho_{\gamma^{-1}}(x)-\frac{m}{N})\cdot (Y+XY^{-1}X)^{-1}(\rho_{\gamma^{-1}}(x)-\frac{m}{N})}\abs{{\rm d}x}.\label{L^p theta-1}
\end{align}
By changing the coordinates as $x'=\rho_{\gamma^{-1}}(x)$,
\begin{align}
\eqref{L^p theta-1}&=CN^n\sum_{\gamma \in \Gamma}\int_{\rho_{\gamma^{-1}}(F)} {\rm e}^{-pN\pi t(x'-\frac{m}{N} )\cdot (Y+XY^{-1}X)^{-1}(x'-\frac{m}{N} )} \abs{{\rm d}x'}\nonumber\\
&=CN^n\int_{\R^n} {\rm e}^{-pN\pi t(x'-\frac{m}{N} )\cdot (Y+XY^{-1}X)^{-1}(x'-\frac{m}{N} )} \abs{{\rm d}x'}.\label{L^p theta-2}
\end{align}
Since $Y+XY^{-1}X$ is positive definite, symmetric, there exists $P\in {\rm O}(n)$ such that
\[
Y+XY^{-1}X={}^tP
\begin{pmatrix}
\lambda_1 & & \\
& \ddots & \\
 & & \lambda_n
\end{pmatrix}P.
\]
Then, we define a positive definite symmetric matrix $\sqrt{Y+XY^{-1}X}$ by
\[
\sqrt{Y+XY^{-1}X}:={}^tP
\begin{pmatrix}
\sqrt{\lambda_1} & & \\
& \ddots & \\
 & & \sqrt{\lambda_n}
\end{pmatrix}P,
\]
and put \smash{$\tau:=\sqrt{\bigl(Y+XY^{-1}X\bigr)^{-1}}\bigl(x'-\frac{m}{N}\bigr)$}. Then,
\begin{align*}
\eqref{L^p theta-2}&=C\sqrt{\det\bigl(Y+XY^{-1}X\bigr)}N^n\int_{\R^n} {\rm e}^{-pN\pi t\norm{\tau}^2} \abs{{\rm d}\tau}\\
&=C\sqrt{\det\bigl(Y+XY^{-1}X\bigr)}N^n\prod_{i=1}^n\int_{-\infty}^{\infty} {\rm e}^{-pN\pi t\tau_i^2} {\rm d}\tau_i \\
&=C\sqrt{\det\bigl(Y+XY^{-1}X\bigr)}N^n\left(\sqrt{\dfrac{1}{pNt}}\right)^n.\tag*{\qed}
\end{align*} \renewcommand{\qed}{}
\end{proof}

We define the section $\delta_{\frac{m}{N}}$ of \smash{$\bigl(L, \nabla^L\bigr)^{\otimes N}|_{\pi^{-1}(\frac{m}{N})}$} by
\begin{equation}\label{delta_b}
\delta_{\frac{m}{N}}(y):=\frac{1}{C}{\rm e}^{2\pi\sqrt{-1}m\cdot y}.
\end{equation}
By Proposition~\ref{BS}, \smash{$\delta_{\frac{m}{N}}$} is a covariant constant section of \smash{$\bigl(L, \nabla^L\bigr)^{\otimes N}|_{\pi^{-1}(\frac{m}{N})}$}. Let $T_{\pi}^*M$ be the cotangent bundle along the fiber of $\pi$. On $\left(\wedge^nT_{\pi}^*M\right)\otimes \pi^*o(B)^*$, there exists a natural section, i.e., a density along the fiber of $\pi$, say $\abs{{\rm d}y}$, which satisfies $\int_{\pi^{-1}(x)}\abs{{\rm d}y}=1$ on each fiber of $\pi$. Then, we obtain the following theorem.
\begin{Theorem}\label{main3}
Suppose that $Y+XY^{-1}X$ is constant. Then, the section
\[\frac{\vartheta^t_{\frac{m}{N}}}{\big\|\vartheta^t_{\frac{m}{N}}\big\|_{L^1(L^{\otimes N})}}
\]
 converges to a delta-function section supported on the fiber $\pi^{-1}(\frac{m}{N})$ as $t$ goes to $\infty$ in the following sense: for any $L^2$-section $s$ of $L^{\otimes N}$,
\begin{equation*}
\lim_{t\to \infty}\left( s, \frac{\vartheta^t_{\frac{m}{N}}}{\big\|\vartheta^t_{\frac{m}{N}}\big\|_{L^1(L^{\otimes N})}} \right)_{L^2(L^{\otimes N})}=
\int_{\pi^{-1}(\frac{m}{N})}\langle s,\delta_{\frac{m}{N}}\rangle_{L^{\otimes N}} \abs{{\rm d}y}.
\end{equation*}
\end{Theorem}
\begin{proof}
We denote by ${\tilde s}$ the pull-back of $s$ to ${L_0}^{\otimes N}\to M_0$. Since ${\tilde s}$ is $\Gamma$-equivariant, the Fourier series of ${\tilde s}$ can be written as in \eqref{Gamma-equiv section}. Then, by using Proposition~\ref{expression of theta}\,(1),
\begin{gather}
\left( s, \frac{\vartheta^t_{\frac{m}{N}}}{\big\|\vartheta^t_{\frac{m}{N}}\big\|_{L^1(L^{\otimes N})}} \right)_{L^2(L^{\otimes N})}\nonumber \\
\qquad=\int_M\left\langle s, \frac{\vartheta^t_{\frac{m}{N}}}{\big\|\vartheta^t_{\frac{m}{N}}\big\|_{L^1(L^{\otimes N})}} \right\rangle_{L^{\otimes N}} (-1)^{\frac{n(n-1)}{2}}\dfrac{( N\omega)^n}{n!}\nonumber\\
\qquad=\int_B\pi_*\left(\left\langle s, \frac{\vartheta^t_{\frac{m}{N}}}{\big\|\vartheta^t_{\frac{m}{N}}\big\|_{L^1(L^{\otimes N})}} \right\rangle_{L^{\otimes N}} (-1)^{\frac{n(n-1)}{2}}\dfrac{( N\omega)^n}{n!}\right)\nonumber\\
\qquad= \frac{CN^n}{\big\|\vartheta^t_{\frac{m}{N}}\big\|_{L^1(L^{\otimes N})}}\nonumber\\
\phantom{\qquad=}{} \times \sum_{\gamma \in \Gamma}\int_F a_m(\rho_{\gamma^{-1}}(x))\overline{{\rm e}^{-2\pi\sqrt{-1}\sum_{i=1}^n G_m^i(\frac{m_1}{N},\ldots ,\frac{m_{i-1}}{N},(\rho_{\gamma^{-1}}(x))_i,\ldots ,(\rho_{\gamma^{-1}}(x))_n)}}\abs{{\rm d}x}.\label{delta-1}
\end{gather}
Here, we remark that we can interchange the operations to take infinite sums and integrals by Lemma~\ref{convergence in F times T}. By putting $x'=\rho_{\gamma^{-1}}(x)$, we have
\begin{align}
\eqref{delta-1}={}&\frac{CN^n}{\big\|\vartheta^t_{\frac{m}{N}}\big\|_{L^1(L^{\otimes N})}} \sum_{\gamma \in \Gamma}\int_{\rho_{\gamma^{-1}}(F)} a_m(x')\overline{{\rm e}^{-2\pi\sqrt{-1}\sum_{i=1}^n G_m^i(\frac{m_1}{N},\ldots ,\frac{m_{i-1}}{N},x'_i,\ldots ,x'_n)}}\abs{{\rm d}x'}\nonumber \\
={}&\frac{CN^n}{\big\|\vartheta^t_{\frac{m}{N}}\big\|_{L^1(L^{\otimes N})}}\int_{\R^n} a_m(x')\overline{{\rm e}^{-2\pi\sqrt{-1}\sum_{i=1}^n G_m^i(\frac{m_1}{N},\ldots ,\frac{m_{i-1}}{N},x'_i,\ldots ,x'_n)}}\abs{{\rm d}x'}\nonumber \\
={}&\frac{CN^n}{\big\|\vartheta^t_{\frac{m}{N}}\big\|_{L^1(L^{\otimes N})}}
\int_{\R^n} a_m(x'){\rm e}^{2\pi\sqrt{-1}\sum_{i=1}^n \re G_m^i(\frac{m_1}{N},\ldots ,\frac{m_{i-1}}{N},x'_i,\ldots ,x'_n)}\nonumber \\
& \times {\rm e}^{-\pi Nt(x'-\frac{m}{N})\cdot (Y+XY^{-1}X)^{-1}(x'-\frac{m}{N})}\abs{{\rm d}x'}.\label{delta-2}
\end{align}
We put
\[
f(x'):=a_m(x'){\rm e}^{2\pi\sqrt{-1}\sum_{i=1}^n \re G_m^i(\frac{m_1}{N},\ldots ,\frac{m_{i-1}}{N},x'_i,\ldots ,x'_n)}
\]
and
\smash{$\tau:=\sqrt{\bigl(Y+XY^{-1}X\bigr)^{-1}}\left(x'-\frac{m}{N}\right)$}.
By using Lemma~\ref{norm of theta} for $p=1$, \eqref{delta-2} can be written as follows:
\begin{align}
\eqref{delta-2}&=\frac{CN^n}{\big\|\vartheta^t_{\frac{m}{N}}\big\|_{L^1(L^{\otimes N})}}\int_{\R^n}f(x'){\rm e}^{-\pi Nt(x'-\frac{m}{N})\cdot (Y+XY^{-1}X)^{-1}(x'-\frac{m}{N})}\abs{{\rm d}x'}\nonumber\\
&=\frac{CN^n}{\big\|\vartheta^t_{\frac{m}{N}}\big\|_{L^1(L^{\otimes N})}}\sqrt{\det \bigl(Y+XY^{-1}X\bigr)}
\int_{\R^n}f\left(\sqrt{Y+XY^{-1}X}\tau +\frac{m}{N}\right){\rm e}^{-\pi Nt\norm{\tau}^2}\abs{{\rm d}\tau}\nonumber \\
&=(Nt)^{\frac{n}{2}}\int_{\R^n}f\left(\sqrt{Y+XY^{-1}X}\tau +\frac{m}{N}\right){\rm e}^{-\pi Nt\norm{\tau}^2}\abs{{\rm d}\tau}.\label{delta-3}
\end{align}
It is well known that $
\lim_{t\to \infty}\eqref{delta-3}=f\left(\frac{m}{N}\right)=a_m\left(\frac{m}{N}\right)$.
On the other hand, by using the expression
\[
{\tilde s}=\sum_{ (\gamma ,\frac{m'}{N} )\in \Gamma \times (F\cap \frac{1}{N}\Z^n )}a_{N\rho_{\gamma} (\frac{m'}{N} )}(x){\rm e}^{2\pi\sqrt{-1}N\rho_{\gamma} (\frac{m'}{N} )\cdot y},
\]
the right-hand side can be computed as
\begin{align*}
\int_{\pi^{-1}(\frac{m}{N})}\langle s,\delta_{\frac{m}{N}}\rangle_{L^{\otimes N}}\abs{{\rm d}y}
&=\int_{T^n}\big\langle {\tilde s},\delta_{\frac{m}{N}}\big\rangle_{{L_0}^{\otimes N}} \abs{{\rm d}y}\\
&= \sum_{(\gamma ,\frac{m'}{N})\in \Gamma \times (F\cap \frac{1}{N}\Z^n)} a_{N\rho_{\gamma}(\frac{m'}{N})} \left(\frac{m}{N}\right) \int_{T^n}{\rm e}^{2\pi\sqrt{-1}(N\rho_{\gamma}(\frac{m'}{N})-m)\cdot y}\abs{{\rm d}y}.
\end{align*}
The integral
\smash{$\int_{T^n}{\rm e}^{2\pi\sqrt{-1}(N\rho_{\gamma}(\frac{m'}{N})-m)\cdot y}\abs{{\rm d}y}
$}
 vanishes unless $\rho_{\gamma}(\frac{m'}{N})=\frac{m}{N}$. Since both $\frac{m'}{N}$ and~$\frac{m}{N}$ lie in the fundamental domain $F$, this implies $\gamma=e$ and $m'=m$, and in this case,
\[\int_{T^n}{\rm e}^{2\pi\sqrt{-1}(N\rho_{\gamma}(\frac{m'}{N})-m)\cdot y}\abs{{\rm d}y}=1.\]
 Thus,
\smash{$
\int_{\pi^{-1}(\frac{m}{N})}\langle s,\delta_{\frac{m}{N}}\rangle_{L^{\otimes N}}\abs{{\rm d}y}=a_m\left(\frac{m}{N}\right).
$}
This proves the theorem.
\end{proof}

\section{The non-integrable case}\label{approximation}
In this section, let us consider the case where the almost complex structure is not integrable. We still use the same notations introduced in Section~\ref{equiv quantization}.
By Lemma~\ref{basis of kerD}, the equation \eqref{Dirac3} has no smooth solution for $\frac{m}{N}\in F\cap \frac{1}{N}\Z^n$ if and only if $m$ is not integrable. For such $\frac{m}{N}$, instead of \eqref{Dirac3}, let us consider the following equation which is obtained by replacing $\Omega$ with its value~$\Omega_{\frac{m}{N}}$ at~$\frac{m}{N}$ in~\eqref{Dirac3}
\begin{equation}\label{approximated Dirac3}
0=
\begin{pmatrix}
\partial_{x_1}{\tilde a}_m\\
\vdots \\
\partial_{x_n}{\tilde a}_m
\end{pmatrix}
+2\pi\sqrt{-1}{\tilde a}_m\Omega_{\frac{m}{N}}
(m-Nx).
\end{equation}
The equation~\eqref{approximated Dirac3} has a solution of the form
\begin{equation*}
{\tilde a}_m(x)={\tilde a}_m\left(\frac{m}{N}\right){\rm e}^{\pi\sqrt{-1}N(x-\frac{m}{N})\cdot \Omega_{\frac{m}{N}}(x-\frac{m}{N})}.
\end{equation*}
We put the initial condition ${\tilde a}_m(\frac{m}{N})=1$ on the above ${\tilde a}_m$, and for each $\gamma \in \Gamma$, define ${\tilde a}_{N\rho_\gamma (\frac{m}{N})}$ in such a way that it satisfies \eqref{Gamma-invariance}.
\begin{Lemma}
${\tilde a}_{N\rho_\gamma (\frac{m}{N})}$ satisfies the following equality:
\begin{align}
0={}&
\begin{pmatrix}
\partial_{x_1}{\tilde a}_{N\rho_\gamma (\frac{m}{N})}(x)\\
\vdots \\
\partial_{x_n}{\tilde a}_{N\rho_\gamma (\frac{m}{N})}(x)
\end{pmatrix}
+2\pi\sqrt{-1}{\tilde a}_{N\rho_\gamma (\frac{m}{N})}(x)\Omega_x \left({N\rho_\gamma \left(\frac{m}{N}\right)}-Nx\right) \nonumber \\
&+2\pi\sqrt{-1}{\tilde a}_{N\rho_\gamma (\frac{m}{N})}(x){}^tA_{\gamma}^{-1}( \Omega_{\frac{m}{N}}-\Omega_{\rho_{\gamma^{-1}}(x)}) A_{\gamma}^{-1}\left({N\rho_\gamma \left(\frac{m}{N}\right)}-Nx\right).\label{approximated DE for a_m}
\end{align}
\end{Lemma}
\begin{proof}
By the same calculation as in the proof of Lemma~\ref{equiv commutativity}, we have
\begin{gather*}
{}^tA_{\gamma}
\begin{pmatrix}
\partial_{x_1}{\tilde a}_{N\rho_\gamma (\frac{m}{N})}( \rho_\gamma(x))\\
\vdots \\
\partial_{x_n}{\tilde a}_{N\rho_\gamma (\frac{m}{N})}( \rho_\gamma(x))
\end{pmatrix}\\
\qquad=-2\pi\sqrt{-1}{\tilde a}_{N\rho_\gamma (\frac{m}{N})}( \rho_\gamma(x))\bigl(\Omega_{\frac{m}{N}}A^{-1}_{\gamma}+{}^t( \J u_{\gamma})_x\bigr)
\left({N\rho_\gamma \left(\frac{m}{N}\right)}-N\rho_\gamma (x)\right).
\end{gather*}
\eqref{approximated DE for a_m} can be obtained from this equation and \eqref{Gamma-relation of Omega}.
\end{proof}

By using ${\tilde a}_{N\rho_\gamma (\frac{m}{N})}$'s, we define ${\tilde \vartheta}_{\frac{m}{N}}$ in the same manner as $\vartheta_{\frac{m}{N}}$, i.e.,
\[
{\tilde \vartheta}_{\frac{m}{N}}(x,y)=\sum_{\gamma \in \Gamma}{\tilde a}_{N\rho_\gamma (\frac{m}{N})}(x){\rm e}^{2\pi\sqrt{-1}N\rho_\gamma (\frac{m}{N})\cdot y}.
\]
${\tilde \vartheta}_{\frac{m}{N}}$ converges absolutely and uniformly on any compact set and can be written as
$
{\tilde \vartheta}_{\frac{m}{N}}=\sum_{\gamma \in \Gamma}{\vardbtilde \rho}_\gamma \circ s'_m\circ {\tilde \rho}_{\gamma^{-1}},
$
where $s'_m$ is the section defined by
\[
s'_m(x, y):={\rm e}^{\pi\sqrt{-1}N(x-\frac{m}{N})\cdot \Omega_{\frac{m}{N}}(x-\frac{m}{N})+2\pi\sqrt{-1}m\cdot y}.
\]
In particular, \smash{${\tilde \vartheta}_{\frac{m}{N}}$} defines an $L^p$-section of $L^{\otimes N}\to M$. Moreover, \smash{$\big\{{\tilde \vartheta}_{\frac{m}{N}}\big\}_{\frac{m}{N}\in F\cap \frac{1}{N}\Z^n}$} is an orthogonal system of the space of $L^2$-sections of $L^{\otimes N}$.
These can be proved by the same way as Proposition~\ref{expression of theta} and Lemma~\ref{norm of theta}.

Next let us consider the one parameter family of $J^t$ and $g^t$ defined in Section~\ref{adiabatic limit}. Then, corresponding to $J^t$ and $g^t$, we can obtain \smash{${\tilde \vartheta}_{\frac{m}{N}}^t$}, which can be explicitly described as
\begin{equation*}
{\tilde \vartheta}_{\frac{m}{N}}^t(x,y)=\sum_{\gamma\in \Gamma}g_\gamma {\rm e}^{2\pi\sqrt{-1}\{\Theta(\frac{m}{N},\gamma, x) +N\rho_\gamma (\frac{m}{N})\cdot y\}},
\end{equation*}
where
\begin{align*}
\Theta\left(\frac{m}{N}, \gamma ,x\right)={}&
\frac{N}{2}\left(\rho_{\gamma^{-1}}(x)-\frac{m}{N}\right)\cdot \Omega^t_{\frac{m}{N}}\left(\rho_{\gamma^{-1}}(x)-\frac{m}{N}\right)\\
&+N\left\{{\tilde g}_\gamma (\rho_{\gamma^{-1}}(x))-\rho_{\gamma}\left(\frac{m}{N}\right)\cdot u_\gamma (\rho_{\gamma^{-1}}(x))\right\}
\end{align*}
and \smash{$\Omega^t_{\frac{m}{N}}$} is the value of $\Omega^t$ given in \eqref{Omega-t} at $\frac{m}{N}$.
Then, \smash{${\tilde \vartheta}_{\frac{m}{N}}^t$} has the following property. The proof is same as Theorem~\ref{main3}.
\begin{Theorem}\label{main4}
For each $\frac{m}{N}\in F\cap \frac{1}{N}\Z^n$, the section
\[
 \frac{{\tilde \vartheta}^t_{\frac{m}{N}}}{\norm{{\tilde \vartheta}^t_{\frac{m}{N}}}_{L^1(L^{\otimes N})}}
 \]
 converges to a delta-function section supported on the fiber $\pi^{-1}(\frac{m}{N})$ as $t$ goes to $\infty$ in the following sense: for any $L^2$-section $s$ of $L^{\otimes N}$,
\begin{equation*}
\lim_{t\to \infty}\left(s, \frac{{\tilde \vartheta}^t_{\frac{m}{N}}}{\norm{{\tilde \vartheta}^t_{\frac{m}{N}}}_{L^1(L^{\otimes N})}}\right)_{L^2(L^{\otimes N})}
=\int_{\pi^{-1}(\frac{m}{N})}\langle s,\delta_{\frac{m}{N}}\rangle_{L^{\otimes N}} \abs{{\rm d}y}.
\end{equation*}
\end{Theorem}

Finally, let us investigate the behavior of \smash{$D^t {\tilde \vartheta}^t_{\frac{m}{N}}$} as $t$ goes to $\infty$. \smash{$D^t {\tilde \vartheta}^t_{\frac{m}{N}}$} is a section of ${\bigl(TM, J^t\bigr)\otimes_{\C}L^{\otimes N}}$, and $\bigl(TM, J^t\bigr)\otimes_{\C}L^{\otimes N}$ admits a Hermitian metric \smash{$\langle \cdot ,\ \cdot \rangle_{(TM, J^t)\otimes_{\C}L^{\otimes N}}$} induced by the one parameter version of \eqref{Hermitian metric on M} of $\bigl(TM, J^t\bigr)$ and the Hermitian metric of $L$. In terms of this Hermitian metric, the $L^2$-norm is defined as
\begin{equation*}
\big\|D^t {\tilde \vartheta}^t_{\frac{m}{N}}\big\|_{L^2((TM, J^t)\otimes_{\C}L^{\otimes N})}^2:=\int_M \big\langle D^t {\tilde \vartheta}^t_{\frac{m}{N}} , D^t {\tilde \vartheta}^t_{\frac{m}{N}} \big\rangle_{(TM, J^t)\otimes_{\C}L^{\otimes N}}(-1)^{\frac{n(n-1)}{2}}\frac{(N\omega)^n}{n!}.
\end{equation*}
In general, \smash{${\tilde \vartheta}^t_{\frac{m}{N}}$} is no longer a solution of $0=D^ts$, but we can show that \smash{${\tilde \vartheta}^t_{\frac{m}{N}}$} approximates the solution of this equation in the following sense:
\begin{Theorem}\label{main5}
\begin{equation*}
\lim_{t\to \infty}\big\|D^t {\tilde \vartheta}^t_{\frac{m}{N}}\big\|_{L^2((TM, J^t)\otimes_{\C}L^{\otimes N})}=0.
\end{equation*}
\end{Theorem}
\begin{proof}
For $n=1$, it is clear that all $m\in \Z$ are integrable. Thus, it is sufficient to prove the theorem for $n\ge 2$. By the definition of \smash{${\tilde \vartheta}^t_{\frac{m}{N}}$} and \eqref{approximated DE for a_m}, \smash{$D^t {\tilde \vartheta}^t_{\frac{m}{N}}$} can be written as
\begin{align*}
D^t {\tilde \vartheta}^t_{\frac{m}{N}}
={}&-\frac{\sqrt{-1}}{N}\sum_{i=1}^n\partial_{y_i}\otimes_{\C}\left\{\partial_{x_i}{\tilde \vartheta}^t_{\frac{m}{N}}+\sum_{j=1}^n(\Omega^t_x)_{ij}\left(\partial_{y_j}{\tilde \vartheta}^t_{\frac{m}{N}}-2\pi\sqrt{-1}Nx_j {\tilde \vartheta}^t_{\frac{m}{N}}\right)\right\} \\
={}&-\frac{\sqrt{-1}}{N}\sum_{i=1}^n\partial_{y_i}\otimes_{\C}\sum_{\gamma \in \Gamma}\left\{\partial_{x_i}{\tilde a}_{N\rho_\gamma (\frac{m}{N})}(x)\right.\\
&\left. +2\pi\sqrt{-1}{\tilde a}_{N\rho_\gamma (\frac{m}{N})}(x) \left(\Omega^t_x \left( N\rho_{\gamma}\left(\frac{m}{N}\right)-Nx\right)\right)_i \right\} {\rm e}^{2\pi\sqrt{-1}N\rho_\gamma (\frac{m}{N})\cdot y} \\
={}&-2\pi \sum_{i=1}^n\partial_{y_i}\otimes_{\C}\sum_{\gamma \in \Gamma}{\tilde a}_{N\rho_\gamma (\frac{m}{N})}(x)\left( B\left(\frac{m}{N},\gamma, x,t\right)\right)_i {\rm e}^{2\pi\sqrt{-1}N\rho_\gamma (\frac{m}{N})\cdot y},
\end{align*}
where
\begin{align*}
B\left(\frac{m}{N}, \gamma ,x,t\right)={}^tA_{\gamma}^{-1}\bigl( \Omega^t_{\frac{m}{N}}-\Omega^t_{\rho_{\gamma^{-1}}(x)}\bigr)\left(\frac{m}{N}-\rho_{\gamma^{-1}}(x)\right).
\end{align*}
Then,
\begin{gather*}
\big\langle D^t {\tilde \vartheta}^t_{\frac{m}{N}} ,D^t {\tilde \vartheta}^t_{\frac{m}{N}} \big\rangle_{(TM, J^t)\otimes_{\C}L}\\
\qquad=(2\pi)^2\sum_{\gamma_1, \gamma_2 \in \Gamma}\sum_{i_1, i_2}\big\langle {\tilde a}_{N\rho_{\gamma_1} (\frac{m}{N})}(x){\rm e}^{2\pi\sqrt{-1}N\rho_{\gamma_1} (\frac{m}{N})\cdot y}, {\tilde a}_{N\rho_{\gamma_2} (\frac{m}{N})}(x){\rm e}^{2\pi\sqrt{-1}N\rho_{\gamma_2} (\frac{m}{N})\cdot y}\big\rangle_{L^{\otimes N}}\\
\phantom{\qquad=}{}\times \left( B\left(\frac{m}{N}, \gamma_1 ,x,t\right)\right)_{i_1}\overline{ \left( B\left(\frac{m}{N}, \gamma_2 ,x,t\right)\right)_{i_2}}Ng^t(\partial_{y_{i_1}}, \partial_{y_{i_2}})\\
\qquad=(2\pi)^2\frac{N}{t}\sum_{\gamma_1, \gamma_2 \in \Gamma}\big\langle {\tilde a}_{N\rho_{\gamma_1}(\frac{m}{N})}(x){\rm e}^{2\pi\sqrt{-1}N\rho_{\gamma_1} (\frac{m}{N})\cdot y}, {\tilde a}_{N\rho_{\gamma_2}(\frac{m}{N})}(x){\rm e}^{2\pi\sqrt{-1}N\rho_{\gamma_2} (\frac{m}{N})\cdot y}\big\rangle_{L^{\otimes N}}\\
\phantom{\qquad=}{} \times B\left(\frac{m}{N}, \gamma_1 ,x,t\right)\cdot \bigl(Y+XY^{-1}X\bigr)_x \overline{ B\left(\frac{m}{N}, \gamma_2 ,x,t\right)}.
\end{gather*}
For each $x\in F$ and $u\in \C^n$, define the norm of $u$ with respect to $\bigl(Y+XY^{-1}X\bigr)_x$ by
\[
\norm{u}^2_{(Y+XY^{-1}X)_x}:=u\cdot \bigl(Y+XY^{-1}X\bigr)_x\overline{u}.
\]
By \eqref{Gamma-invariance=J2}, for each $\gamma \in \Gamma$, $\norm{u}^2_{(Y+XY^{-1}X)_x}$ satisfies
\[
\big\|{}^tA_{\gamma}u\big\|^2_{(Y+XY^{-1}X)_x}=\norm{u}^2_{(Y+XY^{-1}X)_{\rho_\gamma (x)}}.
\]
By using this norm, we obtain
\begin{gather*}
\big\|D^t {\tilde \vartheta}^t_{\frac{m}{N}}\big\|_{L^2((TM, J^t)\otimes_{\C}L^{\otimes N})}^2 \\
\qquad=(2\pi)^2\frac{C N^{n+1}}{t}\sum_{\gamma \in \Gamma}\int_F {\rm e}^{-2\pi Nt(\rho_{\gamma^{-1}}(x)-\frac{m}{N})\cdot ( Y+XY^{-1}X)^{-1}_{\frac{m}{N}}(\rho_{\gamma^{-1}}(x)-\frac{m}{N})} \\
\phantom{\qquad=}{} \times B\left(\frac{m}{N}, \gamma ,x,t\right)\cdot \bigl(Y+XY^{-1}X\bigr)_x \overline{ B\left(\frac{m}{N}, \gamma ,x,t\right)}\abs{{\rm d}x}\nonumber 
\\
\qquad=(2\pi)^2\frac{C N^{n+1}}{t}\sum_{\gamma \in \Gamma}\int_F {\rm e}^{-2\pi Nt(\rho_{\gamma^{-1}}(x)-\frac{m}{N})\cdot ( Y+XY^{-1}X)^{-1}_{\frac{m}{N}}(\rho_{\gamma^{-1}}(x)-\frac{m}{N})} \\
\phantom{\qquad=}{} \times \left\lVert B\left(\frac{m}{N}, \gamma ,x,t\right)\right\rVert^2_{(Y+XY^{-1}X)_x}\abs{{\rm d}x} \\
\qquad=(2\pi)^2\frac{C N^{n+1}}{t}\sum_{\gamma \in \Gamma}\int_F {\rm e}^{-2\pi Nt(\rho_{\gamma^{-1}}(x)-\frac{m}{N})\cdot ( Y+XY^{-1}X)^{-1}_{\frac{m}{N}}(\rho_{\gamma^{-1}}(x)-\frac{m}{N})} \\
\phantom{\qquad=}{} \times \left\lVert \bigl( \Omega^t_{\frac{m}{N}}-\Omega^t_{\rho_{\gamma^{-1}}(x)}\bigr)\left(\frac{m}{N}-\rho_{\gamma^{-1}}(x)\right)\right\rVert^2_{(Y+XY^{-1}X)_{\rho_{\gamma^{-1}}(x)}}\abs{{\rm d}x} \\
\qquad= (2\pi)^2\frac{C N^{n+1}}{t}\sum_{\gamma \in \Gamma}\int_{\rho_{\gamma^{-1}}(F)} {\rm e}^{-2\pi Nt(x'-\frac{m}{N})\cdot ( Y+XY^{-1}X)^{-1}_{\frac{m}{N}}(x'-\frac{m}{N})} \\
\phantom{\qquad=}{} \times \left\lVert \bigl( \Omega^t_{\frac{m}{N}}-\Omega^t_{x'}\bigr)\left(\frac{m}{N}-x'\right)\right\rVert ^2_{(Y+XY^{-1}X)_{x'}}\abs{{\rm d}x'}\qquad \left(\because x':=\rho_{\gamma^{-1}}(x)\right) \\
\qquad= (2\pi)^2\frac{C N^{n+1}}{t}\int_{\R^n} {\rm e}^{-2\pi Nt(x'-\frac{m}{N})\cdot ( Y+XY^{-1}X)^{-1}_{\frac{m}{N}}(x'-\frac{m}{N})}\\
\phantom{\qquad=}{} \times \left\lVert \bigl( \Omega^t_{\frac{m}{N}}-\Omega^t_{x'}\bigr)\left(\frac{m}{N}-x'\right)\right\rVert ^2_{(Y+XY^{-1}X)_{x'}}\abs{{\rm d}x'}.
\end{gather*}
Since $\Omega^t$ can be described as \eqref{Omega-t},
\begin{gather*}
\left\lVert \bigl( \Omega^t_{\frac{m}{N}}-\Omega^t_{x'}\bigr)\left(\frac{m}{N}-x'\right)\right\rVert ^2_{(Y+XY^{-1}X)_{x'}}\\
\qquad=\left\lVert \bigl(\re \bigl( \Omega_{\frac{m}{N}}-\Omega_{x'}\bigr)\bigr)\left(\frac{m}{N}-x'\right)\right\rVert ^2_{(Y+XY^{-1}X)_{x'}}\\
\phantom{\qquad=}{} +t^2\left\lVert \bigl(\im \bigl( \Omega_{\frac{m}{N}}-\Omega_{x'}\bigr)\bigr)\left(\frac{m}{N}-x'\right)\right\rVert ^2_{(Y+XY^{-1}X)_{x'}}.
\end{gather*}
We put
\begin{gather*}
R(x'):=\left\lVert \bigl(\re \bigl( \Omega_{\frac{m}{N}}-\Omega_{x'}\bigr)\bigr)\left(\frac{m}{N}-x'\right)\right\rVert ^2_{(Y+XY^{-1}X)_{x'}}, \\
I(x'):=\left\lVert \bigl(\im \bigl( \Omega_{\frac{m}{N}}-\Omega_{x'}\bigr)\bigr)\left(\frac{m}{N}-x'\right)\right\rVert ^2_{(Y+XY^{-1}X)_{x'}}.
\end{gather*}
By changing coordinates as
\[
\tau:=\sqrt{\bigl(Y+XY^{-1}X\bigr)_{\frac{m}{N}}^{-1}}\left(x'-\frac{m}{N}\right),
\]
 \smash{$\big\|D^t {\tilde \vartheta}^t_{\frac{m}{N}}\big\|_{L^2((TM, J^t)\otimes_{\C}L)}^2$} can be written by
\begin{gather*}
\big\|D^t {\tilde \vartheta}^t_{\frac{m}{N}}\big\|_{L^2((TM, J^t)\otimes_{\C}L)}^2 \\
\qquad{}
 = 2^{2-\frac{n}{2}}\pi^2 C N^{\frac{n}{2}+1}\sqrt{\det \bigl(Y+XY^{-1}X\bigr)_{\frac{m}{N}}}\\
\phantom{\qquad=}{}\times \left\{ t^{-1-\frac{n}{2}} \int_{\R^n} R\left(\sqrt{\bigl(Y+XY^{-1}X\bigr)_{\frac{m}{N}}}\tau +\frac{m}{N}\right)(2Nt)^\frac{n}{2} {\rm e}^{-2\pi Nt\norm{\tau}^2}\abs{{\rm d}\tau} \right.\\
\phantom{\qquad=}{}\left. +t^{1-\frac{n}{2}} \int_{\R^n} I\left(\sqrt{\bigl(Y+XY^{-1}X\bigr)_{\frac{m}{N}}}\tau +\frac{m}{N}\right)(2Nt)^\frac{n}{2} {\rm e}^{-2\pi Nt\norm{\tau}^2}\abs{{\rm d}\tau}\right\}.
\end{gather*}
It is well known that
\begin{align*}
& \lim_{t\to \infty}\int_{\R^n} R\left(\sqrt{\bigl(Y+XY^{-1}X\bigr)_{\frac{m}{N}}}\tau +\frac{m}{N}\right)(2Nt)^\frac{n}{2} {\rm e}^{-2\pi Nt\norm{\tau}^2}\abs{{\rm d}\tau}=R\left(\frac{m}{N}\right)=0, \\
& \lim_{t\to \infty} \int_{\R^n} I\left(\sqrt{\bigl(Y+XY^{-1}X\bigr)_{\frac{m}{N}}}\tau +\frac{m}{N}\right)(2Nt)^\frac{n}{2} {\rm e}^{-2\pi Nt\norm{\tau}^2}\abs{{\rm d}\tau}=I\left(\frac{m}{N}\right)=0.
\end{align*}
Since $n\ge 2$, this proves Theorem~\ref{main4}.
\end{proof}

\begin{Example}
For Example~\ref{Kodaira--Thurston}, let us consider the compatible almost complex structure associated with
\[
Z:=
\begin{pmatrix}
0 & 0 \\
0 & x_1
\end{pmatrix}
+\sqrt{-1}
\begin{pmatrix}
\frac{1}{x_1^2+1} & 0 \\
0 & 1
\end{pmatrix}.
\]
The corresponding $\Omega$ is
\begin{equation*}
\Omega_x=
\begin{pmatrix}
\sqrt{-1} & 0 \\
0 & x_1+\sqrt{-1}
\end{pmatrix}.
\end{equation*}
This $Z$ does not satisfies $(2)$ in Proposition~\ref{prop commutativity2}, nor
the condition~\eqref{commutativity} for any $m\in \Z^2$. In fact, for any $m\in \Z^2$, $((\partial_{x_1}\Omega)(m-Nx))_2=m_2-Nx_2$ while $((\partial_{x_2}\Omega)(m-Nx))_1=0$. In this case, \smash{${\tilde \vartheta}_{\frac{m}{N}}^t$} can be written as
\begin{align*}
{\tilde \vartheta}_{\frac{m}{N}}^t(x,y)={}&\sum_{\gamma \in \Z^2}g_\gamma \exp\biggl\{2\pi\sqrt{-1}N\Bigl[\frac{1}{2}\Bigl\{t\sqrt{-1}\Bigl(x_1-\gamma_1-\frac{m_1}{N}\Bigl)^2\\
& \hphantom{\sum_{\gamma \in \Z^2}g_\gamma}{}
+\Bigl(\frac{m_1}{N}+t\sqrt{-1}\Bigr)\Bigl(x_2-\gamma_2-\frac{m_2}{N}\Bigr)^2\Bigr\}\\
& \hphantom{\sum_{\gamma \in \Z^2}g_\gamma}{} +(x_2-\gamma_2)\Bigl\{ \frac{1}{2}\gamma_1(x_2+\gamma_2)-\Bigl( \frac{m_2}{N}+\gamma_2\Bigr)\gamma_2\Bigr\}\Bigr]\biggr\}{\rm e}^{2\pi\sqrt{-1}(m+N\gamma)\cdot y}.
\end{align*}
\end{Example}

\begin{Example}
In the case, where $n=2$ of Example~\ref{Jordan block2}, we can take the compatible almost complex structure associated with
\[
Z:=\frac{1}{x_2^2+1}
\begin{pmatrix}
\lambda^2x_2^3 & \lambda x_2^2 \\
\lambda x_2^2 & x_2
\end{pmatrix}
+\dfrac{\sqrt{-1}}{x_2^2+1}
\begin{pmatrix}
\bigl(1+\lambda^2\bigr)x_2^2+1 & \lambda x_2 \\
\lambda x_2 & 1
\end{pmatrix}.
\]
The corresponding $\Omega$ is
\begin{equation*}
\Omega_x=
\begin{pmatrix}
\sqrt{-1} & -\sqrt{-1}\lambda x_2 \\
-\sqrt{-1}\lambda x_2 & x_2+\sqrt{-1}\bigl(\lambda^2x_2^2+1\bigr)
\end{pmatrix}.
\end{equation*}
In this case, $\partial_{x_2}\Omega_{12}=-\sqrt{-1}\lambda$ and $\partial_{x_1}\Omega_{22}=0$. So, $Z$ satisfies $(2)$ in Proposition~\ref{prop commutativity2} if and only if $\lambda=0$, which is the special case of Example~\ref{EX4.8-2}. Equivalently, $Z$
does not satisfy the condition~\eqref{commutativity} for any $m\in \Z^2$ unless $\lambda=0$. In fact, for any $m\in \Z^2$, $((\partial_{x_1}\Omega)(m-Nx))_2=0$ while $((\partial_{x_2}\Omega)(m-Nx))_1=-\sqrt{-1}\lambda (m_2-Nx_2)$. In this case, \smash{${\tilde \vartheta}_{\frac{m}{N}}^t$ }can be written as
\begin{align*}
{\tilde \vartheta}_{\frac{m}{N}}^t(x,y)=\sum_{\gamma \in \Gamma}g_\gamma {\rm e}^{2\pi\sqrt{-1}\Theta\left(\frac{m}{N}, \gamma, x\right)}{\rm e}^{2\pi\sqrt{-1}\left\{(m_1+\gamma_2\lambda m_2+N\gamma_1)y_1+(m_2+N\gamma_2)y_2\right\}},
\end{align*}
where
\begin{align*}
\Theta\left( \frac{m}{N},\gamma ,x\right)={}&
N\left[
\frac{t\sqrt{-1}}{2}\left\{x_1-\gamma_1-\gamma_2\lambda (x_2-\gamma_2)-\frac{m_1}{N}\right\}^2\right.\\
& \left. -t\sqrt{-1}\lambda \frac{m_2}{N}\left\{x_1-\gamma_1-\gamma_2\lambda (x_2-\gamma_2)-\frac{m_1}{N}\right\} \left( x_2-\gamma_2-\frac{m_2}{N}\right)\right.\\
& \left. +\frac{1}{2}\left\{ \frac{m_2}{N}+t\sqrt{-1}\left( \lambda^2\frac{m_2^2}{N^2}+1\right)\right\} \left( x_2-\gamma_2-\frac{m_2}{N}\right)^2\right.\\
&\left.+\frac{1}{2}\gamma_2 (x_2-\gamma_2)(x_2+\gamma_2)-\left(\frac{m_2}{N}+\gamma_2\right) \gamma_2(x_2-\gamma_2)\right].
\end{align*}
\end{Example}

\appendix
\section{Fourier series}\label{appendixA}
Let $\pi\colon (M,\omega)\to B$ be a Lagrangian fibration on a complete $n$-dimensional $B$ with prequantum line bundle $\bigl(L,\nabla^L\bigr)\to (M,\omega)$. Then, it is obtained as the quotient of the $\pi_1(B)$-action on~${\pi_0\colon (M_0,\omega_0)\to \R^n}$ with prequantum line bundle $\bigl({L_0}, \nabla^{L_0}\bigr)\to (M_0,\omega_0)$. We take and fix a~fundamental domain $F$ of the $\pi_1(B)$-action on $\R^n$ as before.
Let $N\in \N$ be a positive integer and $s$ a smooth section of $L^{\otimes N}$. We identify $s$ with a $\pi_1(B)$-equivariant section of ${L_0}^{\otimes N}$. Then, for each $x\in \R^n$, $s(x, \cdot )$ can be expressed as the Fourier series
\begin{equation}\label{Fourier_again}
s(x,y)=\sum_{m\in \Z^n}a_m(x){\rm e}^{2\pi\sqrt{-1}m\cdot y}
\end{equation}
in $L^2\left({L_0}^{\otimes N}|_{\{x\}\times T^n}\right)$, where
\[
a_m(x):=\int_{T^n}s(x,y){\rm e}^{-2\pi\sqrt{-1}m\cdot y}{\rm d}y
\]
for $m\in \Z^n$. Then, we have the following lemma.
\begin{Lemma}\label{Fourier for derivatives}
For $j=1,\ldots , n$, the partial derivatives $\partial_{x_j}s$ and $\partial_{y_j}s$ have the following Fourier series:
\begin{align*}
\partial_{x_j}s(x,y)&=\sum_{m\in \Z^n}\partial_{x_j}a_m(x){\rm e}^{2\pi\sqrt{-1}m\cdot y},\\
\partial_{y_j}s(x,y)&=\sum_{m\in \Z^n}2\pi\sqrt{-1}m_ja_m(x){\rm e}^{2\pi\sqrt{-1}m\cdot y} 
\end{align*}
in $L^2\left({L_0}^{\otimes N}|_{\{x\}\times T^n}\right)$.
\end{Lemma}
\begin{proof}
Suppose that $\partial_{x_j}s$ has the following Fourier series with respect to $y_i$'s:
\begin{equation*}
\partial_{x_j}s(x,y)=\sum_{m\in \Z^n}b_m(x){\rm e}^{2\pi\sqrt{-1}m\cdot y}.
\end{equation*}
Then, $b_m(x)$ is computed by
\begin{align*}
b_m(x)&:=\int_{T^n}\partial_{x_j}s(x,y){\rm e}^{-2\pi\sqrt{-1}m\cdot y}{\rm d}y =\partial_{x_j}\int_{T^n}s(x,y){\rm e}^{-2\pi\sqrt{-1}m\cdot y}{\rm d}y=\partial_{x_j}a_m(x).
\end{align*}
This proves the first equality. Suppose that $\partial_{y_j}s$ has the following Fourier series with respect to~$y_i$'s:
\begin{equation*}
\partial_{y_j}s(x,y)=\sum_{m\in \Z^n}c_m(x){\rm e}^{2\pi\sqrt{-1}m\cdot y}.
\end{equation*}
Then, $c_m(x)$ is computed by
\begin{align*}
c_m(x):={}&\int_{T^n}\partial_{y_j}s(x,y){\rm e}^{-2\pi\sqrt{-1}m\cdot y}{\rm d}y\\
={}&\int_{T^{n-1}}
 \left(\int_{S^1}\partial_{y_j}s(x,y){\rm e}^{-2\pi\sqrt{-1}m\cdot y}{\rm d}y_j\right)
{\rm d}y_1\cdots \hat{{\rm d}y_j}\cdots {\rm d}y_n \\
={}& \int_{T^{n-1}} \left(\big[s(x,y){\rm e}^{-2\pi\sqrt{-1}m\cdot y}\big]_{y_j=0}^{y_j=1}\right.\\
&\left.- \int_{S^1}s(x,y)\partial_{y_j}{\rm e}^{-2\pi\sqrt{-1}m\cdot y}{\rm d}y_j\right){\rm d}y_1\cdots \hat{{\rm d}y_j}\cdots {\rm d}y_n \\
={}&\int_{T^{n-1}} \left(2\pi\sqrt{-1}m_j\int_{S^1}s(x,y){\rm e}^{-2\pi\sqrt{-1}m\cdot y}{\rm d}y_j\right){\rm d}y_1\cdots \hat{{\rm d}y_j}\cdots {\rm d}y_n \\
={}&2\pi\sqrt{-1}m_j\int_{T^n}s(x,y){\rm e}^{-2\pi\sqrt{-1}m\cdot y}{\rm d}y \\
={}& 2\pi\sqrt{-1}m_ja_m(x).
\end{align*}
This proves the second equality.
\end{proof}

\begin{Lemma}\label{convergence in F times T}
If $s$ is in $L^2\bigl(L^{\otimes N}\bigr)$, then the following formulae hold:
\begin{gather}
 \lim_{k\to \infty}\bigg\lVert s-\sum_{\abs{m}\le k}a_m(x){\rm e}^{2\pi\sqrt{-1}m\cdot y}\bigg\rVert_{L^2({L_0}^{\otimes N}|_{F\times T^n})}=0, \label{Fourier series for s}\nonumber\\
\norm{s}_{L^2(L^{\otimes N})}^2=\sum_{m\in \Z^n}\norm{a_m(x)}_{L^2(F)}^2,\label{norm of s}
\end{gather}
where $\abs{m}:=m_1+\cdots +m_n$ for $m=(m_1,\ldots ,m_n)\in \Z^n$. Namely, the right-hand side of \eqref{Fourier_again} also converges to $s$ in $L^2\bigl({L_0}^{\otimes N}|_{F\times T^n}\bigr)$.
\end{Lemma}
\begin{proof}
For each $x\in \R^n$, we define $s_x\in L^2\bigl({L_0}^{\otimes N}|_{\{x\}\times T^n}\bigr)$ by $s_x(y):=s(x,y)$. Then, \eqref{Fourier_again}~implies
\begin{gather}
\lim_{k\to \infty}\bigg\lVert s_x-\sum_{\abs{m}\le k}a_m(x){\rm e}^{2\pi\sqrt{-1}m\cdot y}\bigg\rVert_{L^2({L_0}^{\otimes N}|_{\{x\}\times T^n})}=0, \label{Fourier1}
\end{gather}
and $s_x$ satisfies
\begin{gather}
\norm{s_x}_{L^2({L_0}^{\otimes N}|_{\{x\}\times T^n})}^2=\sum_{m\in \Z^n}\abs{a_m(x)}^2.\label{Fourier2}
\end{gather}
By using \eqref{Fourier2} and the monotone convergence theorem, \eqref{norm of s} can be obtained as follows:
\begin{align*}
\norm{s}_{L^2(L^{\otimes N})}^2&=\int_M\langle s, s\rangle_{L^{\otimes N}} (-1)^{\frac{n(n-1)}{2}}\dfrac{(N\omega)^n}{n!}=\int_B\pi_*\left(\langle s, s\rangle_{L^{\otimes N}} (-1)^{\frac{n(n-1)}{2}}\dfrac{(N\omega)^n}{n!}\right) \\
&=\int_F\norm{s_x}_{L^2({L_0}^{\otimes N}|_{\{x\}\times T^n})}^2\abs{{\rm d}x} =\int_F\sum_{m\in \Z^n}\abs{a_m(x)}^2\abs{{\rm d}x} \\
&=\sum_{m\in \Z^n}\int_F\abs{a_m(x)}^2\abs{{\rm d}x} =\sum_{m\in \Z^n}\norm{a_m(x)}_{L^2(F)}^2.
\end{align*}
Next, let us prove that
\[\sum_{m\in \Z^n}a_m(x){\rm e}^{2\pi\sqrt{-1}m\cdot y}
\]
 converges with respect to the norm of $L^2\bigl({L_0}^{\otimes N}|_{F\times T^n}\bigr)$. For each $k\in \N$, we put
 \[s_k(x,y):=\sum_{\abs{m}\le k}a_m(x){\rm e}^{2\pi\sqrt{-1}m\cdot y}.
 \]
 To prove it, it is sufficient to show $\{s_k\}_{k\in \N}$ is a Cauchy sequence in $L^2\bigl({L_0}^{\otimes N}|_{F\times T^n}\bigr)$. For $k<l$ in $\N$,
\begin{align}
\norm{s_l-s_k}_{L^2({L_0}^{\otimes N}|_{F\times T^n})}^2&=\bigg\lVert\sum_{\abs{m}\le l}a_m(x){\rm e}^{2\pi\sqrt{-1}m\cdot y}-\sum_{\abs{m}\le k}a_m(x){\rm e}^{2\pi\sqrt{-1}m\cdot y}\bigg\rVert_{L^2({L_0}^{\otimes N}|_{F\times T^n})}^2 \nonumber\\
&=\bigg\lVert\sum_{k<\abs{m}\le l}a_m(x){\rm e}^{2\pi\sqrt{-1}m\cdot y}\bigg\rVert_{L^2({L_0}^{\otimes N}|_{F\times T^n})}^2\nonumber\\
&=\int_F\sum_{k<\abs{m}\le l}\abs{a_m(x)}^2\abs{{\rm d}x} \nonumber\\
&=\bigg\lvert \sum_{\abs{m}\le k}\norm{a_m(x)}_{L^2(F)}^2-\sum_{\abs{m}\le l}\norm{a_m(x)}_{L^2(F)}^2\bigg\rvert.\label{Cauchy for sN}
\end{align}
Since $s$ is square integrable, as we showed above, \smash{$\sum_{m\in \Z^n}\norm{a_m(x)}_{L^2(F)}^2$} converges to \smash{$\norm{s}_{L^2(L^{\otimes N})}$}. In particular, the sequence \smash{$\big\{\sum_{\abs{m}\le k}\norm{a_m(x)}_{L^2(F)}^2\big\}_{k\in \N}$} is a Cauchy sequence. Thus, by \eqref{Cauchy for sN}
\[
\lim_{k\to \infty, l\to \infty}\norm{s_l-s_k}_{L^2({L_0}^{\otimes N}|_{F\times T^n})}=0.
\]
Let ${\tilde s}\in L^2\bigl({L_0}^{\otimes N}|_{F\times T^n}\bigr)$ be the limit of $\{s_k\}_{k\in \N}$. Then, $\{s_k\}_{k\in \N}$ pointwise converges to ${\tilde s}$. But, by \eqref{Fourier1}, $\{s_k\}_{k\in \N}$ also pointwise converges to $s$. This implies ${\tilde s}=s$.
\end{proof}

\begin{Remark}
By the continuity of the inner product of the Hilbert space, Lemma~\ref{convergence in F times T} enable us to interchange operations to take limits and integrals for $L^2$-sections on $L^{\otimes N}$.
\end{Remark}

\section{Proof of Proposition~\ref{prop commutativity2}}\label{appendixB}

If all $m\in \Z^n$ are integrable, then by putting $m=0$, we have $ ( (\partial_{x_i}\Omega )_x x )_j= ( (\partial_{x_j}\Omega )_x x )_i$. By substituting this to \eqref{commutativity}, we can see the condition $ ( (\partial_{x_i}\Omega )_x m )_j= ( (\partial_{x_j}\Omega )_x m )_i$ holds for all~${m\in \Z^n}$. In particular, by substituting $m=e_k$ to this condition for each $k=1,\ldots , n$, we can obtain $(2)$. $(2) \Rightarrow (1)$ is trivial.

We show $(2)\Leftrightarrow (3)$. $(2)$ is equivalent to the following two conditions:
\begin{align}
& \bigl(\bigl(Y+XY^{-1}X\bigr)^{-1}\partial_{x_i}\bigl(XY^{-1}\bigr)\bigr)_{jk}=\bigl(\bigl(Y+XY^{-1}X\bigr)^{-1}\partial_{x_j}\bigl(XY^{-1}\bigr)\bigr)_{ik} \label{ReOmega}\\
& \partial_{x_i}\bigl(Y+XY^{-1}X\bigr)^{-1}_{jk}=\partial_{x_j}\bigl(Y+XY^{-1}X\bigr)^{-1}_{ik}\label{ImOmega}
\end{align}
for $i, j, k=1,\ldots , n$. For $i=1, \ldots ,2n$, we set
\[
\Gamma_i:=
\begin{pmatrix}
\Gamma^1_{i 1} & \cdots & \Gamma^1_{i 2n} \\
\vdots & & \vdots \\
\Gamma^{2n}_{i 1} & \cdots & \Gamma^{2n}_{i 2n}
\end{pmatrix},
\]
where $\Gamma_{ij}^k$ is the Christoffel symbol. Then, $(3)$ is equivalent to
\[
0=\partial_iJ+\Gamma_iJ-J\Gamma_i, \qquad i=1,\ldots ,2n,
\]
where
\[
\partial_i=
\begin{cases}
\partial_{x_i}, & i=1,\ldots ,n, \\
\partial_{y_{i-n}}, & i=n+1,\ldots , 2n.
\end{cases}
\]
It is also equivalent to the following conditions:
\begin{align}
& XY^{-1}
\begin{pmatrix}
\partial_{x_1}\bigl(XY^{-1}\bigr)_{1i} & \cdots & \partial_{x_n}\bigl(XY^{-1}\bigr)_{1i} \\
\vdots & & \vdots \\
\partial_{x_1}\bigl(XY^{-1}\bigr)_{ni} & \cdots & \partial_{x_n}\bigl(XY^{-1}\bigr)_{ni}
\end{pmatrix} \nonumber \\
&\qquad-\bigl(Y+XY^{-1}X\bigr)
\begin{pmatrix}
\partial_{x_1}\bigl(Y^{-1}\bigr)_{1i}-\partial_{x_1}\bigl(Y^{-1}\bigr)_{1i} & \cdots & \partial_{x_n}\bigl(Y^{-1}\bigr)_{1i}-\partial_{x_1}\bigl(Y^{-1}\bigr)_{ni} \\
\vdots & & \vdots \\
\partial_{x_1}\bigl(Y^{-1}\bigr)_{ni}-\partial_{x_n}\bigl(Y^{-1}\bigr)_{1i} & \cdots & \partial_{x_n}\bigl(Y^{-1}\bigr)_{ni}-\partial_{x_n}\bigl(Y^{-1}\bigr)_{ni}
\end{pmatrix}\nonumber \\
&\phantom{\qquad-}{}=
\begin{pmatrix}
\partial_{x_1}\bigl(XY^{-1}\bigr)_{1i} & \cdots & \partial_{x_n}\bigl(XY^{-1}\bigr)_{1i} \\
\vdots & & \vdots \\
\partial_{x_1}\bigl(XY^{-1}\bigr)_{ni} & \cdots & \partial_{x_n}\bigl(XY^{-1}\bigr)_{ni}
\end{pmatrix}XY^{-1},\label{K1}
\\
& Y^{-1}
\begin{pmatrix}
\partial_{x_1}\bigl(XY^{-1}\bigr)_{1i} & \cdots & \partial_{x_n}\bigl(XY^{-1}\bigr)_{1i} \\
\vdots & & \vdots \\
\partial_{x_1}\bigl(XY^{-1}\bigr)_{ni} & \cdots & \partial_{x_n}\bigl(XY^{-1}\bigr)_{ni}
\end{pmatrix}\nonumber \\
&\qquad-Y^{-1}X
\begin{pmatrix}
\partial_{x_1}\bigl(Y^{-1}\bigr)_{1i}-\partial_{x_1}\bigl(Y^{-1}\bigr)_{1i} & \cdots & \partial_{x_n}\bigl(Y^{-1}\bigr)_{1i}-\partial_{x_1}\bigl(Y^{-1}\bigr)_{ni} \\
\vdots & & \vdots \\
\partial_{x_1}\bigl(Y^{-1}\bigr)_{ni}-\partial_{x_n}\bigl(Y^{-1}\bigr)_{1i} & \cdots & \partial_{x_n}\bigl(Y^{-1}\bigr)_{ni}-\partial_{x_n}\bigl(Y^{-1}\bigr)_{ni}
\end{pmatrix}\nonumber \\
&\phantom{\qquad-}{}=
\begin{pmatrix}
\partial_{x_1}\bigl(Y^{-1}\bigr)_{1i}-\partial_{x_1}\bigl(Y^{-1}\bigr)_{1i} & \cdots & \partial_{x_n}\bigl(Y^{-1}\bigr)_{1i}-\partial_{x_1}\bigl(Y^{-1}\bigr)_{ni} \\
\vdots & & \vdots \\
\partial_{x_1}\bigl(Y^{-1}\bigr)_{ni}-\partial_{x_n}\bigl(Y^{-1}\bigr)_{1i} & \cdots & \partial_{x_n}\bigl(Y^{-1}\bigr)_{ni}-\partial_{x_n}\bigl(Y^{-1}\bigr)_{ni}
\end{pmatrix}
XY^{-1}\nonumber \\
&\phantom{\qquad-=}{}+
\begin{pmatrix}
\partial_{x_1}\bigl(Y^{-1}X\bigr)_{i1} & \cdots & \partial_{x_1}\bigl(Y^{-1}X\bigr)_{in} \\
\vdots & & \vdots \\
\partial_{x_n}\bigl(Y^{-1}X\bigr)_{i1} & \cdots & \partial_{x_n}\bigl(Y^{-1}X\bigr)_{in}
\end{pmatrix}
Y^{-1},\label{K2}
\\
& \bigl(Y+XY^{-1}X\bigr)
\begin{pmatrix}
\partial_{x_1}\bigl(Y^{-1}X\bigr)_{i1} & \cdots & \partial_{x_1}\bigl(Y^{-1}X\bigr)_{in} \\
\vdots & & \vdots \\
\partial_{x_n}\bigl(Y^{-1}X\bigr)_{i1} & \cdots & \partial_{x_n}\bigl(Y^{-1}X\bigr)_{in}
\end{pmatrix}\nonumber \\
&\qquad=
\begin{pmatrix}
\partial_{x_1}\bigl(XY^{-1}\bigr)_{1i} & \cdots & \partial_{x_n}\bigl(XY^{-1}\bigr)_{1i} \\
\vdots & & \vdots \\
\partial_{x_1}\bigl(XY^{-1}\bigr)_{ni} & \cdots & \partial_{x_n}\bigl(XY^{-1}\bigr)_{ni}
\end{pmatrix}
\bigl(Y+XY^{-1}X\bigr), \label{K3}
\\
& Y^{-1}X
\begin{pmatrix}
\partial_{x_1}\bigl(Y^{-1}X\bigr)_{i1} & \cdots & \partial_{x_1}\bigl(Y^{-1}X\bigr)_{in} \\
\vdots & & \vdots \\
\partial_{x_n}\bigl(Y^{-1}X\bigr)_{i1} & \cdots & \partial_{x_n}\bigl(Y^{-1}X\bigr)_{in}
\end{pmatrix}\nonumber \\
&\qquad=
\begin{pmatrix}
\partial_{x_1}\bigl(Y^{-1}\bigr)_{1i}-\partial_{x_1}\bigl(Y^{-1}\bigr)_{1i} & \cdots & \partial_{x_n}\bigl(Y^{-1}\bigr)_{1i}-\partial_{x_1}\bigl(Y^{-1}\bigr)_{ni} \\
\vdots & & \vdots \\
\partial_{x_1}\bigl(Y^{-1}\bigr)_{ni}-\partial_{x_n}\bigl(Y^{-1}\bigr)_{1i} & \cdots & \partial_{x_n}\bigl(Y^{-1}\bigr)_{ni}-\partial_{x_n}\bigl(Y^{-1}\bigr)_{ni}
\end{pmatrix}
(Y+XY^{-1})\nonumber \\
&\phantom{\qquad=}{}+
\begin{pmatrix}
\partial_{x_1}\bigl(Y^{-1}X\bigr)_{i1} & \cdots & \partial_{x_1}\bigl(Y^{-1}X\bigr)_{in} \\
\vdots & & \vdots \\
\partial_{x_n}\bigl(Y^{-1}X\bigr)_{i1} & \cdots & \partial_{x_n}\bigl(Y^{-1}X\bigr)_{in}
\end{pmatrix}
Y^{-1}X, \label{K4}
\\
& XY^{-1}
\begin{pmatrix}
\partial_{x_1}\bigl(Y+XY^{-1}X\bigr)_{1i} & \cdots & \partial_{x_n}\bigl(Y+XY^{-1}X\bigr)_{1i} \\
\vdots & & \vdots \\
\partial_{x_1}\bigl(Y+XY^{-1}X\bigr)_{ni} & \cdots & \partial_{x_n}\bigl(Y+XY^{-1}X\bigr)_{ni}
\end{pmatrix}+\bigl(Y+XY^{-1}X\bigr)\nonumber \\
&\qquad\times
\begin{pmatrix}
-\partial_{x_1}\bigl(Y^{-1}X\bigr)_{1i}+\partial_{x_1}\bigl(XY^{-1}\bigr)_{i1} & {\hspace*{-3mm}}\cdots & {\hspace*{-3mm}}-\partial_{x_n}\bigl(Y^{-1}X\bigr)_{1i}+\partial_{x_1}\bigl(XY^{-1}\bigr)_{in} \\
\vdots & & \vdots \\
-\partial_{x_1}\bigl(Y^{-1}X\bigr)_{ni}+\partial_{x_n}\bigl(XY^{-1}\bigr)_{i1} & {\hspace*{-3mm}}\cdots & {\hspace*{-3mm}}-\partial_{x_n}\bigl(Y^{-1}X\bigr)_{ni}+\partial_{x_n}\bigl(XY^{-1}\bigr)_{in}
\end{pmatrix}\nonumber \\
&\phantom{\qquad+}{}=\begin{pmatrix}
\partial_{x_1}\bigl(Y+XY^{-1}X\bigr)_{1i} & \cdots & \partial_{x_n}\bigl(Y+XY^{-1}X\bigr)_{1i} \\
\vdots & & \vdots \\
\partial_{x_1}\bigl(Y+XY^{-1}X\bigr)_{ni} & \cdots & \partial_{x_n}\bigl(Y+XY^{-1}X\bigr)_{ni}
\end{pmatrix}
XY^{-1},\label{K5}
\\
&Y^{-1}
\begin{pmatrix}
\partial_{x_1}\bigl(Y+XY^{-1}X\bigr)_{1i} & \cdots & \partial_{x_n}\bigl(Y+XY^{-1}X\bigr)_{1i} \\
\vdots & & \vdots \\
\partial_{x_1}\bigl(Y+XY^{-1}X\bigr)_{ni} & \cdots & \partial_{x_n}\bigl(Y+XY^{-1}X\bigr)_{ni}
\end{pmatrix}\nonumber \\
&\qquad+Y^{-1}X
\begin{pmatrix}
-\partial_{x_1}\bigl(Y^{-1}X\bigr)_{1i}+\partial_{x_1}\bigl(XY^{-1}\bigr)_{i1} & {\hspace*{-3mm}}\cdots{\hspace*{-3mm}} & -\partial_{x_n}\bigl(Y^{-1}X\bigr)_{1i}+\partial_{x_1}\bigl(XY^{-1}\bigr)_{in} \\
\vdots & & \vdots \\
-\partial_{x_1}\bigl(Y^{-1}X\bigr)_{ni}+\partial_{x_n}\bigl(XY^{-1}\bigr)_{i1} &{\hspace*{-3mm}} \cdots {\hspace*{-3mm}}& -\partial_{x_n}\bigl(Y^{-1}X\bigr)_{ni}+\partial_{x_n}\bigl(XY^{-1}\bigr)_{in}
\end{pmatrix}\nonumber \\
&\phantom{\qquad+}{}=-
\begin{pmatrix}
-\partial_{x_1}\bigl(Y^{-1}X\bigr)_{1i}+\partial_{x_1}\bigl(XY^{-1}\bigr)_{i1} &{\hspace*{-3mm}} \cdots {\hspace*{-3mm}}& -\partial_{x_n}\bigl(Y^{-1}X\bigr)_{1i}+\partial_{x_1}\bigl(XY^{-1}\bigr)_{in} \\
\vdots & & \vdots \\
-\partial_{x_1}\bigl(Y^{-1}X\bigr)_{ni}+\partial_{x_n}\bigl(XY^{-1}\bigr)_{i1} & {\hspace*{-3mm}}\cdots{\hspace*{-3mm}} & -\partial_{x_n}\bigl(Y^{-1}X\bigr)_{ni}+\partial_{x_n}\bigl(XY^{-1}\bigr)_{in}
\end{pmatrix}XY^{-1}\nonumber \\
&\phantom{\qquad+=}{}+
\begin{pmatrix}
\partial_{x_1}\bigl(Y+XY^{-1}X\bigr)_{i1} & \cdots & \partial_{x_1}\bigl(Y+XY^{-1}X\bigr)_{in} \\
\vdots & & \vdots \\
\partial_{x_n}\bigl(Y+XY^{-1}X\bigr)_{i1} & \cdots & \partial_{x_n}\bigl(Y+XY^{-1}X\bigr)_{in}
\end{pmatrix}
Y^{-1},\label{K6}
\\
& \bigl(Y+XY^{-1}X\bigr)
\begin{pmatrix}
\partial_{x_1}\bigl(Y+XY^{-1}X\bigr)_{i1} & \cdots & \partial_{x_1}\bigl(Y+XY^{-1}X\bigr)_{in} \\
\vdots & & \vdots \\
\partial_{x_n}\bigl(Y+XY^{-1}X\bigr)_{i1} & \cdots & \partial_{x_n}\bigl(Y+XY^{-1}X\bigr)_{in}
\end{pmatrix}\nonumber \\
&\qquad=
\begin{pmatrix}
\partial_{x_1}\bigl(Y+XY^{-1}X\bigr)_{1i} & \cdots & \partial_{x_n}\bigl(Y+XY^{-1}X\bigr)_{1i} \\
\vdots & & \vdots \\
\partial_{x_1}\bigl(Y+XY^{-1}X\bigr)_{ni} & \cdots & \partial_{x_n}\bigl(Y+XY^{-1}X\bigr)_{ni}
\end{pmatrix}
\bigl(Y+XY^{-1}X\bigr),\label{K7}
\\
&Y^{-1}X
\begin{pmatrix}
\partial_{x_1}\bigl(Y+XY^{-1}X\bigr)_{i1} & \cdots & \partial_{x_1}\bigl(Y+XY^{-1}X\bigr)_{in} \\
\vdots & & \vdots \\
\partial_{x_n}\bigl(Y+XY^{-1}X\bigr)_{i1} & \cdots & \partial_{x_n}\bigl(Y+XY^{-1}X\bigr)_{in}
\end{pmatrix}\nonumber\\
&\qquad= -
\begin{pmatrix}
 - \partial_{x_1}\bigl(Y^{-1}X\bigr)_{1i} + \partial_{x_1}\bigl(XY^{-1}\bigr)_{i1} & {\hspace*{-3mm}}\cdots & {\hspace*{-3mm}}- \partial_{x_n}\bigl(Y^{-1}X\bigr)_{1i} + \partial_{x_1}\bigl(XY^{-1}\bigr)_{in} \\
\vdots & & \vdots \\
 - \partial_{x_1}\bigl(Y^{-1}X\bigr)_{ni} + \partial_{x_n}\bigl(XY^{-1}\bigr)_{i1} & {\hspace*{-3mm}}\cdots & {\hspace*{-3mm}}- \partial_{x_n}\bigl(Y^{-1}X\bigr)_{ni} + \partial_{x_n}\bigl(XY^{-1}\bigr)_{in}
\end{pmatrix}
\nonumber \\
&\phantom{\qquad=}{}\times \bigl(Y + XY^{-1}X\bigr)\nonumber\\
&\phantom{\qquad=\times}{}+
\begin{pmatrix}
\partial_{x_1}\bigl(Y+XY^{-1}X\bigr)_{i1} & \cdots & \partial_{x_1}\bigl(Y+XY^{-1}X\bigr)_{in} \\
\vdots & & \vdots \\
\partial_{x_n}\bigl(Y+XY^{-1}X\bigr)_{i1} & \cdots & \partial_{x_n}\bigl(Y+XY^{-1}X\bigr)_{in}
\end{pmatrix}
Y^{-1}X.\label{K8}
\end{align}
for $i=1,\ldots ,n$. It is easy to see that \eqref{K4} and \eqref{K8} are obtained by transposing \eqref{K1} and~\eqref{K5}, respectively. First, we show that \eqref{ReOmega} is equivalent to \eqref{K3}. In fact, \eqref{ReOmega} implies
\[
\begin{pmatrix}
\partial_{x_1}\bigl(XY^{-1}\bigr)_{1k} & \cdots & \partial_{x_1}\bigl(XY^{-1}\bigr)_{nk} \\
\vdots & & \vdots \\
\partial_{x_n}\bigl(XY^{-1}\bigr)_{1k} & \cdots & \partial_{x_n}\bigl(XY^{-1}\bigr)_{nk}
\end{pmatrix}
\bigl(Y+XY^{-1}X\bigr)^{-1}
\]
is symmetric for $k=1,\ldots , n$. Since $X$, $Y$ is symmetric, this implies \eqref{K3}. Next, we show \eqref{K7} is equivalent to \eqref{ImOmega}. \eqref{K7} is equivalent to
\begin{align}
&\begin{pmatrix}
\partial_{x_1}\bigl(Y+XY^{-1}X\bigr)_{i1} & \cdots & \partial_{x_1}\bigl(Y+XY^{-1}X\bigr)_{in} \\
\vdots & & \vdots \\
\partial_{x_n}\bigl(Y+XY^{-1}X\bigr)_{i1} & \cdots & \partial_{x_n}\bigl(Y+XY^{-1}X\bigr)_{in}
\end{pmatrix}
\bigl(Y+XY^{-1}X\bigr)^{-1}\nonumber\\
&\qquad=\bigl(Y+XY^{-1}X\bigr)^{-1}
\begin{pmatrix}
\partial_{x_1}\bigl(Y+XY^{-1}X\bigr)_{1i} & \cdots & \partial_{x_n}\bigl(Y+XY^{-1}X\bigr)_{1i} \\
\vdots & & \vdots \\
\partial_{x_1}\bigl(Y+XY^{-1}X\bigr)_{ni} & \cdots & \partial_{x_n}\bigl(Y+XY^{-1}X\bigr)_{ni}
\end{pmatrix}.\label{K7'}
\end{align}
By computing the $(j, k)$-components of the both sides of \eqref{K7'}, we obtain
\begin{gather*}
\sum_{l=1}^n\bigl(\partial_{x_j}\bigl(Y+XY^{-1}X\bigr)^{-1}_{kl}\bigr)\bigl(Y+XY^{-1}X\bigr)_{li}=\sum_{l=1}^n\bigl(\partial_{x_k}\bigl(Y+XY^{-1}X\bigr)^{-1}_{jl}\bigr)\bigl(Y+XY^{-1}X\bigr)_{li}
\end{gather*}
for $i, j, k=1,\ldots , n$. Here, we used
\begin{align*}
0={}&\partial_{x_j}\bigl(\bigl(Y+XY^{-1}X\bigr)\bigl(Y+XY^{-1}X\bigr)^{-1}\bigr) \\
={}&\bigl(\partial_{x_j}\bigl(Y+XY^{-1}X\bigr)\bigr)\bigl(Y+XY^{-1}X\bigr)^{-1}+\bigl(Y+XY^{-1}X\bigr)\partial_{x_j}\bigl(Y+XY^{-1}X\bigr)^{-1}
\end{align*}
and so on. Thus,
\begin{align*}
\partial_{x_j}\bigl(Y+XY^{-1}X\bigr)^{-1}_{km}&=\sum_{i=1}^n\sum_{l=1}^n\partial_{x_j}\bigl(\bigl(Y+XY^{-1}X\bigr)^{-1}_{kl}\bigr)\bigl(Y+XY^{-1}X\bigr)_{li}\bigl(Y+XY^{-1}X\bigr)^{-1}_{im} \\
&=\sum_{i=1}^n\sum_{l=1}^n\bigl(\partial_{x_k}\bigl(Y+XY^{-1}X\bigr)^{-1}_{jl}\bigr)\bigl(Y+XY^{-1}X\bigr)_{li}\bigl(Y+XY^{-1}X\bigr)^{-1}_{im} \\
&=\partial_{x_k}\bigl(Y+XY^{-1}X\bigr)^{-1}_{jm}.
\end{align*}
This implies \eqref{ImOmega}. In particular, this means $(3)\Rightarrow (2)$.

We show \eqref{K1}, \eqref{K2}, \eqref{K5}, and \eqref{K6} are obtained from $(2)$. 
To show \eqref{K5}, it is sufficient to show
\begin{align}
0={}& \bigl(Y+XY^{-1}X\bigr)^{-1}XY^{-1}
\begin{pmatrix}
\partial_{x_1}\bigl(Y+XY^{-1}X\bigr)_{1i} & \cdots & \partial_{x_n}\bigl(Y+XY^{-1}X\bigr)_{1i} \\
\vdots & & \vdots \\
\partial_{x_1}\bigl(Y+XY^{-1}X\bigr)_{ni} & \cdots & \partial_{x_n}\bigl(Y+XY^{-1}X\bigr)_{ni}
\end{pmatrix} \nonumber \\
&-
\begin{pmatrix}
\partial_{x_1}\bigl(Y^{-1}X\bigr)_{1i} & \cdots & \partial_{x_n}\bigl(Y^{-1}X\bigr)_{1i} \\
\vdots & & \vdots \\
\partial_{x_1}\bigl(Y^{-1}X\bigr)_{ni} & \cdots & \partial_{x_n}\bigl(Y^{-1}X\bigr)_{ni}
\end{pmatrix} \nonumber \\
&-\bigl(Y+XY^{-1}X\bigr)^{-1}
\begin{pmatrix}
\partial_{x_1}\bigl(Y+XY^{-1}X\bigr)_{1i} & \cdots & \partial_{x_n}\bigl(Y+XY^{-1}X\bigr)_{1i} \\
\vdots & & \vdots \\
\partial_{x_1}\bigl(Y+XY^{-1}X\bigr)_{ni} & \cdots & \partial_{x_n}\bigl(Y+XY^{-1}X\bigr)_{ni}
\end{pmatrix}
XY^{-1}\nonumber \\
&+
\begin{pmatrix}
\partial_{x_1}\bigl(XY^{-1}\bigr)_{i1} & \cdots & \partial_{x_1}\bigl(XY^{-1}\bigr)_{in} \\
\vdots & & \vdots \\
\partial_{x_n}\bigl(XY^{-1}\bigr)_{i1} & \cdots & \partial_{x_n}\bigl(XY^{-1}\bigr)_{in}
\end{pmatrix}.\label{K5'}
\end{align}
Since $\Omega$ is symmetric, so is its real part $\re \Omega=\bigl(Y+XY^{-1}X\bigr)^{-1}XY^{-1}$. By taking the real part of~$(2)$, we also have
\[
\partial_{x_i}\bigl(\bigl(Y+XY^{-1}X\bigr)^{-1}XY^{-1}\bigr)_{jk}=\partial_{x_j}\bigl(\bigl(Y+XY^{-1}X\bigr)^{-1}XY^{-1}\bigr)_{ik}.
\]
By using these as well as \eqref{ReOmega} and \eqref{ImOmega}, the $(j, k)$-component of the first two terms of the right-hand side of \eqref{K5'} can be computed as
\begin{gather*}
 \sum_{l}\bigl(\bigl(Y+XY^{-1}X\bigr)^{-1}XY^{-1}\bigr)_{jl}\partial_{x_k}\bigl(Y+XY^{-1}X\bigr)_{li}-\partial_{x_k}\bigl(Y^{-1}X\bigr)_{ji}\\
\qquad= \sum_{l}\bigl(Y^{-1}X\bigl(Y+XY^{-1}X\bigr)^{-1}\bigr)_{jl}\partial_{x_k}\bigl(Y+XY^{-1}X\bigr)_{li}-\partial_{x_k}\bigl(Y^{-1}X\bigr)_{ji} \\
\qquad= \partial_{x_k} \Big(\sum_{l}\bigl(Y^{-1}X\bigl(Y+XY^{-1}X\bigr)^{-1}\bigr)_{jl}\bigl(Y+XY^{-1}X\bigr)_{li}\Big)\\
\qquad\phantom{=}{} -\sum_{l}\bigl(\partial_{x_k} \bigl(Y^{-1}X\bigl(Y+XY^{-1}X\bigr)^{-1}\bigr)_{jl}\bigr)\bigl(Y+XY^{-1}X\bigr)_{li}-\partial_{x_k}\bigl(Y^{-1}X\bigr)_{ji} \\
\qquad= -\sum_{l}\bigl(\partial_{x_k} \bigl(Y^{-1}X\bigl(Y+XY^{-1}X\bigr)^{-1}\bigr)_{jl}\bigr)\bigl(Y+XY^{-1}X\bigr)_{li}\\
\qquad= -\sum_{l}\bigl(\partial_{x_j} \bigl(\bigl(Y+XY^{-1}X\bigr)^{-1}XY^{-1}\bigr)_{kl}\bigr)\bigl(Y+XY^{-1}X\bigr)_{li}.
\end{gather*}
On the other hand, the $(j, k)$-component of the last two terms of the right-hand side of \eqref{K5'} can be computed as
\begin{gather*}
-\sum_{m, l} \bigl(Y+XY^{-1}X\bigr)^{-1}_{jl}\bigl(\partial_{x_m}\bigl(Y+XY^{-1}X\bigr)_{li}\bigr)\bigl(XY^{-1}\bigr)_{mk}+\partial_{x_j}\bigl(XY^{-1}\bigr)_{ik}\\
\qquad=\sum_{m, l} \bigl(\partial_{x_m}\bigl(Y+XY^{-1}X\bigr)^{-1}_{jl}\bigr)\bigl(Y+XY^{-1}X\bigr)_{li}\bigl(XY^{-1}\bigr)_{mk}+\partial_{x_j}\bigl(XY^{-1}\bigr)_{ik}\\
\qquad=\sum_{m, l}\bigl(Y+XY^{-1}X\bigr)_{li}\bigl(\partial_{x_j}\bigl(Y+XY^{-1}X\bigr)^{-1}_{ml}\bigr) \bigl(XY^{-1}\bigr)_{mk}\\
\qquad\phantom{=}{} +\sum_{m, l}\bigl(Y+XY^{-1}X\bigr)_{li}\bigl(Y+XY^{-1}X\bigr)^{-1}_{ml}\partial_{x_j}\bigl(XY^{-1}\bigr)_{mk}
\\
\qquad=\sum_{l}\bigl(\partial_{x_j}\bigl(\bigl(Y+XY^{-1}X\bigr)^{-1}XY^{-1}\bigr)_{kl}\bigr) \bigl(Y+XY^{-1}X\bigr)_{li}.
\end{gather*}
This proves \eqref{K5'}. We show \eqref{K6}. We put
\[
W:=
\begin{pmatrix}
\partial_{x_1}\bigl(Y+XY^{-1}X\bigr)_{1i} & \cdots & \partial_{x_n}\bigl(Y+XY^{-1}X\bigr)_{1i} \\
\vdots & & \vdots \\
\partial_{x_1}\bigl(Y+XY^{-1}X\bigr)_{ni} & \cdots & \partial_{x_n}\bigl(Y+XY^{-1}X\bigr)_{ni}
\end{pmatrix}.
\]
By \eqref{K5} and \eqref{K7}, we obtain
\begin{align*}
&\begin{pmatrix}
-\partial_{x_1}\bigl(Y^{-1}X\bigr)_{1i}+\partial_{x_1}\bigl(XY^{-1}\bigr)_{i1} & \cdots & -\partial_{x_n}\bigl(Y^{-1}X\bigr)_{1i}+\partial_{x_1}\bigl(XY^{-1}\bigr)_{in} \\
\vdots & & \vdots \\
-\partial_{x_1}\bigl(Y^{-1}X\bigr)_{ni}+\partial_{x_n}\bigl(XY^{-1}\bigr)_{i1} & \cdots & -\partial_{x_n}\bigl(Y^{-1}X\bigr)_{ni}+\partial_{x_n}\bigl(XY^{-1}\bigr)_{in}
\end{pmatrix}\\
&\qquad=\bigl(Y+XY^{-1}X\bigr)^{-1}WXY^{-1}-\bigl(Y+XY^{-1}X\bigr)^{-1}XY^{-1}W
\end{align*}
and
\[
\bigl(Y+XY^{-1}X\bigr){}^tW=W\bigl(Y+XY^{-1}X\bigr).
\]
In order to show \eqref{K6}, it is sufficient to check
\begin{align}
0={}&Y^{-1}W+Y^{-1}X\bigl(Y+XY^{-1}X\bigr)^{-1}WXY^{-1}-Y^{-1}X\bigl(Y+XY^{-1}X\bigr)^{-1}XY^{-1}W\nonumber\\
&+\bigl(Y+XY^{-1}X\bigr)^{-1}WXY^{-1}XY^{-1}\nonumber\\
&-\bigl(Y+XY^{-1}X\bigr)^{-1}XY^{-1}WXY^{-1}-{}^tWY^{-1}.\label{K6'}
\end{align}
By using above equalities, the right-hand side of \eqref{K6'} can be computed as
\begin{gather*}
Y^{-1}W-Y^{-1}X\bigl(Y+XY^{-1}X\bigr)^{-1}XY^{-1}W+\bigl(Y+XY^{-1}X\bigr)^{-1}WXY^{-1}XY^{-1}-{}^tWY^{-1}\\
\qquad=Y^{-1}W-\bigl(Y+XY^{-1}X\bigr)^{-1}XY^{-1}XY^{-1}W+\bigl(Y+XY^{-1}X\bigr)^{-1}WXY^{-1}XY^{-1}\\
\phantom{\qquad=}{}-\bigl(Y+XY^{-1}X\bigr)^{-1}W\bigl(Y+XY^{-1}X\bigr)Y^{-1}\\
\qquad=Y^{-1}W-\bigl(Y+XY^{-1}X\bigr)^{-1}XY^{-1}XY^{-1}W+\bigl(Y+XY^{-1}X\bigr)^{-1}WXY^{-1}XY^{-1}\\
\phantom{\qquad=}{}-\bigl(Y+XY^{-1}X\bigr)^{-1}W-\bigl(Y+XY^{-1}X\bigr)^{-1}WXY^{-1}XY^{-1}\\
\qquad=Y^{-1}W-\big\{\bigl(Y+XY^{-1}X\bigr)^{-1}XY^{-1}X +\bigl(Y+XY^{-1}X\bigr)^{-1}Y \big\}Y^{-1}W=0.
\end{gather*}
This proves \eqref{K6}.

We show \eqref{K1}. To see this, we show
\begin{align}
0={}& \bigl(Y+XY^{-1}X\bigr)^{-1}XY^{-1}
\begin{pmatrix}
\partial_{x_1}\bigl(XY^{-1}\bigr)_{1i} & \cdots & \partial_{x_n}\bigl(XY^{-1}\bigr)_{1i} \\
\vdots & & \vdots \\
\partial_{x_1}\bigl(XY^{-1}\bigr)_{ni} & \cdots & \partial_{x_n}\bigl(XY^{-1}\bigr)_{ni}
\end{pmatrix} \nonumber \\
&-
\begin{pmatrix}
\partial_{x_1}\bigl(Y^{-1}\bigr)_{1i}-\partial_{x_1}\bigl(Y^{-1}\bigr)_{1i} & \cdots & \partial_{x_n}\bigl(Y^{-1}\bigr)_{1i}-\partial_{x_1}\bigl(Y^{-1}\bigr)_{ni} \\
\vdots & & \vdots \\
\partial_{x_1}\bigl(Y^{-1}\bigr)_{ni}-\partial_{x_n}\bigl(Y^{-1}\bigr)_{1i} & \cdots & \partial_{x_n}\bigl(Y^{-1}\bigr)_{ni}-\partial_{x_n}\bigl(Y^{-1}\bigr)_{ni}
\end{pmatrix}\nonumber \\
&-\bigl(Y+XY^{-1}X\bigr)^{-1}
\begin{pmatrix}
\partial_{x_1}\bigl(XY^{-1}\bigr)_{1i} & \cdots & \partial_{x_n}\bigl(XY^{-1}\bigr)_{1i} \\
\vdots & & \vdots \\
\partial_{x_1}\bigl(XY^{-1}\bigr)_{ni} & \cdots & \partial_{x_n}\bigl(XY^{-1}\bigr)_{ni}
\end{pmatrix}XY^{-1}.\label{K1'}
\end{align}
The $(j,k)$-component of the right-hand side of \eqref{K1'} is
\begin{gather*}
 \sum_l \bigl(\bigl(Y+XY^{-1}X\bigr)^{-1}XY^{-1}\bigr)_{jl}\partial_{x_k}\bigl(XY^{-1}\bigr)_{li}-\partial_{x_k}Y^{-1}_{ji}+\partial_{x_j}Y^{-1}_{ki}\\
 \qquad-\sum_{l,m}\bigl(Y+XY^{-1}X\bigr)^{-1}_{jm}\partial_{x_l}\bigl(XY^{-1}\bigr)_{mi} \bigl(XY^{-1}\bigr)_{lk} \\
\phantom{\qquad-}{}= \bigl(\bigl(Y+XY^{-1}X\bigr)^{-1}XY^{-1}\partial_{x_k}\bigl(XY^{-1}\bigr)\bigr)_{ji}-\partial_{x_k}Y^{-1}_{ji}+\partial_{x_j}Y^{-1}_{ki}\\
\phantom{\qquad-=}{} -\sum_{l,m}\bigl(Y+XY^{-1}X\bigr)^{-1}_{lm}\partial_{x_j}\bigl(XY^{-1}\bigr)_{mi} \bigl(XY^{-1}\bigr)_{lk} \\
\phantom{\qquad-}{}= \bigl(\bigl(Y+XY^{-1}X\bigr)^{-1}\big\{\partial_{x_k}\bigl(\bigl(Y+XY^{-1}X\bigr)Y^{-1}\bigr)-\partial_{x_k}\bigl(XY^{-1}\bigr)XY^{-1}\big\}\bigr)_{ji}\\
\phantom{\qquad-=}{} -\partial_{x_k}Y^{-1}_{ji}+\partial_{x_j}Y^{-1}_{ki}-\sum_{l,m}\bigl(Y+XY^{-1}X\bigr)^{-1}_{ml}\bigl(XY^{-1}\bigr)_{lk}\partial_{x_j}\bigl(XY^{-1}\bigr)_{mi} \\
\phantom{\qquad-}{}= \bigl(\bigl(Y+XY^{-1}X\bigr)^{-1}\bigl(\partial_{x_k}\bigl(Y+XY^{-1}X\bigr)\bigr)Y^{-1}+\partial_{x_k}Y^{-1}\bigr)_{ji}\\
 \phantom{\qquad-=}{} -\bigl(\bigl(Y+XY^{-1}X\bigr)^{-1}\partial_{x_k}\bigl(XY^{-1}\bigr)XY^{-1}\bigr)_{ji}-\partial_{x_k}Y^{-1}_{ji}+\partial_{x_j}Y^{-1}_{ki}\\
\phantom{\qquad-=}{} -\sum_m\bigl(\bigl(Y+XY^{-1}X\bigr)^{-1}XY^{-1}\bigr)_{mk}\partial_{x_j}\bigl(XY^{-1}\bigr)_{mi} \\
\phantom{\qquad-}{}= \bigl(\bigl(Y+XY^{-1}X\bigr)^{-1}\bigl(\partial_{x_k}\bigl(Y+XY^{-1}X\bigr)\bigr)Y^{-1}\bigr)_{ji}\\
\phantom{\qquad-=}{} -\bigl(\bigl(Y+XY^{-1}X\bigr)^{-1}\partial_{x_k}\bigl(XY^{-1}\bigr)XY^{-1}\bigr)_{ji}+\partial_{x_j}Y^{-1}_{ki}\\
\phantom{\qquad-=}{} -\sum_m\bigl(\bigl(Y+XY^{-1}X\bigr)^{-1}XY^{-1}\bigr)_{km}\partial_{x_j}\bigl(XY^{-1}\bigr)_{mi} \\
\phantom{\qquad-}{}= \bigl(\bigl(Y+XY^{-1}X\bigr)^{-1}\bigl(\partial_{x_k}\bigl(Y+XY^{-1}X\bigr)\bigr)Y^{-1}\bigr)_{ji}\\
\phantom{\qquad-=}{} -\bigl(\bigl(Y+XY^{-1}X\bigr)^{-1}\partial_{x_k}\bigl(XY^{-1}\bigr)XY^{-1}\bigr)_{ji}\\
\phantom{\qquad-=}{} +\bigl(\partial_{x_j}Y^{-1}-\bigl(Y+XY^{-1}X\bigr)^{-1}XY^{-1}\partial_{x_j}\bigl(XY^{-1}\bigr)\bigr)_{ki} \\
\phantom{\qquad-}{}= \bigl(-\bigl(\partial_{x_k}\bigl(Y+XY^{-1}X\bigr)^{-1}\bigr)\bigl(Y+XY^{-1}X\bigr)Y^{-1}\bigr)_{ji}\\
\phantom{\qquad-=}{} -\sum_l\bigl(\bigl(Y+XY^{-1}X\bigr)^{-1}\partial_{x_k}\bigl(XY^{-1}\bigr)\bigr)_{jl}XY^{-1}_{li}\\
\phantom{\qquad-=}{} +\bigl(\partial_{x_j}Y^{-1}-\bigl(Y+XY^{-1}X\bigr)^{-1}XY^{-1}\partial_{x_j}\bigl(XY^{-1}\bigr)\bigr)_{ki}\\
\phantom{\qquad-}{}= -\sum_l\partial_{x_k}\bigl(Y+XY^{-1}X\bigr)^{-1}_{jl}\bigl(\bigl(Y+XY^{-1}X\bigr)Y^{-1}\bigr)_{li}\\
\phantom{\qquad-=}{} -\sum_l\bigl(\bigl(Y+XY^{-1}X\bigr)^{-1}\partial_{x_j}\bigl(XY^{-1}\bigr)\bigr)_{kl}XY^{-1}_{li}\\
\phantom{\qquad-=}{} +\bigl(\partial_{x_j}Y^{-1}-\bigl(Y+XY^{-1}X\bigr)^{-1}XY^{-1}\partial_{x_j}\bigl(XY^{-1}\bigr)\bigr)_{ki} \\
\phantom{\qquad-}{}= -\sum_l\partial_{x_j}\bigl(Y+XY^{-1}X\bigr)^{-1}_{kl}\bigl(\bigl(Y+XY^{-1}X\bigr)Y^{-1}\bigr)_{li}\\
\phantom{\qquad-=}{} -\bigl(\bigl(Y+XY^{-1}X\bigr)^{-1}\partial_{x_j}\bigl(XY^{-1}\bigr)XY^{-1}\bigr)_{ki}\\
\phantom{\qquad-=}{} +\bigl(\partial_{x_j}Y^{-1}-\bigl(Y+XY^{-1}X\bigr)^{-1}XY^{-1}\partial_{x_j}\bigl(XY^{-1}\bigr)\bigr)_{ki} \\
\phantom{\qquad-}{}= \bigl(-\bigl(\partial_{x_j}\bigl(Y+XY^{-1}X\bigr)^{-1}\bigr)\bigl(Y+XY^{-1}X\bigr)Y^{-1}\bigr)_{ki}\\
\phantom{\qquad-=}{} -\bigl(\bigl(Y+XY^{-1}X\bigr)^{-1}\partial_{x_j}\bigl(XY^{-1}\bigr)XY^{-1}\bigr)_{ki}\\
\phantom{\qquad-=}{} +\bigl(\partial_{x_j}Y^{-1}-\bigl(Y+XY^{-1}X\bigr)^{-1}XY^{-1}\partial_{x_j}\bigl(XY^{-1}\bigr)\bigr)_{ki} \\
\phantom{\qquad-}{}= \bigl(\bigl(Y+XY^{-1}X\bigr)^{-1}\bigl(\partial_{x_j}\bigl(Y+XY^{-1}X\bigr)\bigr)Y^{-1}\bigr)_{ki}\\
\phantom{\qquad-=}{} -\bigl(\bigl(Y+XY^{-1}X\bigr)^{-1}\partial_{x_j}\bigl(XY^{-1}\bigr)XY^{-1}\bigr)_{ki}\\
\phantom{\qquad-=}{} +\bigl(\partial_{x_j}Y^{-1}-\bigl(Y+XY^{-1}X\bigr)^{-1}XY^{-1}\partial_{x_j}\bigl(XY^{-1}\bigr)\bigr)_{ki}
\\
\phantom{\qquad-}{}=\bigl( \bigl(Y+XY^{-1}X\bigr)^{-1}\big\{\bigl(\partial_{x_j}\bigl(Y+XY^{-1}X\bigr)\bigr)Y^{-1}+\bigl(Y+XY^{-1}X\bigr)\partial_{x_j}Y^{-1}\big\}\bigr)_{ki}\\
\phantom{\qquad-=}{}-\bigl( \bigl(Y+XY^{-1}X\bigr)^{-1}\big\{\partial_{x_j}\bigl(XY^{-1}\bigr)XY^{-1}+XY^{-1}\partial_{x_j}\bigl(XY^{-1}\bigr)\big\}\bigr)_{ki}\\
\phantom{\qquad-}{}=\bigl(\bigl(Y+XY^{-1}X\bigr)^{-1}\big\{\partial_{x_j}(\bigl(Y+XY^{-1}X\bigr)Y^{-1})-\partial_{x_j}(XY^{-1}XY^{-1})\big\}\bigr)_{ki}=0.
\end{gather*}
This proves \eqref{K1}.

Finally, we show \eqref{K2}. We put
\[
V:=
\begin{pmatrix}
\partial_{x_1}\bigl(XY^{-1}\bigr)_{1i} & \cdots & \partial_{x_n}\bigl(XY^{-1}\bigr)_{1i} \\
\vdots & & \vdots \\
\partial_{x_1}\bigl(XY^{-1}\bigr)_{ni} & \cdots & \partial_{x_n}\bigl(XY^{-1}\bigr)_{ni}
\end{pmatrix}.
\]
By \eqref{K1} and \eqref{K3}, we obtain
\begin{align*}
&\begin{pmatrix}
\partial_{x_1}\bigl(Y^{-1}\bigr)_{1i}-\partial_{x_1}\bigl(Y^{-1}\bigr)_{1i} & \cdots & \partial_{x_n}\bigl(Y^{-1}\bigr)_{1i}-\partial_{x_1}\bigl(Y^{-1}\bigr)_{ni} \\
\vdots & & \vdots \\
\partial_{x_1}\bigl(Y^{-1}\bigr)_{ni}-\partial_{x_n}\bigl(Y^{-1}\bigr)_{1i} & \cdots & \partial_{x_n}\bigl(Y^{-1}\bigr)_{ni}-\partial_{x_n}\bigl(Y^{-1}\bigr)_{ni}
\end{pmatrix}\\
&\qquad=\bigl(Y+XY^{-1}X\bigr)^{-1}XY^{-1}V-\bigl(Y+XY^{-1}X\bigr)^{-1}VXY^{-1}
\end{align*}
and
\[
\bigl(Y+XY^{-1}X\bigr){}^tV=V\bigl(Y+XY^{-1}X\bigr).
\]
In order to show \eqref{K2}, it is sufficient to check
\begin{align}
0={}&Y^{-1}V-Y^{-1}X\bigl(Y+XY^{-1}X\bigr)^{-1}XY^{-1}V+Y^{-1}X\bigl(Y+XY^{-1}X\bigr)^{-1}VXY^{-1}\nonumber\\
&+\bigl(Y+XY^{-1}X\bigr)^{-1}VXY^{-1}XY^{-1}-\bigl(Y+XY^{-1}X\bigr)^{-1}XY^{-1}VXY^{-1}\nonumber\\
&-{}^tVY^{-1}.\label{K2'}
\end{align}
Then, \eqref{K2'} can be checked in the same way as \eqref{K6'}.

\subsection*{Acknowledgments}
The many part of this work was done while the author stayed in McMaster university. The author would like to thank the department of Mathematics and Statistics, McMaster university and especially Megumi Harada for their hospitality. The author would also like to express our sincere gratitude to the referees who carefully read the manuscript and helped him improve it. This work is supported by Grant-in-Aid for Scientific Research (C) 15K04857 and 19K03479.

\pdfbookmark[1]{References}{ref}
\LastPageEnding

\end{document}